\documentclass[12pt]{article}
\usepackage{amssymb}
\usepackage{hyperref}
\usepackage[active]{srcltx}
\input{xy}
\xyoption{all}

\newtheorem{theo}{\sc \thesection.\arabic{abz}. Theorem}

\newtheorem{lemm}{\sc \thesection.\arabic{abz}. Lemma}

\newtheorem{coro}{\sc \thesection.\arabic{abz}. Corollary}

\newtheorem{defi}{\sc \thesection.\arabic{abz}. Definition}

\newtheorem{exa}{\sc \thesection.\arabic{abz}. Example}

\newtheorem{rema}{\sc \thesection.\arabic{abz}. Remark}

\newtheorem{conj}{\sc \thesection.\arabic{abz}. Conjecture}

\newcounter{abz}[section]
\newcounter{equ}[abz]
\newcommand{\abz}{\refstepcounter{abz}}
\newcommand{\equ}{\refstepcounter{equ}}

\def\a{\mathfrak{a}}
\def\ad{\mathrm{ad\,}}

\def\Aut{\mathrm{Aut}}
\def\b{\mathfrak{b}}
\def\c{\mathfrak{c}}

\def\d{\mathfrak{d}}
\def\diag{\mathrm{diag}}
\def\Der{\mathrm{Der}}
\def\e{\mathfrak{e}}
\def\End{\mathrm{End}}

\def\f{\mathfrak{f}}
\def\G{\mathcal{G}}
\def\g{\mathfrak{g}}
\def\gl{\mathfrak{gl}}
\def\h{\mathfrak{h}}

\def\hei{\mathrm{ht}}
\def\Id{\mathrm{Id}}
\def\im{\mathop{\mathrm{im}}}

\def\k{\mathfrak{k}}

\def\m{\mathfrak{m}}
\def\mod{\mathrm{mod}}

\def\p{\mathfrak{p}}
\def\pr{\mathrm{pr}}

\def\sl{\mathfrak{sl}}
\def\so{\mathfrak{so}}
\def\su{\mathfrak{su}}
\def\s{\mathrm{Symm}}
\def\ss{\mathfrak{s}}

\def\Span{\mathop{\mathrm{Span}}}

\def\sp{\mathop{\mathfrak{sp}}}
\def\Tr{\mathop{\mathrm{Tr}}}
\def\Tor{\mathop{\mathrm{Tor}}}
\def\uu{\mathfrak{u}}

\def\z{\mathfrak{z}}
\def\Z{\mathbb{Z}}

\def\R{\mathbb{R}}

\def\L{\mathcal{L}}
 \def\C{\mathbb{C}}

\def\T{\mathcal{T}}

\def\al{\alpha}
\def\be{\beta}
\def\Ga{\Gamma}
\def\ga{\gamma}

\def\ep{\varepsilon}
\def\ka{\kappa}
\def\la{\lambda}

\def\si{\sigma}

\def\qed{$\square$}
\renewcommand\phi{\varphi}

\newcommand{\lbr}{\linebreak[0]}
\date{}
\title{Lie--Poisson pencils related to semisimple Lie algebras: towards classification}
\author{Andriy Panasyuk\\
 \\
Faculty of Mathematics and Computer Science\\
University of Warmia and Mazury\\
ul. S\l oneczna 54, 10-710 Olsztyn,
Poland\\
and\\
Pidstryhach Institute for Applied Problems\\ of Mechanics and Mathematics, NASU\\
3-b, Naukova Str., 79060 L'viv, Ukraine }
\begin{document}
\maketitle
\tableofcontents
\vspace{3cm}
\hspace{9cm}
\begin{minipage}[b]{8cm}
{\em
To the memory of my parents\\

}
\end{minipage}

\vspace{1cm}
\section{Introduction}
\label{intro}

Since its invention by Magri  \cite{m} the  bihamiltonian property proved to be very useful in the study of integrable systems. The second hamiltonian structure was discovered for majority of important examples. In most of the cases it has been first found "by hand" but afterwards it was seen as a subcase of some general algebraic scheme. For instance, the bihamiltonian structure of the KdV \cite{m}  later was understood as the so-called argument shift method \cite[Ch. VI]{arnKhesin}, \cite{km} on the Virasoro algebra.

Recall that a bihamiltonian structure is a pair of Poisson structures on a manifold which are compatible, i.e. such that their sum is again a Poisson structure. There is a lot of the above mentioned algebraic schemes of appearing bihamiltonian structures.  However, using a simple characteristic, the order of the coefficients of the corresponding Poisson bivectors, one can organize these schemes in the following hierarchy:
\begin{itemize}
\item constant+constant;
\item linear+constant;
\item linear+linear;
\item quadratic+linear;
\item etc.
\end{itemize}
The first case (when the coefficient of both the brackets are constant) is rather not interesting. The second one (here one of the bivectors is the Lie--Poisson structure on a Lie algebra and another one is related to a cocycle on this Lie algebra) is already very important, for instance, due to the argument shift  method mentioned, see also \cite{bols'}. The third case is the subject of this paper and we will comment it below. The higher steps in the hierarchy are usually related to different forms of the classical Yang--Baxter equations \cite{mks89}, \cite[\S 1--2]{karasevMaslov}, \cite{doninGurevich}, \cite{gurevichPanyushev}, \cite{gurevichRubtsov} (however it is worth mentioning that the  argument shift of arbitrary quadratic Poisson structure  in the vanishing direction always gives a "quadratic+linear" pair).

The "linear+linear" case is represented by pairs of  compatible Lie--Poisson structures, or, equivalently,
pairs of compatible structures of Lie algebras on a vector space (cf. Definition \ref{0.10}). The situation here is essentially more complicated than in the "linear+constant" case (note that it is easy to construct examples of  pairs related to the last case on any Lie algebra considering a trivial cocycle). There are no "automatic" "linear+linear" examples related to arbitrary Lie algebra (with the exception of the Lie bracket being the antisymmetrization of an associative multiplication $x\cdot y$: for any fixed element $a$ the operation $x\cdot_ay:= x\cdot a \cdot y$ is a new associative product and the corresponding commutator is compatible with the initial one).

One has to distinguish two main fields on which "linear+linear" pairs appeared.
 In finite-dimensio-\linebreak nal context and within the "purely" bihamiltonian theory the pair $[,],[,_A]$ of compatible Lie brackets on $\so(n)$  (here $[,]$ is the standard commutator and $[,_A]$ is the modified one, $xAy-yAx$, with  a diagonal matrix $A$) was applied by Bolsinov \cite{bols'}, \cite{bolsLiePencils} to prove the integrability  of the Euler--Manakov top and to show its relation to the so-called Clebsch--Perolomov  case (see also \cite{morosip}, \cite{yanovski}). Another Lie--Poisson pencil on $\so(n)\times\so(n)$ was used for study of generalized Steklov--Lyapunov systems  \cite{bolsFedorov}. A Poisson pencil which incorporate, on one hand, a Lie--Poisson pair (either the $\so(n)$ pair mentioned, or another one related to a $\Z_2$-grading on a Lie algebra), and, on the other hand, the argument shift  method, was used in the theory of a series of systems (in particular for building  Lax representations) such as Brun--Bogoyavlenskii and Zhukovsky--Volterra systems, the Kovalevski top, etc. \cite{bolsBorisov}, \cite[Ch. 2]{borisovMamaev}.

Another area where  pairs of compatible Lie brackets play important role is the classical $r$-matrix formalism. The main idea (which in implicit form is contained already in the paper of Holod \cite{holod}) that explains this role is as follows. We say that a Lie algebra $(\hat\g,[,])$ with a decomposition $\hat\g=\oplus_{n\in\Z}\g_n$ is {\em quasigraded}  of degree $1$ if $[\g_i,\g_j]\subset\g_{i+j}\oplus \g_{i+j+1}$ for any $i,j\in\Z$. Given such a Lie algebra,  one observes that $\g_+:=\oplus_{n\ge 0}\g_n$ and $\g_-:=\oplus_{n< 0}\g_n$ are subalgebras, so the kernel of the difference of the corresponding projectors $P_+-P_-$ is the classical $r$-matrix and one can apply the standard Adler--Costant scheme \cite{rst}. On the other hand, if $\g$ is a vector space with two compatible Lie brackets $[,]_0,[,]_1$, one obtains a  quasigraded Lie algebra of degree $1$ by putting $\hat\g:=\g[\la,1/\la],[,]:=[,]_0+\la[,]_1$ and extending this bracket by bilinearity to $\hat\g$ (the grading is usual, by the degree of a Laurent polynomial).
 This idea was consequently developed by Skrypnyk \cite{skrypnyk2002}, \cite{skrypnyk2004}, \cite{skrypnyk2006}, and, in  a slightly different form, by  Golubchik and Sokolov \cite{golsok}, \cite{golSok2005}  generating many new results on finite- and infinite-dimensional integrable systems (see also  \cite{lombardoMikhailov}, \cite[Ch. 16]{GeViYanovsky}). Note that the correspondence \{compatible Lie brackets on $\g$\} $\to$ \{Lie algebra structures $(\g[\la,1/\la],[,])$ such that $\g_+,\g_-$ are subalgebras\} is  one-to-one under some reasonable restrictions on the bracket $[,]$.

The above mentioned results are enough motivating for considering the classification matter of pairs of compatible Lie brackets on finite dimensional spaces. However, this matter (like many other questions of Lie theory) is itself a beautiful mathematical problem which probably inspired Kantor and Persits in their pioneering work in this direction. In \cite{KantorPersits} they announced (unfortunately, the proof never appeared) the following result.  Let $\g,V$ be finite-dimensional vector spaces and let $\{[,]^v\}_{v\in V}$ be a family of Lie brackets on $\g$ parametrized by $V$ (note that any pair of compatible  Lie brackets generates such a family with $V$ $1$- or $2$-dimensional). We say that such a family is {\em irreducible} if the Lie algebras $(\g,[,]^v)$ do not have common ideals and is {\em closed} if for any $a\in\g,v,w\in V$ there exists $z\in V$ such that $[,]^v_{\ad^wa}=[,]^z$, where $\ad^wa(b):=[a,b]^w$ and we used   notation (\ref{e10.10}). The result of Kantor and Persits gives a complete list of closed irreducible $\C$-vector spaces of Lie brackets of dimension greater than one:
\begin{enumerate}\item $\g=\so(n,\C)$ is the space of $n\times n$ antisymmetric matrices, $V$ is the space of $n\times n$ symmetric matrices, the commutator $[,]^X,X\in V$, is given by $[a,b]^X=aXb-bXa=:[a,_Xb]$;
\item $\g$ is the space of $n\times n$ symmetric matrices, $V=\so(n,\C),[,]^X=[,_X]$;
\item $\g$ is the space of $n\times m$  matrices, $V$ is the space of $m\times n$  matrices, $[,]^X=[,_X]$;
\item $\g=V$ is an even-dimensional  vector space with a nondegenerate skew-symmetric form $(,)$, the commmutator is given by $[a,b]^v:=(v,a)b-(v,b)a-(a,b)v$;
\end{enumerate}
Compatible Lie brackets from this list were studied in \cite[\S 44]{TFbols} and \cite{yanovski}.

Another result which should be mentioned here is due to Odesski and Sokolov \cite{odesskisok} who classified associative multiplications on the complex $n\times n$ matrices compatible with the usual one (in particular, by antisymmetrization one can deduce from this result a lot of new examples of compatible Lie brackets on $\gl(n,\C)$, see also \cite{cgm3}).

The main objective of this paper is to outline a systematic approach to classification of pairs of compatible complex Lie brackets (we call such pairs {\em bi-Lie structures} in this paper)  one of which is semisimple.  As a byproduct we obtain new examples of such pairs which are not contained in the following
list of known examples (KE) of this type (in Items 1,2 $[,]$ stands for the standard commutator of matrices).
\begin{enumerate}\item[KE1] $(\so(n,\C),[,],[,_A])$, where $A$ is a symmetric $n\times n$-matrix (cf. the first example from the Kantor--Persits list and Example \ref{0.30}).
\item[KE2] $(\sp(n,\C),[,],[,_A])$, see Example \ref{0.35} (this example can be derived from the second item of the Kantor--Persits list in the case when the matrix $X$ is  the standard antisymmetric matrix of maximal rank equal to $2n$ \cite[\S 44]{TFbols}).
\item[KE3] Let $(\g,[,])$ be semisimple. There exists a  bi-Lie structure related  to any $\Z_m  $-grading on $(\g,[,])$ and to decomposition of the subalgebra $\g_0$ to two subalgebras, see Examples \ref{10.40}, \ref{140.95}, Theorem \ref{130.70} (in general case first appeared in \cite{golsok}, for $m=2$ known much earlier, see discussion above).
\item[KE4]  Let $(\g,[,])$ be semisimple. There exists a  bi-Lie structure related  to any parabolic subalgebra $\g_0 \subset\g$, see Example \ref{10.35} (in case when $\g_0$ is a Borel subalgebra it appeared in \cite{p5}).
\item[KE5] Examples on $\sl(3,\C),\so(4,\C)$ related to $\Z_2 \times\Z_2 $-gradings \cite{golsok}.
\end{enumerate}
Another byproduct is establishing of isomorphisms between some of bi-Lie structures from Items KE3 and KE4, which are not obvious.

It is known \cite{mks},\cite{golsok} that, given a semisimple Lie algebra $(\g,[,])$ and a bilinear bracket $[,]'$, the triple $(\g,[,],[,]')$ is a bi-Lie structure if and only if $[,]'$ is of the form $[,]'=[,]_W$ (see \ref{e10.10}) for some linear operator $W\in\End(\g)$ such that there is another operator $P\in\End(\g)$  satisfying  the "main identity" $T_W(,)=[,]_P$, where $T_W$ is the {\em Nijenhuis torsion} of $N$ (see Section \ref{prelim}).  We call an operator $W$ satisfying $T_W(,)=[,]_P$ for some $P$ a {\em weak Nijenhuis operator} (WNO for short) and $P$ itself a {\em primitive}  of $W$. The problem is that these operators are defined by the bi-Lie structure nonuniquely, up to adding a differentiation of the bracket $[,]$ (which is inner). Our first main result (see Section \ref{40}) proposes  a special fixing of these two operators and establishes some important properties of them. Namely, among all the WNOs corresponding to a bi-Lie structure $(\g,[,],[,]')$ we choose that (uniquely defined) orthogonal to $\ad\g$ in $\End(\g)$ with respect to the "trace form" (we call it {\em principal}). We further prove that two bi-Lie structures are isomorphic (see Definition \ref{0.15}) if and only if so are the corresponding principal WNOs    (Theorem \ref{40.nonum}). We also show that the principal operator $W$ of a  bi-Lie structure $(\g,[,],[,]')$ possesses the unique primitive (called {\em principal}) which is symmetric with respect to the Killing form on $(\g,[,])$ (Theorem \ref{4}). Thus the problem of classification of bi-Lie structures with semisimple $(\g,[,])$ is reduced to the problem of classification  of the pairs $(W,P)$, where $W$ is a principal WNO on $(\g,[,])$ and $P$ is its principal primitive, with respect to the natural action of the automorphisms of $(\g,[,])$.

The last problem in general seems to be very complicated.
We are trying to solve it under the following reasonable restriction: a principal WNO preserves some  $\Ga$-grading on $(\g,[,])$, where $\Ga$ is an abelian group. This restriction is motivated by the fact that the majority of examples from the list above are in fact of this type (for Examples KE3, KE5 this is obvious, one needs a while to see a $\Z_2\times \cdots\times\Z_2$-grading in Examples KE1, KE2 if $A$ is diagonalizable, for Example KE4 see Theorem \ref{130.70}). We first prove (see Theorem \ref{100.20}) that if a principal WNO preserves a grading, then so does its principal primitive. Further on, we make one more restriction and consider principal WNOs $W$ such that:
\begin{enumerate}\item[(1)]
$W$ preserves the root grading $\g=\h+\sum_{\al\in R}\g_\al$ related to some Cartan subalgebra $\h$ in $(\g,[,])$; \item[(2)] $W|_\h$ is semisimple.
\end{enumerate}
 For such operators the "main identity" $T_N(,)=[,]_P$ becomes a system of $1$-dimensional algebraic equations of second order (see Theorem \ref{110.20}), that is why is in principle solvable. Next we propose a method of solving it.

More precisely, this  system consists of two parts. First is a system of quadratic equations indexed by the roots and relating the eigenvalues of the operators $W,P$ with another  invariants, which we call the {\em times} of the corresponding bi-Lie structure $(\g,[,],[,]')$. These are the complex numbers $t_i$ such that the Lie algebra $(\g,[,]'-t_i[,])$ is not semisimple (in the pencil $\{[,]'-t[,]\}_{t\in\C}$ the generic bracket is isomorphic to $[,]$ but for a finite number of values of $t$ the corresponding brackets are nonsemisimple; we call them {\em exceptional}). The second part, indexed by positive roots (with respect to any basis of roots) relates the eigenvectors and eigenvalues of $W|_\h$ with the times. As a result we obtain a family of unordered pairs $\{t_{i,\al}t_{j,\al}\}_{\al\in R}$ of complex numbers indexed by the roots (we call such a family a {\em pairs diagram}) and a system of restrictions on the operator $W|_\h$ (see Section \ref{110}). It turns out that the "main identity" implies some rules of behavior of the pairs $\{t_{i,\al}t_{j,\al}\}$ with respect to the addition of the roots (see Theorem \ref{110.30}, Item 1). This relates bi-Lie structures with the geometry of the corresponding root systems and makes it possible to classify these first.

Our next result (Theorem \ref{120.10}) says that the set of  pairs diagrams splits in two disjoint classes, I and II. Further on, we study bi-Lie structures corresponding to these classes.

The rules of behavior of unordered pairs $\{t_{i,\al}t_{j,\al}\}$ in the pairs diagrams of Class I with respect to the addition of roots remind very much the operation rules of the famous pair groupoid and induce (together with the operator $W|_\h$ subject to the restrictions mentioned) on $(\g,[,])$ a special type of grading which we call a {\em toral symmetric pairoid quasigrading of Class I} (see Section \ref{140}). It turns out that this quasigrading together with the grading set $\{t_1,\ldots,t_n\}$ contains all the information about the corresponding bi-Lie structure (see Theorem \ref{140.50}; note that the corresponding WNO has the symmetric restriction to $\h^\bot=\sum_{\al\in R}\g_\al$). Thus the classification problem of bi-Lie structures of Class I is equivalent to the classification problem of toral symmetric pairoid qusigradings of Class I. The last remains open, however we present a list of examples (which contains known and new ones) and conjecture that this list is in a sense complete (see Conjecture \ref{140.170}). We remark that any such quasigrading is related to a specific type of $\Z_2 \times \cdots\times\Z_2 $-grading on $(\g,[,])$ (see Lemma \ref{140.70}).

In the case of pairs diagrams of Class II the corresponding bi-Lie structures also induce on $(\g,[,])$ a special type of gradings which we call a {\em toral symmetric pairoid quasigradings of Class II}. They are much simpler than that  of Class I, however, they  contain only a part of information about the initial bi-Lie structures. The crucial difference here, in contrast to Class I structures, lies in the fact that the corresponding WNOs have also antisymmetric parts on $\h^\bot$ which contain the rest of the information about the bi-Lie structures (in Theorem \ref{130.50} we show how to rebuild the bi-Lie structure from the corresponding quasigrading and the antisymmetric part of the WNO on $\h^\bot$). These antisymmetric parts are subject to some restrictions implied by the "main identity" (see Theorem \ref{110.30}, Item 2) and in principle can be classified. The problem of classification of  bi-Lie structures of Class II is now reduced to the problem of classification of toral symmetric pairoid qusigradings of Class II and the corresponding antisymmetric operators. Although this last is not solved in this paper in full extent, we present a list of examples (containing known and new ones) and conjecture the completeness of this list (see Conjecture \ref{130.90}). We also prove this conjecture for  $\g=\sl(n,\C)$ thus obtaining a complete classification of bi-Lie structures of Class II in this case (Theorem \ref{130.100}). Note that bi-Lie structures corresponding to pairs diagrams of Class II are also related to some gradings on $(\g,[,])$ which are coarsenings of the root grading (see Theorem \ref{130.50}).

Let us exhibit the examples which we get within our theory and compare them with that from the list above.
\begin{enumerate}\item Class I
\begin{enumerate}\item We recover Example KE1 in the particular case of $A=\diag(t_1,t_1,\ldots,t_n,t_n,t_{n+1})$ for $\so(2n+1,\C)$ and $A=\diag(t_1,t_1,\ldots,t_n,t_n)$ for $\so(2n,\C)$ (see Examples \ref{140.90}, \ref{140.100}).
\item We recover Example KE2 in the particular case of $A=\diag(t_1,\ldots,t_n,t_1,\ldots,t_n)$ (see Example \ref{140.110}).
\item We obtain new examples of bi-Lie structures, which are analogues of the examples above, for $\sl(n,\C)$ and which, surprisingly, were not known in the literature (see Examples \ref{140.120}, \ref{140.140}).
\item We recover Example KE3 for $m=2$ (see Example \ref{140.95}). \end{enumerate}
\item Class II
\begin{enumerate}\item We recover Example KE3 for $m>2$ in the case when the $\Z_m$-grading is related to an inner automorphism (see Theorem \ref{130.70}).
\item We recover and generalize Example KE4  (see Theorem \ref{130.60}).
\item We obtain a new example of a bi-Lie structure on the exceptional Lie algebra $\e_7$ related to $\Z_3 \times\Z_3 $-grading (see Example \ref{130.80}).
\end{enumerate}
\end{enumerate}

Moreover, we prove that all the bi-Lie structures from Example KE4 are isomorphic to special bi-Lie structures from Example KE3 (see Theorem \ref{130.70}). We  also deduce from Examples \ref{140.120}, \ref{140.140} new examples of real bi-Lie structures on $\su(n)$ and, in particular, on $\so(6,\R)\cong\su(4)$ (see Appendix II).

Remark that a specific form of the matrix $A$ which appeared in Example Class I (a) is related to the fact that the above mentioned restrictions (1), (2) on the WNO imply the following condition on the corresponding bi-Lie structure: the sum of the centres of the exceptional brackets contains the Cartan subalgebra $\h$ (see Theorem \ref{110.10}). In the case of generic $A$ (diagonalizable with simple spectrum) in Example KE1 these centres are trivial and this case does not fit our theory. However, this case is related to a specific pairoid quasigrading which is not toral (see Remark \ref{140.200}).

Summarizing, we can say that although the theory of the second part of the paper (Sections \ref{110}--\ref{130}) recover only examples related to toral gradings and quasigradings (cf.   Appendix I and Definition \ref{140.20}), or in other words, to inner automorphisms, the results of the first part give a hope that similar theory can be also built in general, i.e. including gradings coming from outer automorphisms too. The author also hopes that some of these results (for instance Theorems \ref{40.th}, \ref{4}, \ref{40.nonum}, \ref{t10.1}, \ref{overcentral}, \ref{100.20}), which deal with general properties of bi-Lie structures, are new and useful.

Now let us overview the content of the paper in details. In Section \ref{basic} we give basic definitions (of bi-Lie structures and their isomorphisms) and examples. We also prove a lemma, which will be used in the proof of Theorem \ref{40.th}. In Section \ref{cohom} we recall basic formulae from the theory of cohomology of Lie algebras, mainly for fixing notations.

Section \ref{prelim} is devoted to WNOs and their primitives. We give definitions, basic facts and examples. We also prove a result (Theorem \ref{10.34}) which shows how the primitive changes when we add a differentiation to the WNO.

In Section \ref{15} we introduce a notion of a {\em parial Nijenhuis operator} (PNO for short) which in fact provides (besides the WNO) an alternative way of description of bi-Lie structures (in the semisimple case). The notion of a PNO is auxiliary but convenient for understanding the relation between bi-Lie structures and WNOs. The relation between PNOs and WNOs is considered in Section \ref{30}.

This relation is exhibited in the full extent in Section \ref{40}, where we prove the existence of a WNO related to any bi-Lie structure $(\g,[,],[,]')$ such that $(\g,[,])$ is semisimple, see Theorem \ref{40.th} (note that this proof is independent of the Whitehead lemmata on the triviality of the cohomology of semisimple Lie algebras). We also introduce thea notion of a principal WNO and prove the existence of a symmetric primitive for the principal WNO (again without a reference to the Whitehead lemmata, see Theorem \ref{4}). Finally, we show how the principal WNO determines the centre of the bracket $[,]'=[,]_N$ (Corollary \ref{coro2}) and prove the equivalence of the category of bi-Lie structures $(\g,[,],[,]')$ with the same semisimple Lie algebra $(\g,[,])$ and the category of principal WNOs on $(\g,[,])$ (Theorem \ref{40.nonum}).

In Section \ref{90} we study the pencil of brackets  $[,]'-t[,],t\in\C$, in particular, the kernels of the Killing forms of the exceptional Lie brackets. We introduce two important subalgebras in $(\g,[,])$, a {\em central} subalgebra $\z \subset \g$, which is the sum of the centres of the exceptional brackets, and  a subalgebra $\hat{\z}\supset\z$, which roughly speaking is the maximal subspace in $\g$ such that $W|_{\hat\z}$ is a {\em Nijenhuis} operator, i.e. such that its Nijenhuis torsion vanishes.

Section \ref{100} is devoted to WNOs preserving a group grading on $(\g,[,])$. The main result (Item 1 of Corollary \ref{100.25}) says that if a principal WNO preserves a grading, then so does its principal primitive. We also present Theorem \ref{100.260} which, given a WNO $W$ preserving the root grading with respect to a Cartan subalgebra, describes the principal WNO $\pr(W)$ with the property $[,]_W=[,]_{\pr(W)}$. Although we apply  this result only partially (in the proof of Theorem \ref{130.50}), it is useful in calculations and we retain it in full extent.

Section \ref{110} is central. Here we study the bi-Lie structures related to the WNOs satisfying conditions (1), (2). In Theorem \ref{110.20} we obtain the system of equations on the operators $W,P$ implied by the "main identity" $T_W(,)=[,]_P$, describe the central subalgebra and the set of times. We introduce the notions of {\em admissible operator} on a Cartan subalgebra implied by the preceding theorem (Definition \ref{110.255}). Further on, we introduce a notion of a pairs diagram (see Definition \ref{110.50}), and show that conditions (1), (2) imply  existence of a pairs diagram (Item 1 of Theorem \ref{110.30}). In this theorem we also obtain restrictions on the antisymmetric part of $W|_{\h^\bot}$ (Item 2) and describe an important subalgebra of $(\g,[,])$ called the {\em basic subalgebra} (Items 3--6) of a bi-Lie structure (the last is related to the above mentioned subalgebra $\hat\z$, see Remark \ref{remaover}). Finally we formalize the restrictions on the antisymmetric part of $W|_{\h^\bot}$ in a notion of an {\em operator obeying the triangle rule} (Definitions \ref{130.30}, \ref{110.60}).

In Section \ref{s120} we prove that the set of pairs diagrams splits into two disjoint parts: Class I and Class II.

In Section \ref{140} we study bi-Lie structures related to pairs diagrams of Class I. We introduce a notion of toral symmetric pairoid quasigrading of Class I on $(\g,[,])$ (Definition \ref{140.20}) and show that it is equivalent to the pair consisting of a pairs diagram and the corresponding admissible operator (Lemma \ref{140.40}). Theorem \ref{140.50} shows how to construct a WNO by a quasigrading and asserts that all WNOs related with bi-Lie structures of Class I are of this kind. It also gives necessary and sufficient conditions for two such  bi-Lie structures to be isomorphic. Then we show that the structure of a toral symmetric pairoid quasigrading of Class I implies a specific type $\Z_2 \times \cdots \times\Z_2 $-grading on $(\g,[,])$ called {\em admissible} (Lemma \ref{140.70}). We further present examples of gradings which are not admissible (Example \ref{140.90}), of admissible grading not related with  toral symmetric pairoid quasigradings of Class I (Example \ref{140.80}), and a series of examples of these last (Examples \ref{140.95}--\ref{140.150}). In Lemma \ref{140.160} we propose a method of constructing a toral symmetric pairoid quasigrading of Class I from a given one. Finally, we conjecture that all bi-Lie structures of Class I are obtained from the listed examples or as reductions from them with the help of Lemma \ref{140.160} (see Conjecture \ref{140.170}).

Section \ref{130} is devoted to bi-Lie structures related to pairs diagrams of Class II. Similarly to the previous case we start from introducing of the notion of a toral symmetric pairoid quasigrading of Class II and showing its equivalence to the corresponding pair consisting of an admissible operator and a pairs diagram (Lemma \ref{130.20}). In Theorem \ref{130.50} we show how to construct a WNO $W$ by a quasigrading (which in fact is determined by two reductive subalgebras, see Remark \ref{130.255}) and an antisymmetric operator on $\h^\bot$ obeying the triangle  rule. Note that from the quasigrading itself we get $W|_\h$ and the symmetric part of $W|_{\h^\bot}$. Again, all bi-Lie structures of Class II can be obtained by this construction (Item 2 of Theorem \ref{130.50}). Then we deal with the question of uniqueness (the answer is given by Theorem \ref{130.500}). In Theorem \ref{130.60} we construct a series of examples generalizing Example  KE4. In Theorem \ref{130.70} we obtain a series of examples recovering Example KE3 and prove the existence of some isomorphisms between them and with examples from the previous series. Then we present an example of bi-Lie structure of Class II on the Lie algebra $\e_7$ using a construction combining that from Theorem  \ref{130.50} and an observation related to the "triangle rule" from preceding examples. We conclude the section by conjecturing that this construction is universal, i.e. all the bi-Lie structures of Class II can be obtained by means of it (Conjecture \ref{130.90}). We also prove this conjecture for  $\g=\a_n$, in particular obtaining the classification of bi-Lie structures of Class II for this Lie algebra.

The paper contains two appendices. In Appendix I we discuss toral gradings of semisimple Lie algebras $\g$ corresponding to regular reductive subalgebras $\g_0$. In particular, we recall \cite{ostapenko} that the family of irreducible components of the natural representation of $\g_0$ in $\g$ forms a group grading, which we call the {\em toral irreducible grading related to $\g_0$}. We also recall \cite{dynkin}, \cite{doninOst} the classification of regular reductive subalgebras.

Appendix II is devoted to examples of real bi-Lie structures related to the complex ones of Class I. We show that the WNO constructed in Theorem \ref{140.50} under some additional conditions can be restricted to the compact real form of the complex semisimple Lie algebra $(\g,[,])$, see Theorem \ref{ApReal.10}. We also deduce, by means of this theorem, examples of real bi-Lie structures  from complex ones, in particular for the Lie algebra $\su(n)$.

\section{Preliminaries and basic examples}
\label{basic}

We assume throughout the paper that the ground field is $\C$ (with the exception of Appendix II).

\abz\label{0.10}
\begin{defi}\rm A bi-Lie structure is a triple $(\g,[,],[,]')$, where $\g$ is a vector space and $[,],[,]'$ are two Lie algebra structures on $\g$ which are {\em compatible}, i.e.  so that any their linear combination $[,]^\la:=\la[,]+\la'[,]'$ is a Lie algebra structure.
\end{defi}

Note that, given  Lie brackets $[,],[,]'$  on $\g$, they are compatible if and only if $[,]+[,]'$ is a Lie bracket.

\abz\label{0.20}
\begin{exa}\rm Let $\g=\gl(n,\C)$, $ A\in\g$ be a fixed matrix. Put $[x,y]'=xAy-yAx=[x,_Ay]$. Then it is easy to see that the triple $(\g,[,],[,]')$, where $[,]$ is the standard commutator, is a bi-Lie structure.
\end{exa}

\abz\label{0.30}
\begin{exa}\rm
Let $\g=\so(n,\C)$, $ A\in\s(n,\C)$, a fixed symmetric matrix. Then again $(\g,[,],[,]')$, where $[,]'$ is given by the same formula as above,  is a bi-Lie structure. Note that symmetric matrices form the orthogonal complement to $\so(n,\C)$ in $\gl(n,\C)$ with respect to the trace form $\Tr(XY)$.
\end{exa}

\abz\label{0.35}
\begin{exa}\rm
Let $\g=\sp(n,\C)=\{X\in\gl(2n)\mid XJ+JX^T=0\}$, here
$J=\left[
\begin{array}{cc}
0 & I \\
 -I & 0 \\
 \end{array}
 \right]$ is the matrix of the standard symplectic form, and let $A$ belong to the orthogonal complement $\sp(n,\C)^\bot$ to the symplectic subalgebra in $\gl(2n)$ with respect to the trace form. This subspace can be described by the formula $\{X\in\gl(2n)\mid XJ-JX^T=0\}$.
  Then again $(\g,[,],[,_A])$, is a bi-Lie structure.
\end{exa}

The following lemma will be used later.

\abz\label{0.40}
\begin{lemm}
Let $(\g,[,])$ and $(\g,[,]')$ be Lie algebras and let $\ad,\ad'$ be the corresponding $\ad$-operators, $\ad x(y)=[x,y],\ad' x(y)=[x,y]'$. The triple $(\g,[,],[,]')$ is a bi-Lie structure if and only if
$$
\ad [x,y]'+\ad'[x,y]=[\ad'x,\ad y]+[\ad x,\ad' y],x,y\in\g.
$$
\end{lemm}

\noindent Indeed, $[,],[,]'$ are compatible if and only if $[,]'':=[,]+[,]'$ is a Lie bracket if and only if $\ad''[x,y]''=[\ad'' x,\ad'' y]$, i.e. $(\ad +\ad')([x,y]+[x,y]')=[\ad x+\ad' x,\ad y+\ad'y]$. Simplifying the last condition with the account of $\ad[x,y]=[\ad x,\ad y]$ and $\ad'[x,y]'=[\ad' x,\ad'y]$ gives the result. \qed

\abz\label{0.15}
\begin{defi}\rm We say that bi-Lie structures $(\g,[,],[,]')$ and $(\g_1,[,]_1,[,]_1')$ are {\em strongly isomorphic (isomorphic)} if there exists a linear map $\phi:\g\to\g_1$ transforming isomorphically the bracket $[,]$ to $[,]_1$ and $[,]'$ to $[,]_1'$ (transforming isomorphically the brackets $[,]$ to $[,]_1$ to some linear combinations of the brackets $[,]'$ to $[,]_1'$).
\end{defi}

\abz\label{0.17}
\begin{rema}\rm The notion of isomorphism of bi-Lie structures is  more natural in the context of studying the  {\em pencils } $\{[,]^\la\}$ and $\{[,]_1^\la\}$ of Lie brackets generated by the pairs $([,],[,]')$ and $([,]_1,[,]_1')$.
\end{rema}

\section{Lie algebra cohomology}
\label{cohom}

In this section we recall basic definitions from the theory of (low-dimensiomnal) cohomology of a Lie algebra $(\g,[,])$ with coefficients in a $\g$-module $\m$. The space of {\em $j$-cochains} $C^j(\g,\m)$ consists of $j$-linear skew-symmetric maps $\g^j\to\m$. By definition $C^0(\g,\m)=\m$.

The coboundary operator $\d^j:C^j(\g,\m)\to C^{j+1}(\g,\m)$ is defined by
$$
\begin{array}{l}
\d \al(h_1):=h_1\al\ \mbox{for}\ j=0;\\ \d\al(h_1,h_2):=h_1\al(h_2)-h_2\al(h_1)-\al([h_1,h_2])\ \mbox{for}\ j=1;\\ \d\al(h_1,h_2,h_3):=\sum_{\mathrm{c.p.}i,j,k}(h_i\al(h_j,h_k)+\al(h_i,[h_j,h_k]))\ \mbox{for}\ j=2
\end{array}
$$
(here $\al\in C^j(\g,\m), h_i\in\g$). The cohomology space is defined standardly: $H^j(\g,\m):=\ker\d^j/\im\d^{j-1}$. The Whitehead lemmata assert  that $H^j(\g,\m)=0$ if $\g$ is semisimple and $j=1,2$.

For the adjoint module $\m=\g$ we will also use notations $L(\g):=C^1(\h,\m)$ (later we will also write $\End(\g)$ for $L(\g)$), $L^2(\g):=C^2(\g,\m)$, and
\abz\equ\label{e10.10}
\begin{equation}
[\cdot,\cdot]_{W}:=\d W(\cdot,\cdot)=[W \cdot,\cdot]+[\cdot,W \cdot]-W[\cdot,\cdot], W\in L(\g).
 \end{equation}

\section{Weak  Nijenhuis operators}
\label{prelim}

Given  two Lie algebras $(\g,[,])$ and $(\g,[,]')$, we observe that $(\g,[,],[,]')$ is a bi-Lie structure if and only if $[,]'$ is a 2-cocycle with respect to $[,]$ with coefficients in the adjoint module. In particular, if $(\g,[,],[,]')$ is a bi-Lie structure such that $(\g,[,])$ is semisimple, then $[,]'=[,]_W$ for some $W\in L(\g)$ (note that in Theorem \ref{40.th} we will prove existence of such $W$ without using the Whitehead lemmata).

Vice versa, let an operator $W\in L(\g)$ be given such that $[,]_W$ is a Lie bracket. Then, since $[,]_W$ is a cocycle, it is automatically compatible with $[,]$ and $(\g,[,],[,]_W)$ is a bi-Lie structure.

We come to the following definitions.

\abz\label{10.10}
\begin{defi}\rm Let $(\g,[,])$ be a Lie algebra.
An operator $W\in L(\g)$ is called a {\em weak Nijenhuis operator} (WNO for short) if the bracket
$[\cdot,\cdot]_{W}$ defined by formula (\ref{e10.10})
 is a Lie bracket on $\g$ (the last will be called the {\em modified bracket} on $\g$).
\end{defi}

\abz\label{10.15}
\begin{rema}\rm A similar structure in the context of Lie algebroids appeared in \cite{cgm2}, we borrowed the terminology from this paper.
\end{rema}

\abz\label{10.20}
\begin{defi}\rm
The Nijenhuis torsion  of an operator
$W\in L(\g)$ is defined as
$$
T_{W}(x,y):=[Wx,Wy]-W[x,y]_{W}.
$$
\end{defi}

The following lemma justifies the term "WNO" (in the sense that it generalizes the well known notion of a Nijenhuis operator, i.e. an operator $W$ with vanishing $T_{W}$).

\abz\label{0}
\begin{lemm} \cite{mks}
An operator $W\in L(\g)$ is a WNO if and only if $T_{W}$ is a 2-cocycle on $\g$ with coefficients in the adjoint module.
\end{lemm}

\abz\label{10.25}
\begin{defi}\rm Let $W\in L(\g)$ be a WNO. An operator $P\in L(\g)$ is said to be a {\em primitive } of $W$ if $T_W(x,y)=[x,y]_P,x,y\in\g$ (in other words $T_W=\d P$ in the sense of the cohomology with coefficients in the adjoint module).

\end{defi}

If $(\g,[,])$ is semisimple, its second cohomology is trivial, hence any WNO possesses a primitive. Below we will prove the existence of such an operator directly (see Theorem \ref{4}). Note that the primitive of a WNO is defined up to addition of a derivation of the Lie algebra $(\g,[,])$.

Similarly, given a WNO $W\in L(\g)$ and any derivation $d\in \Der(\g)$ of the Lie algebra $(\g,[,])$, one observes that $W+d$ is again a WNO (since $[,]_{W+d}=[,]_W$). Assuming $W$ has a primitive $P$, what is a primitive of $W+d$?  Below we will need to know how to express a primitive of $W+d$ in terms of $P$ in a particular case of inner derivation. To find such an expression let us introduce another useful notion.

\abz\label{10.33}
\begin{defi}\rm A formal power series $R(t)=I+tR_1+t^2R_2+\cdots\in L(\g)[[t]]$ is called a {\em formal resolving function of order $n$} for a WNO $W$ if
\begin{equation}\equ\label{resolv1}
[R(t)\cdot,R(t)\cdot]\stackrel{n}{=}R(t) ([\cdot,\cdot]+
t[\cdot,\cdot]_W),
\end{equation}
 where $\stackrel{n}{=}$ stands for the equality modulo terms of order $> n$ in $t$ (here we understand the left hand side and the right hand side of the equality as elements of the space $L^2(\g)[[t]]$ of the formal power series with coefficients in $L^2(\g)$ and assume that $[,]$ is bilinear with respect to  $t$).

 We say that $\Phi(t)=I+t\Phi_1+t^2\Phi_2+\cdots\in L(\g)[[t]]$ is {\em a formal automorphism}  of the Lie algebra $(\g,[,])$ if
$$
[\Phi(t)\cdot,\Phi(t)\cdot]{=}\Phi(t) [\cdot,\cdot].
$$
\end{defi}

Clearly, if $\Phi(t)=I+t\Phi_1+\cdots$ is a formal automorphism, then $\Phi_1\in\Der(\g)$.

\abz\label{10.34}
\begin{theo}
\begin{enumerate}
\item A series $R(t)=I+tR_1+\cdots$ is a formal resolving function  of order  $1$ for a WNO $W$ if and only if $R_1=W+d$ for some $d\in\Der(\g)$.
     \item A formal series $R(t)=I+t(W+d)+t^2R_2+\cdots$ is a formal resolving function  of order  $2$ for a WNO $W$ if and only if $-R_2$ is a primitive of the WNO $W+d$.
\item Let $R(t)$ be a formal resolving function of order $n$ for a WNO $W$. If $\Phi(t)=I+t\Phi_1+t^2\Phi_2+\cdots$ is a formal automorphism of the Lie algebra $(\g,[,])$, then
    $    \tilde R(t)=\Phi(t)R(t)$
    is also a formal resolving function of order $n$ for $W$
\item  In particular, if $P$ is a primitive for $W$, then $P-\Phi_1W-\Phi_2$ is a primitive for $W+\Phi_1$.
\end{enumerate}
\end{theo}

\medskip \noindent{\em Ad  1, 2.} The proof follows by comparing of coefficients of $t$ and $t^2$ in the left hand side  and the right hand side of (\ref{resolv1}).

\medskip\noindent{\em Ad  3.}  The proof follows from the definitions.

\medskip\noindent{\em Ad  4.} If $R(t)=I+tW+t^2R_2+\cdots$, then $\tilde R(t)=I+t(W+\Phi_1)+t^2(R_2+\Phi_1W+\Phi_2)+\cdots$ and the result follows from Item 2. \qed

\abz\label{10.33a}
\begin{coro}
If $P$ is a primitive of a WNO $W$, the operator $P-\ad x_0W-(\ad^2 x_0)/2$ is a  primitive of the WNO $W+\ad x_0,x_0\in\g$. In particular, the WNO $\ad x_0$ has a primitive $-(\ad^2 x_0)/2$.
\end{coro}

\noindent Use the theorem and the fact that $\Phi(t):=\exp(t\ad x_0)=I+t\ad x_0+t^2(\ad^2 x_0)/2+\cdots$ is a formal automorphism of the Lie algebra  $(\g,[,])$. \qed

 \abz\label{10.35}
 \begin{exa}\rm Let $\g$ be decomposed to a sum of two subalgebras $\g=\g_1\oplus\g_2$. Put $W|_{\g_i}=\la_i \Id_{\g_i}, i=1,2$, where $\la_{1,2}$ are any scalars. Then it is easy to see that the Nijenhuis torsion of $W$ vanishes (cf. \cite{p5}). This way one gets a lot of examples of bi-Lie structures (for instance, taking $\g$ to be semisimple, $\g_1$ to be a parabolic subalgebra and $\g_2$ its "complement", cf. Theorem \ref{130.60}).
\end{exa}

\abz\label{10.40}
 \begin{exa}\rm Let $\g=\g_0\oplus \cdots\oplus\g_{n-1}$ be a $\Z_n$-grading on $\g$. Put $W|_{\g_i}=i \Id_{\g_i}, i=0,\ldots,n-1$. We will show that $W$ is a WNO calculating explicitly its primitive $P$.

Put $P|_{\g_i}=\frac{1}{2}i(n-i)\Id_{\g_i}$. Show that it is indeed a primitive of $W$. Let  $x_i\in\g_i,x_j\in\g_j$.  Then, if $i+j\le n-1$, we get
 $
 T_{W}(x_i,x_j)=[Wx_i,Wx_j]=ij[x_i,x_j]
 $
 (the term $W[x_i,x_j]_{W}$ vanishes). On the other hand,  we have $[x_i,x_j]_{P}=\frac{1}{2}(i(n-i)+j(n-j)-(i+j)(n-i-j))[x_i,x_j]=ij[x_i,x_j]$. If
$i+j\ge n$, we get
$
 T_{W}(x_i,x_j)=[Wx_i,Wx_j]-W[x_i,x_j]_{W}=(ij-(n-i-j)(i+j-(n-i-j)))
 [x_i,x_j]=(ij+(n-i-j)n)[x_i,x_j]
$
and $[x_i,x_j]_{P}=\frac{1}{2}(i(n-i)+j(n-j)-(i+j-n)(2n-i-j))[x_i,x_j]=
(ij-(i+j-n)n)[x_i,x_j]$.

\end{exa}
 \abz\label{10.50}
 \begin{exa}\rm Consider $\g=\gl(n,\C)$ and the operator $L_A, L_Ax:=Ax$, of the left multiplication by a fixed matrix $A\in \g$. Then the Nijenhuis torsion of $L_A$ vanishes, i.e. it is a WNO, and the modified bracket $[,]_{L_A}$ coincides with the bracket $[,_A]$ from Example \ref{0.20}.
\end{exa}

\abz\label{10.60}
 \begin{exa}\rm Consider $\g=\so(n,\C)$ and $A\in\s(n,\C)$. The operator $L_A$ is not correctly defined on $\g$ (since $L_A\g\not \subset\g$).  However, the operator $W:=W'|_\g$, where $W':=(1/2)(L_A+R_A)\in\End(\gl(n,\C))$ and $R_A$ stands for the operator of the right multiplication by $A$, is a WNO and it easy to see that again $[,]_{W}=[,_A]$ (see Example \ref{0.30}). The torsion of $W$ does not vanish, but it is a cocycle.  Indeed, since $W'=L_A-(1/2)(L_A-R_A)=L_A-\ad (A/2)$ and $0$ is a primitive of $L_A$, by Corollary \ref{10.33a} we conclude that the operator $\ad (A/2)\circ L_A-(1/8)\ad^2 A$ is a primitive of $W'$. Thus $T_{W'}$ is a cocycle and so is $T_W$.
\end{exa}

\section{Partial Nijenhuis operators}
\label{15}

Let us start from  a different look at Example \ref{10.60}. Put $\G:=\gl(n,\C)$ and consider the operator  $L_A$ as an operator acting from  $\g=\so(n,\C)$ to $\G$. It is easy to see that the formula similar to (\ref{e10.10}), i.e.
$$
[\cdot,\cdot]_{L_A}:=[L_A \cdot,\cdot]+[\cdot,L_A \cdot]-L_A[\cdot,\cdot],
$$
still defines a bilinear operation on $\g$ (in fact $[x,y]_{L_A}=[x,_Ay]$)  and, moreover, $T_{L_A}(x,y)=[L_Ax,L_Ay]-L_A[x,y]_{L_A}$ vanishes for any $x,y\in\g$. This motivates the following definition.

\abz\label{PNO}
\begin{defi}\rm Let $\G$ be a Lie algebra and $\g \subset \G$  a Lie subalgebra. We say that a pair $(\g,N)$, where $N\colon \g\to\G$ is a linear  operator,  is a {\em partial Nijenhuis operator on $\G$} (PNO for short) if the following two conditions hold:
\begin{enumerate}
\item[(i)] $[x,y]_N\in\g$ for any $x,y\in\g$;
\item[(ii)] $T_N(x,y)=0$ for any $x,y\in\g$.
\end{enumerate}
(Here $[,]_N$ and $T_N$ are  given by the same formulae as in Definitions \ref{10.10}, \ref{10.20}; note that it follows from condition (i) that the  term $N[x,y]_N$ which appears in the definition of $T_N$ is correctly defined.)
\end{defi}

\abz\label{cgm}
\begin{rema}\rm A similar structure in the context of Lie algebroids appeared in \cite{cgm2} under the name "outer Nijenhuis tensor".
\end{rema}

\abz\label{PNOlem}
\begin{lemm}
If $(\g,N)$ is a PNO on $\G$, then:
\begin{enumerate}
\item $N\g$ is a Lie subalgebra in $\G$;
\item $(\g,N^t)$, $N^t:=N-tI$, is a partial Nijenhuis operator on $\G$ for any $t\in\overline\C:=\C\cup\{\infty\}$,  here $I\colon \g\to\G$ is the natural embedding and $N^\infty=I$ by definition;
\item $N^t\g$ is a Lie subalgebra in $\G$ for any $t$;
\item for any $t\in\overline\C$ the operation $[,]_{N^t}$ is a Lie algebra structure on $\g$ and $N^t\colon \g\to\G$ is a homomorphism between Lie algebras $(\g,[,]_{N^t})$ and $(\G,[,])$;
\item the Lie bracket $[,]_N$ is compatible with the Lie bracket $[,]$.
\end{enumerate}
\end{lemm}
Indeed, Item 1 follows from the definition of a PNO. Item 2 is due to the equality $[,]_{N^t}=[,]_N-t[,]$ and to the equality $T_{N-tI}=T_N$. Item 3 follows from Item 2.

Now Item 4  follows easily from the equality $[x,y]_{N-tI}=(N-tI)^{-1}[(N-tI)x,(N-tI)y]$, which makes sense for almost all $t$, and  Item 5 is a consequence of Item 4.  \qed

\medskip
\abz\label{exaPNO}
\begin{exa}\rm The basic example here is $L_A:\so(n,\C)\to\gl(n,\C), A\in \s(n,\C), [,]_{L_A}=[,_A]$.
\end{exa}
Below we will associate a PNO with any bi-Lie structure $(\g,[,],[,]')$ with a semisimple bracket $[,]$.

\section{PNO and WNO in the presence of orthogonal decomposition}
\label{30}

Let $\G$ be a Lie algebra, $\tilde\g \subset \G$ a Lie subalgebra and let $B$ be an invariant form on $\G$ (we put tilde over $\g$ to be consistent with the notations of the next section).
\abz\label{1}
\begin{lemm} $[\g_0,\g_1]\subset \g_1$, here $\g_0:=\tilde\g, \g_1:=\tilde\g^\bot$ is the orthogonal complement to $\tilde\g$
with respect to $B$.
\end{lemm}

\noindent Let $x,z\in\g_0,y\in\g_1$. Then  $B([x,y],z)=-B(y,[x,z])=0$. Since $z\in\g_0$ is arbitrary, we get the result. \qed

\medskip

From now on we will assume that
\begin{equation}\equ\label{eq1}
\g_0\oplus\g_1=\G.
 \end{equation}
Such decomposition occurs for instance when $B$ is nondegenerate and so is $B|_{\tilde\g\times\tilde\g}$. The main example which we have in mind is the pair $\tilde\g=\ad\g \subset\End(\g)$, where $\g$ is semisimple and $B$ is the trace form on $\End(\g)$ (see Section \ref{40}; cf. also Example \ref{30.10} for an instance of different kind).

Now let  $N:\tilde\g\to\G$, be a linear operator. Put $p_i:\G\to\g_i$ for the corresponding orthogonal projection, $N_i:=p_iN,i=1,2$. Then
$N=N_0+N_1$ and $[,]_N=[,]_{N_0}+[,]_{N_1}$. We obviously have
$[\g_0,\g_0]_{N_0}\subset\g_0$. Assuming that $[\g_0,\g_0]_{N}\subset\g_0$ (the first condition of the definition of PNO),  we conclude from the lemma that
$$
[\g_0,\g_0]_{N_1}=\{[N_1x,y]+[x,N_1y]+N_1[x,y]\mid x,y\in\g_0\}\subset \g_1
$$
and $[,]_N=[,]_{N_0}$ on $\g_0$. Hence, if also the second condition of the definition of PNO holds for $N$, the operator  $N_0$ is a WNO on $\tilde\g$ by Item 4 of Lemma \ref{PNOlem}, which is uniquely
(!) defined by $N$ and which gives the same modified bracket on $\tilde\g$ as $N$. Let us summarize this information in the following lemma.

\abz\label{30.05}
\begin{lemm}
Let $N:\tilde\g\to\G$ be a PNO. Assume that an invariant form $B$ on $\G$ is given and condition  (\ref{eq1}) is satisfied, where $\g_0=\tilde\g,\g_1=\tilde\g^\bot$. Then the operator $N_0:=p_0N$ (here $p_0:\G\to\g_0$ is the orthogonal projection) is a WNO on $\tilde\g$ and moreover $[,]_N=[,]_{N_0}$.
\end{lemm}

\abz\label{30.10}
\begin{exa}\rm Let $\tilde\g=\g_0=\so(n,\C),\G=\gl(n,\C),\g_1=\s(n,\C), N=L_A|_{\tilde\g},A\in\s(n,\C)$ (see Example \ref{exaPNO}), and let $B$ be the trace form on $\G$. Then $N_0=(1/2)(L_A+R_A)|_{\tilde\g},N_1=(1/2)(L_A-R_A)|_{\tilde\g}$. Thus we get the WNO $W=N_0$ from Example \ref{10.60}.
\end{exa}

Let us study the pair of operators $(N,N_0)$, where $N$ is a PNO,  in more details.
First, we are able to express the Nijenhuis torsion of the WNO $N_0$ by means of the operator $N_1$. Indeed, using Lemma \ref{1} we get
$0=p_0T_N(x,y)=p_0([Nx,Ny]-N[x,y]_N)=[N_0x,N_0y]+p_0[N_1x,N_1y]-N_0([x,y]_{N_0})=
T_{N_0}+p_0[N_1x,N_1y]$ for $x,y\in\g_0$, hence
\begin{equation}\equ\label{tor}
T_{N_0}(x,y)=-p_0[N_1x,N_1y],\ x,y\in \tilde\g.
\end{equation}

Second, the condition $[,]_N=[,]_{N_0}$, shows that $[,]_{N_1}=0$ on $\tilde\g$. In other words, the map $N_1:\tilde\g\to\g_1$ is a 1-cocycle on $\tilde\g$ with coefficients in the $\tilde\g$-module $\g_1$.
Thus, if we assume  additionally that $\tilde\g$ is semisimple, by the Whitehead lemma $N_1=\d w$ for some (uniquely defined) element $w\in\g_1$, i.e.
$$
N_1=-(\ad_\G w)|_{\tilde\g}.
$$
By Corollary \ref{10.33a} we have $[\ad w (x),\ad w(y)]=T_{\ad w}(x,y)=[x,y]_{-(\ad^2w)/2}$, hence
$$
T_{N_0}(x,y)=-p_0[x,y]_{-(\ad_\G^2w)/2}=[x,y]_{p_0(\ad_\G^2w)/2},\ x,y\in\tilde\g.
$$
In other words,  $P:=p_0(\ad_\G^2w)/2$ is a primitive of the WNO $N_0$. Note that $P$ is a symmetric operator with respect to the form $B|_{\tilde\g\times\tilde\g}$. Indeed, $\ad_\G w$ is antisymmetric, hence, given $x,y\in\tilde\g$, we have $B(Px,y)=(1/2)B((\ad_\G^2w)(x),y)=
(1/2)B(x,(\ad_\G^2w)(y))=B(x,Py)$. Let us summarize this in

\abz\label{30.20}
\begin{lemm}
Retain the hypotheses of Lemma \ref{30.05} and assume moreover that $\tilde\g$ is semisimple. Then there exists a uniquely defined element $w\in \g_1$ such that $N_1:=p_1N=-\ad_\G w$. Moreover the operator $P:=(1/2)p_0(\ad_\G^2w):\tilde\g\to\tilde\g$ is symmetric with respect to $B|_{\tilde\g\times\tilde\g}$ and is a primitive of the WNO $N_0$.
\end{lemm}

\section{Semisimple bi-Lie structures and their principal WNO}
\label{40}

\abz\label{40.0}
\begin{defi}\rm
We say that a bi-Lie structure $(\g,[,],[,]')$ is {\em semisimple} if $(\g,[,])$ is a semisimple Lie algebra (which will be called the {\em underlying Lie algebra} of the bi-Lie structure).
\end{defi}
We  want to apply the construction from Lemma \ref{30.05} to a situation in which:
\begin{enumerate}
\item $(\g,[,],[,]')$ is a semisimple bi-Lie structure;
     \item $\G:=\End(\g)$;
     \item $\tilde\g:=\ad{\g}\subset\G$;
     \item $N:=\ad'(\ad)^{-1}:\tilde\g\to\G$;
     \item $B(a,b):=\Tr(ab),\ a,b\in\G$;
\end{enumerate}
 (here $\ad x(y):=[x,y],\ad'x(y)=[x,y]',x,y\in\g$).
 It is easy to see  that $B|_{\tilde\g\times\tilde\g}$ is the Killing form of $(\tilde\g,[,])$ (here $[,]$ is the commutator of endomorphisms). It is nondegenerate (due to the semisimplicity of $\g$) and,
moreover, $B$ is nondegenerate itself, hence  condition (\ref{eq1}) is satisfied.

The following lemma shows that $N$ is a PNO on $\G$.

\abz\label{3}
\begin{lemm} \begin{enumerate}\item $[\ad x,\ad y]_{N}=\ad {[x,y]'},x,y\in\g$;\item $T_N(\ad x,\ad y)=0,x,y\in\g$. \end{enumerate}
 \end{lemm}

\noindent{\em Ad  1.}  By definition, $N(\ad x)=\ad' x$. Thus $[\ad x,\ad y]_{N}=[N(\ad x),\ad y]+[\ad x,N(\ad y)]- \linebreak N[\ad x,\ad y]
=[\ad' x,\ad y]+[\ad x,\ad' y]-N\ad {[x,y]}=
[\ad' x,\ad y]+[\ad x,\ad' y]-\ad' {[x,y]}$. The last expression  is equal to $\ad {[x,y]'}$ due to the compatibility of the brackets $[,],[,]'$ (see Lemma \ref{0.40}).

\medskip

\noindent{\em Ad  2.} $T_N(\ad x,\ad y)=[N\ad x,N\ad y]-N[\ad x,\ad y]_N=[\ad'x,\ad' y]-N\ad[x,y]'=[\ad'x,\ad' y]-\ad'[x,y]'=0$. \qed

\medskip

The construction from the previous section gives a WNO $\widetilde W:\tilde\g\to\tilde\g, \widetilde W:=N_0=p_0\circ N$.
Let us  denote by $W$ the corresponding WNO on $\g$ induced by $\widetilde W$ and by the isomorphism of Lie algebras $\ad:(\g,[,])\to(\tilde\g,[,])$, in other words, $\widetilde W(\ad x)=\ad {Wx},x\in \g$. It follows from the lemma above  that under the identification "$\ad$" the second bracket $[,]'$ corresponds to the bracket $[,]_N=[,]_{\widetilde W}$, hence $[,]'=[,]_{W}$.

Summarizing we get the first assertion of the following theorem.

\abz\label{40.th}
\begin{theo} \begin{enumerate}\item
Given a  semisimple bi-Lie structure $(\g,[,],[,]')$, the operator $W=\ad^{-1}\circ p_0\circ\ad':\g\to\g$, where $p_0:\G\to\tilde\g$ is the orthogonal projection with respect to $B$, is a WNO  on the Lie algebra $(\g,[,])$ such that
$$
[,]'=[,]_{W}.
$$
\item Among all WNOs $V\in\End(\g)$ such that $[,]'=[,]_{V}$ there is a unique WNO $ V'$ satisfying the condition
    $ V'\in\tilde\g^\bot$.
\item $W=V'$.
    \end{enumerate}
\end{theo}

\noindent We have to justify only items 2,3. If $V_1:\g\to\g$ is a WNO with the property $[,]_{V_1}=[,]_{V}$, then the difference $V-V_1$ obviously should be a derivation of the bracket $[,]$. Since all the derivations of $[,]$ are inner, $V$ and $V_1$ should differ by  some $\ad x, x\in\g$, i.e. by an element of $\tilde\g=\g_0$. Now the decomposition property (\ref{eq1}) implies the uniqueness of WNO belonging to $\tilde\g^\bot=\g_1$.

Item 3 follows from Corollary \ref{5} below. \qed

\medskip
\abz\label{defiMain}
\begin{defi}\rm Let   $(\g,[,])$ be a semisimple Lie algebra. An operator $V\in\End(\g)$ is said to be {\em principal} if $V\in\tilde\g^\bot$, here $\tilde\g=\ad \g$.

Given a semisimple bi-Lie structure $(\g,[,],[,]')$,  any   $W\in\End(\g)$ such that $[,]'=[,]_W$ will be called a {\em WNO of the bi-Lie structure} and the unique (in view of Theorem \ref{40.th}) WNO of the bi-Lie structure belonging to $\g^\bot$ will be called the {\em principal WNO} of this bi-Lie structure.
\end{defi}

The next result gives an explicit expression of the orthogonal decomposition $N=N_0+N_1$  (cf. Section \ref{30}) for $N=\ad'(\ad)^{-1}$ in terms of the operator $W$  and calculates one of the primitives of the  WNO ${W}$.

\abz\label{4}
\begin{theo}\begin{enumerate}\item
The action of the operator $N:\tilde\g\to\G$ can be given by the following expression: $N(\ad x)=\ad {Wx}+[\ad x,W]$; moreover, $N_0(\ad x)=\widetilde{ W}(\ad x)=\ad {Wx},N_1(\ad x)=[\ad x,W],x\in\g$ (in other words $N_1=L|_{\tilde\g}$, where $L:\G\to\G$ is equal to $-\ad_\G W$).
\item
The operator $\widetilde{ P}:=(1/2)p_0\circ\ad_\G^2 {W}:\tilde\g\to\tilde\g$ is symmetric with respect to the Killing form $B_{\tilde\g}$ and is a primitive of the WNO $\widetilde{ W}$. The formula $P=\ad^{-1}\circ \widetilde{ P}\circ\ad$ gives a primitive of the WNO $W$, which is symmetric with respect to the Killing form $B_\g$.
\end{enumerate}
\end{theo}

\noindent By the definition of $N$ for any $x,y\in\g$ we have $(N(\ad x))y=\ad' x(y)=[x,y]'$. Moreover, we know that $[,]'=[,]_{W}$, hence $(N(\ad x))y=[{W}x,y]+[x,{W}y]-{W}[x,y]$. Rewrite this equality as $(N(\ad x))y=\ad {Wx}(y)+[\ad x,W]y$.
Since the first term of this expression is equal to $(\widetilde W(\ad x))y$ (by the definition of ${W}$), the second one has to be equal $((N-\widetilde W)(\ad x))y=(N_1(\ad x))y$. The rest of the proof follows from Lemma \ref{30.20}. \qed

\abz\label{prpr}
\begin{defi}\rm The operator $P$ defined in Theorem \ref{4} will be called the {\em principal primitive} of the principal WNO $W$ (note that the condition of being symmetric distinguishes the principal primitive uniquely among all primitives of the WNO $W$).
\end{defi}

\abz\label{coro2}
\begin{coro} An element $x\in\g$ belongs to the centre of a Lie bracket $[,]'$ if and only if $x\in\ker W$ and $[\ad x,W]=0$, where $W$ is the principal WNO of the bi-Lie structure $(\g,[,],[,]')$.
\end{coro}

\noindent Indeed, by definition the images of the operators $\widetilde W,N_1$ are mutually orthogonal. Thus $0=\ad'x=N(\ad x)$ if and only if $\widetilde W(\ad x)=0$ and $N_1(\ad x)=0$ if and only if $Wx=0$ and $[\ad x,W]=0$. \qed

\abz\label{5}
\begin{coro} The endomorphism $W\in\End(\g)$ defined before Theorem \ref{40.th}  belongs to $\g_1=\tilde\g^\bot$.
\end{coro}

\noindent We know from Theorem \ref{4} that $[\ad x,W]\in\g_1$ for any $x\in\g$, i.e. $B(\ad y,[\ad x,W])=0$ for any $x,y\in\g$. By the invariance property of $B$ we have $B(W,[\ad x,\ad y])=0$. However, $(\g,[,])$ is semisimple and coincides with its commutant, hence the last equality means that $W\in\g_1$. \qed

\medskip

\abz\label{coroX}
\begin{coro} Let $(\g,[,],[,]_1')$ and $(\g,[,],[,]_2')$ be two bi-Lie structures with semisimple $(\g,[,])$ and let $W_1,W_2$ be the corresponding principal WNOs on $\g$. Then the bi-Lie structures are strongly isomorphic (see Definition \ref{0.15}) if and only if there exists an automorphism $\phi$ of the Lie algebra $(\g,[,])$ with the property $\phi\circ W_1=W_2\circ \phi$.
\end{coro}

\noindent Assume first that $\phi$ is an automorphism with the prescribed property. Then direct calculation shows that it transforms the bracket $[,]_{W_1}$ to the bracket $[,]_{W_2}$.

Now assume $\phi$ is an automorphism of $(\g,[,])$ with the property $\phi[\cdot,\cdot]_{W_1}=[\phi\cdot,\phi\cdot]_{W_2}$.

Note that, given Lie algebras $(\g,[,]_1)$ and $(\g,[,]_2)$, a linear automorphism $\phi\in\End(\g)$ transforms the bracket $[,]_1$ to $[,]_2$ if and only if $\phi\circ\ad_1 x\circ\phi^{-1}=\ad_2 \phi x$ for any $x\in\g$, here $\ad_ix(y):=[x,y]_i$.

Thus, by the assumption, $\phi\circ\ad x\circ\phi^{-1}=\ad\phi x$ and
\begin{equation}\equ\label{equi}
 \phi\circ \ad^{W_1}x\circ\phi^{-1}=\ad^{W_2}\phi x,
  \end{equation}
  here $\ad^{W_i}x(y):=[x,y]_{W_i}$.

 We know that $\ad^{W_2}\phi x=\ad W_2\phi x+[\ad\phi x,W_2]$, where the first term belongs to $\tilde\g$ and the second one to $\tilde\g^\bot$.

 On the other hand, $\phi\circ \ad^{W_1}x\circ\phi^{-1}=\phi\circ \ad W_1x\circ\phi^{-1}+\phi\circ [\ad x,{W_1}]\circ\phi^{-1}=
\ad \phi W_1x+[\ad\phi x,\phi\circ W_1\circ\phi^{-1}]$. We claim that the first term belongs to $\tilde\g$ (which is obvious) and the second one to $\tilde\g^\bot$. Indeed $B(\ad x,\phi\circ W_1\circ\phi^{-1})=\Tr(\ad x\circ\phi\circ W_1\circ\phi^{-1})=\Tr(\phi^{-1}\circ\ad x\circ\phi\circ W_1)=B(\ad\phi^{-1}x,W_1)$. The last expression is equal to zero by Corollary \ref{5}. Thus $\phi\circ W_1\circ\phi^{-1}\in\tilde\g^\bot$ and by Lemma \ref{1} $[\ad\phi x,\phi\circ W_1\circ\phi^{-1}]\in\tilde\g^\bot$.

Hence equality (\ref{equi}) implies the equality of $0$-components, $\ad W_2\phi x=\ad \phi W_1x$, for any $x\in\g$, which, in view of the injectivity of $\ad$, completes the proof. \qed

\medskip

Let us summarize the discussion above in the following theorem.

\abz\label{40.nonum}
\begin{theo} There is a one-to-one correspondence between semisimple bi-Lie structures $(\g,[,], \lbr [,]')$ and principal WNOs $W\in\End(\g)$ given by $[,]'=[,]_W$. Two semisimple bi-Lie structures $(\g,[,],[,]')$ and $(\g,[,],[,]_1')$ with the same underlying  Lie algebra are (strongly) isomorphic (see Definition \ref{0.15}) if and only if the corresponding principal WNOs $W,W_1$ are (strongly) equivalent in the sense of the following definition.
\end{theo}

\abz\label{isoWNO}
\begin{defi}\rm Let $(\g,[,])$ be a semisimple Lie algebra and let $W,W_1\in\End(\g)$ be two principal WNOs. We say that they are {\em strongly equivalent (equivalent)} if there exists an automorphism $\phi$ of the Lie algebra $(\g,[,])$ such that $\phi\circ W\circ \phi^{-1}=W_1$ (respectively $\phi\circ W\circ \phi^{-1}=\la W_1+\la'\Id_\g$ for some $\la,\la'\in\C$).
\end{defi}

We  conclude this section by an example of a principal WNO.

\abz\label{40.exa}
\begin{exa}\rm Consider the operator $W$ from Example \ref{10.60}. We claim that it is the principal WNO for the corresponding  bi-Lie structure $(\so(n,\C),[,],[,_A])$. Indeed $\Tr(W(X)Y)=(1/2)\Tr((AX+XA)Y)=(1/2)\Tr(X(AY+YA))
=\Tr(XW(Y))$, i.e. the operator $W$ is symmetric with respect to the Killing form $B_{\so(n,\C)}$ which is proportional to the trace form. It remains to use the following lemma.
\end{exa}

\abz\label{40.lemm}
\begin{lemm} Let an operator $W\in\End(\g)$ be given on a Lie algebra $(\g,[,])$ such that $W$ is symmetric with respect to the Killing form $B_\g$. Then $W$ is orthogonal to any antisymmetric operator with respect to the trace form $B$ on $\End(\g)$, in particular $W$ is principal.
\end{lemm}

\noindent Let  $A:\g\to\g$ be an antisymmetric operator. If $e_1,\ldots,e_n$ is an orthonormal with respect to $B_\g$ basis in $\g$, the trace of any operator $L$ can be calculated as $\Tr(L)=\sum_{i=1}^nB_\g(Le_i,e_i)$. Now, $\Tr(WA)=\sum_{i=1}^nB_\g(WAe_i,e_i)=-\sum_{i=1}^nB_\g(e_i,AWe_i)=
-\Tr(AW)=-\Tr(WA)$. Hence $\Tr(WA)=0$.
\qed

\section{Kernels of the Killing forms of the exceptional brackets and their centres}
\label{90}

In this section we assume that  $(\g,[,],[,]')$ is a semisimple bi-Lie structure and  $W:\g\to\g$ is its principal WNO. As usual, we will write $B$ for the trace form on $\G:=\End(\g)$ and $B_\g$ for the Killing form of the Lie algebra $(\g,[,])$.

Consider the bracket $[,]^t:=[,]'-t[,]=[,]_{W^t},t\in \C$, here $W^t:=W-t\Id_\g$. Let $\ad^t:\g\to\g$ stand for the corresponding $\ad$-operator, $\ad^tx(y)=[x,y]^t$. This operator  can be also written in the form $\ad^tx=\ad W^tx+[\ad x,W]$ (since $[\ad x,W^t]=[\ad x, W]$ for any $x\in\g$, see Theorem \ref{4}).

\abz\label{90.10}
\begin{theo} The Killing form $B^t$ of the Lie algebra $(\g,[,]^t)$ is given by the formula
$$
B^t(x,y)=B_\g(W^tx,W^ty)+b(x,y),\ x,y\in\g,
$$
where  $b$ is a symmetric form on $\g$ equal to
$$
b(x,y)=B([\ad x,W],[\ad y,W])=-2B_\g(Px,y);
$$
here  $P$ is the principal  primitive of the principal WNO $W$ (see Definition \ref{prpr}).
\end{theo}

\noindent By definition $B^t(x,y)=B(\ad^tx,\ad^ty)=B(\ad W^tx+[\ad x,W],\ad W^ty+[\ad x,W])=B(\ad W^tx, \lbr \ad W^ty)+B([\ad x,W],[\ad y,W])=B_\g(W^tx,W^ty)+b(x,y)$. Here we used the fact that the endomor- \linebreak phisms $\ad W^tz,z\in\g$, and $[\ad z,W]$ are mutually orthogonal with respect to $B$ (see Corollary \ref{5} and Lemma \ref{1}). The equality $
b(x,y)=-2B_\g(Px,y)
$
follows from the definition of $P$ and from the invariance of the form $B$.
 \qed

The theorem above implies (due to the nondegeneracy of $B_\g$) that an element $x\in\g$ belongs to the kernel $\ker B^t$ of the Killing form $B^t$ for some $t\in\C$ if and only if $x$ belongs to the kernel of the operator
$
\hat B^t:=(W^t)^*W^t-2P\in\End(\g)
$, or in more details,
\begin{equation}\equ\label{killf}
\hat B^t=W^*W-2P-t(W+W^*)+t^2\Id_\g;
\end{equation}
here $L^*$ stands for the operator adjoint to an operator $L$ with respect to $B_\g$. In particular, the set $T=\{t\in\C\mid \ker B^t\not=\{0\}\}$ is finite.

\abz\label{defiTimes}
\begin{defi}\rm The elements of the set  $T=\{t_1,t_2,\ldots\}$ will be called the {\em times} of a bi-Lie structure $(\g,[,],[,]')$ and the corresponding brackets $[,]^{t_i}$ will be called {\em exceptional}.
\end{defi}

The following result describes the centres of the exceptional brackets.

\abz\label{t10.1}
\begin{theo}\begin{enumerate}
\item An element $x\in\g$ belongs  to the centre $\z^t$ of the bracket $[,]^t:=[,]'-t[,]=[,]_{W^t}$ if and only if
$x$ is an eigenvector of the principal WNO $W$ corresponding to the eigenvalue $t$ and $[\ad x,W]=0$;
\item The subset $\z^t$ is a subalgebra in $(\g,[,])$ for any $t$ (this result is known \cite[\S 44]{TFbols});
\item $\z^{t_1}\cap\z^{t_2}=\{0\}$ if $t_1\not=t_2$;
\item $[\z^{t_1},\ker W^{t_2}]\subset\ker W^{t_2}$ for any $t_1,t_2\in\C$;
    \item $[\z^{t_1},\z^{t_2}]=0$ if $t_1\not=t_2$; in particular, the set $\z:=\z^{\theta_1}\oplus \cdots\oplus\z^{\theta_m}$ is a subalgebra in $(\g,[,])$  which is a direct sum of its ideals $\z^{\theta_i}$, here $\Theta:=\{\theta_1,\ldots,\theta_m\}=\{\theta \in\C\mid \z^\theta\not=\{0\}\}\subset T$.
\end{enumerate}
\end{theo}

\noindent{\em Ad  1.} The proof is a direct consequence  of  Corollary \ref{coro2} and the fact that $[\ad x,W^t]=[\ad x, W]$.

\noindent{\em Ad  2.}  For any $x,y\in\z^\theta \subset\ker W^\theta$ we have $0=[x,y]_{W^\theta}=-W^\theta[x,y]$, hence $[x,y]\in\ker W^\theta$. On the other hand, since $[\ad x,W]=0=[\ad y,W]$, then $[\ad [x,y],W]=[[\ad x,\ad y],W]=-[[W,\ad x],\ad y]-[[\ad y,W],\ad x]=0$. Thus by Item 1 $[x,y]\in\z^\theta$.

\medskip

\noindent{\em Ad  3.} The proof follows from the inclusion $\z^{t}\subset\ker W^{t}$.

\medskip

\noindent{\em Ad  4.} Let $x\in\z^{t_1},y\in\ker W^{t_2}$. Then $0=[x,y]_{W^{t_1}}=[x,W^{t_1}y]-W^{t_1}[x,y]=[x,W^{t_2}y]+
(t_2-t_1)[x,y]-W^{t_1}[x,y]=(t_2-t_1)[x,y]-W^{t_1}[x,y]
=t_2[x,y]-W[x,y]$. Hence $[x,y]\in\ker W^{t_2}$.

\medskip

\noindent{\em Ad  5.} We have $[\z^{t_1},\z^{t_2}]\subset [\z^{t_1},\ker W^{t_2}]\subset\ker W^{t_2}$ by Item 4 and, analogously, $[\z^{t_1},\z^{t_2}]\subset [\ker W^{t_1},\z^{t_2}]\subset\ker W^{t_1}$. Hence $[\z^{t_1},\z^{t_2}]\subset\ker W^{t_1}\cap\ker W^{t_2}=\{0\}$. The inclusion $\Theta \subset T$ is obvious, since the centre of any Lie algebra is contained in the kernel of the Killing form. \qed

\abz\label{defCentr}
\begin{defi}\rm The subalgebra $\z \subset\g$ equal to the sum of the centres of the exceptional brackets will be called the {\em central subalgebra} of a bi-Lie structure $(\g,[,],[,]')$.
\end{defi}
We conclude this section by introducing a series of another subalgebras of $(\g,[,])$.

\abz\label{overcentral}
\begin{theo} Let $(\g,[,],[,]')$ be a semisimple bi-Lie structure, $W:\g\to\g$ be its principal WNO, and $P$ its pricipal primitive.
\begin{enumerate}\item The sets $\g_P:=\{x\in\g\mid [\ad x,P]=0\}, \hat\z:=\g_P\cap\ker P$ are subalgebras in $(\g,[,])$ and $\z \subset\hat\z \subset\g_P$.
\item If $W$ is diagonalizable, then $\hat\z^t:=\hat\z\cap\ker W^t$ also is a subalgebra for any $t\in\C$, $\z^t \subset \hat\z^t$ and $\hat\z^{t_1}\cap\hat\z^{t_2}=\{0\}$ if $t_1\not=t_2$. If, moreover, $W\hat\z \subset\hat\z$, then $\hat\z=\hat\z^{t_1}\oplus \cdots\oplus\hat\z^{t_k}$ for some $t_i\in T$ and $\hat\z^{t_i}+\hat\z^{t_j}$ is a subalgebra for any $i,j$.
\end{enumerate}
\end{theo}

\noindent {\em Ad 1.} By the Jacobi identity  $[\ad[x,y],P]=[[\ad x,P],\ad y]-[[\ad y,P],\ad x]$. Thus $x,y\in\g_P$ implies $[x,y]\in\g_P$. If $x,y\in\hat\z$, then $P[x,y]=P\ad x(y)=\ad x(Py)=0$ and $[x,y]\in\hat\z$.

To prove the inclusion $\z \subset \hat\z$ rewrite the equality $T_W(x,y)=[x,y]_P$ in a different way. Namely, fix $x$ and consider left and right hand sides of this equality as operators on the second argument. Then we get $T_W(x,\cdot)=[Wx,W \cdot]-W([x,\cdot]_W)=(\ad Wx\circ W-W\circ (\ad Wx+[\ad x,W]))(\cdot)=([\ad Wx,W]-W\circ[\ad x,W])(\cdot)$. Analogously, $[x,\cdot]_P=(\ad Px+[\ad x,P])(\cdot)$. Finally,
\begin{equation}\equ\label{torsio}
[\ad Wx,W]-W\circ[\ad x,W]=\ad Px+[\ad x,P].
\end{equation}
Note also that $T_W=T_{W^t}$ for any $t\in\C$, hence $[\ad W^tx,W]-W^t\circ[\ad x,W]=\ad Px+[\ad x,P]$.
If $x\in\z$, then by Theorem \ref{t10.1} $x\in\ker W^t$ for some $t$ and $[\ad x,W]=0$, whence $\ad Px+[\ad x,P]=0$. Now it is enough to use Lemmata \ref{1}, \ref{40.lemm} to deduce that $\ad Px\in\ad\g,[\ad x,P]\in(\ad\g)^\bot$ have to be zero and that $x\in\hat\z$. The same argument proves the inclusion $\z^t \subset \hat\z^t$.

\medskip

\noindent {\em Ad 2.} Let $x,y\in\hat\z^t$. Then $T_W(x,y)=t^2[x,y]-W(2t[x,y]-W[x,y])=(W^t)^2[x,y]$ (since $x,y\in\ker W^t$). On the other hand, $T_W(x,y)=\ad Px(y)+[\ad x,P]y=0$ (since $x\in\hat\z$). The diagonalizability of $W$ implies $[x,y]\in\ker W^t$.

The equality $\hat\z^{t_1}\cap\hat\z^{t_2}=\{0\}$ follows from the equality $\ker W^{t_1}\cap\ker W^{t_2}=\{0\}$.

If $\hat\z$ is $W$-invariant, we have a direct decomposition $\hat\z=\hat\z^{t_1}\oplus \cdots\oplus\hat\z^{t_k}$. To prove the last assertion it is enough to observe that $W|_{\hat\z}$ due to formula (\ref{torsio}) has zero torsion, i.e. is Nijenhuis. The needed property is a characteristic property of any diagonalizable Nijenhuis operator \cite{p5}. \qed

\section{Operators preserving a grading}
\label{100}

We will consider  $\Ga$-{\em gradings} of a Lie algebra $\g$, i.e. direct decompositions $\g=\bigoplus_{i\in\Ga}\g_i$ such that $[\g_i,\g_j] \subset\g_{i+j}$ for any $i,j\in\Ga$, here $\Ga$ is an abelian group and we use the additive notation.

\abz\label{100.10}
\begin{defi}\rm We say that a linear operator $W\in\End(\g)$ preserves a $\Ga$-grading $\g=\bigoplus_{i\in\Ga}\g_i$ if $W{\g_i}\subset{\g_i}$ for any $i\in\Ga$.
\end{defi}

The first two items of the following theorem are standard. Let $\g$ be a semisimple Lie algebra with the Killing form $B_\g$ and let $\g=\bigoplus_{i\in\Ga}\g_i$ be a $\Ga$-grading.

\abz\label{100.20}
\begin{theo}
\begin{enumerate}\item $B_\g(\g_i,\g_j)=0$ if $i+j\not=0$.
\item The restriction of $B_\g$ to $\g_i\times\g_j$ is a nondegenerate pairing for all $i,j\in \Ga,i+j=0$. In particular so is the restriction of $B_\g$ to $\g_0$.
\item If $W\in\G=\End(\g)$ preserves the grading,  then the operator $\widetilde{ P}\in\End(\tilde\g),\tilde\g:=\ad\g$, given by the formula $\widetilde{ P}:=(1/2)p_0\circ\ad^2_\G W$ preserves the induced grading  on $\tilde\g$, here $p_0:\G\to\tilde\g$ is the orthogonal  projection onto $\ad\g$ with respect to the trace form $B$. Consequently, the operator $P:=\ad^{-1}\circ \widetilde{P}\circ\ad$ preserves the initial grading on $\g$.
\item If $W\in\End(\g)$ and, moreover,  $W|_{\g_i}=\la_i\Id_{\g_i},\la_i\in\C$, for any $i\in\Ga$, then the adjoint  with respect to $B_\g$ operator $W^*$ is given by $W^*|_{\g_i}=\la_{-i}\Id_{\g_i}$.
\end{enumerate}
\end{theo}

\noindent {\em Ad  1.} Let $x_i\in\g_i,x_j\in\g_j$. Consider a basis of $\g$ which can be divided to parts each of which forms a basis of one of the subspaces $\g_k$. Let $e^1,\ldots,e^m$ be such a part which is a basis of $\g_k$. Then for $A:=\ad x_i\circ\ad x_j$ we have $A(e^l)\not\in\g_k$ as soon as $i+j\not=0$. Hence the matrix of $A$ in such a basis has zero diagonal entries and $\Tr(A)=B_\g(x_i,x_j)=0$.

\medskip

\noindent {\em Ad  2.} Now let $x\in\g_i$ and $i+j=0$. Then there exists an element $y\in\g_j$ such that $B_\g(x,y)\not=0$. Indeed, otherwise $x$ would be orthogonal to all the space $\g$ by Item 1.

\medskip

\noindent {\em Ad  3.} Let $x_{i'}\in\g_{i'},x_{i''}\in\g_{i''}$. Since $W$ preserves the grading, we have  $\ad^2_\G W(\ad x_{i'})(x_{i''})=(W^2\circ\ad x_{i'}-2W\circ\ad x_{i'}\circ W+\ad x_{i'}\circ W^2)(x_{i''})\in\g_{i'+i''}$.

By Items 1,2 we can choose a hyperbolic basis $e_i^1,\ldots,e_i^{m_i},e_{-i}^1,\ldots,e_{-i}^{m_i}, e_i^j\in\g_i$, in each $\g_i+\g_{-i}$ such that $2i\not=0$. In $\g_i$ with $2i=0$ we will choose an orthonormal basis $e_i^1,\ldots,e_i^{m_i}$. The union of all these bases  forms a  basis of $\g$. If $X\in\G$, then $p_0X=\sum_{ij}x^j_i\ad e_i^j$, where $x^j_i=B(\ad e_{-i}^j,X)$. Put $X_{i'}^{j'}=\ad^2_\g W(\ad e_{i'}^{j'})$. Then $\ad e_{-i}^j\circ X_{i'}^{j'}(e_{i''}^{j''})\in\g_{-i+i'+i''}$. The nonzero diagonal entries of the matrix of the endomorphism $\ad e_{-i}^j\circ X_{i'}^{j'}$ in the basis $\{e_i^j\}$ can appear only if $i'=i$. Thus $x_{ii'}^{jj'}:=B(\ad e_{-i}^j, X_{i'}^{j'})=0$ if $i'\not=i$.

Summarizing,  we get $\widetilde{P}(\ad e_{i'}^{j'})=(1/2)\sum_{ij}x_{ii'}^{jj'}\ad e_i^j$, where $x_{ii'}^{jj'}=0$ for $i\not=i'$, which shows the invariance of $\g_i$ with respect to $\widetilde{P}$.

\medskip

\noindent {\em Ad  4.} It is enough to consider the operator $W$ in the basis built in the proof of Item 3.
\qed

\abz\label{100.25}
\begin{coro}
Let $W$ be a principal WNO of a semisimple bi-Lie structure (see Definition \ref{defiMain}).
\begin{enumerate}\item
Assume that $W$ preserves a $\Ga$-grading on $\g$. Then its principal primitive $P$ (see Definition \ref{prpr}) also preserves this grading.
\item If the initial grading is  the root decomposition grading with respect to a Cartan subalgebra $\h \subset\g$,
    \begin{equation}\equ\label{grading}
\g=\h+\sum_{\al\in R}\g_\al,
 \end{equation}
     (here $R$ stands for the set of roots) and $W$ preserves this grading, then $P|_{\g_\al+\g_{-\al}}=\pi_\al\Id_{\g_\al+\g_{-\al}}, \lbr \al\in R$, for some $\pi_\al\in\C$ and $P\h \subset\h$.
     \end{enumerate}
\end{coro}

\abz\label{100.260}
\begin{theo} Let $W\in\G=\End(\g)$ preserve the root grading (\ref{grading}), $W|_{\g_\al}=\la_\al\Id_{\g_\al},
\kappa_\al:=(1/2)(\la_\al-\la_{-\al})$. Then the operator $\pr(W):=p_1(W)\in\End(\g)$, where $p_1:\G\to\tilde\g^\bot$ is the orthogonal projection onto $\tilde\g^\bot$ along $\tilde\g$, is given by the formula
$$
\pr(W)=W- \ad (\sum_{\al\in R}\la_\al H_\al)=W- \ad (2\sum_{\al\in R^+}\ka_\al H_\al),
$$
here $R^+$ stands for the set of positive roots with respect to any basis and $H_\al\in\h$ is given by $B_\g(H_\al,H)=\al(H),H\in\h$.
The operator $\pr(V)$ will  be called the {\em principal projection} of $V$. In particular, $\pr(W)=W$ (i.e. $W$ is the principal operator in the sense of Definition \ref{defiMain}) if and only if $\sum_{\al\in R^+}\ka_\al\al=0$.
\end{theo}

\noindent We have to prove that the operator $\ad (2\sum_{\al\in R^+}\ka_\al H_\al)$ is equal to $p_0(W)$, where $p_0:\G\to\tilde\g$ is the orthogonal projection. To this end choose a basis in $\g$ as in the proof of Item 3 of Theorem \ref{100.20}, i.e. an orthonormal basis $e_1,\ldots,e_n$ in $\h$ and an element $E_\al\in\g_\al$ for any $\al\in R$ such that $B_\g(E_\al,E_{-\al})=1$. Moreover, let $f_1,\ldots,f_k$ be any basis in $\tilde\g^\bot$. Then we have to calculate the coefficients $a_i,b_\al$ of the decomposition  $W=\sum_ia_i\ad e_i+\sum_\al b_\al\ad E_\al+\sum_j c_jf_j$. Since $\ad e_i, \ad E_\al,f_j$ are mutually orthogonal, we have $a_i=B(\ad e_i,W), b_\al=B(\ad E_{-\al},W)$. The operator $\ad E_\al$ does not preserve any of the subspaces $\h, \g_\al$, hence $b_\al=\Tr(\ad E_{-\al}\circ W)=0$.

To calculate $a_i$, remark that $\ad e_i|_\h=0,\ad e_i|_{\g_\al}=\al(e_i)\Id_{\g_\al}$, whence $a_i=\Tr(\ad e_i\circ W)=\sum_{\al\in R}\al(e_i)\la_\al$. Thus $p_0(W)=\sum_ia_i\ad e_i=\sum_i\sum_{\al\in R}\al(e_i)\la_\al\ad e_i=\ad \sum_{\al\in R}\la_\al \sum_i B_\g(H_\al,e_i)e_i=\ad \sum_{\al\in R}\la_\al H_\al=2\ad \sum_{\al\in R^+}\ka_\al H_\al$. The last  equality is due the fact that $H_\al=-H_{-\al}$. \qed

\section{Regular semisimple bi-Lie structures}

\label{110}

Let $(\g,[,],[,]')$ be a semisimple bi-Lie structure, $W\in\End(\g)$ its principal WNO (see Definition \ref{defiMain}).  Let $\h \subset\g$ be a Cartan subalgebra of the semisimple Lie algebra  $(\g,[,])$ and $R=R(\g,\h) \subset\h^*_\R$ be the corresponding root system.

\abz\label{110.10}
\begin{theo} The following two conditions are equivalent:
\begin{enumerate}
\item The operator $W\in\End(\g)$ preserves the root grading (\ref{grading})
  and
the operator $W|_\h$ is diagonalizable.
\item $\h \subset\z$, where $\z$ is the central subalgebra (see Definition \ref{defCentr}).
  \end{enumerate}
If the conditions above are satisfied the central subalgebra $\z \subset\g$ is homogeneous with respect to the decomposition (\ref{grading}), i.e.
$$
\z=\h+\sum_{\al\in R}\z\cap\g_\al.
$$

\end{theo}

\noindent Assume first that $\h \subset\z$.  Item 1 of Theorem \ref{t10.1} implies that $[\ad x,W]=0$ for any $x\in\h$. Let $R=\{\al_1,\ldots,\al_k\}$ (arbitrary ordering of nonzero roots). Choose a basis in $\g$ of the form $e_1,\ldots,e_l,E_{\al_1},\ldots, \lbr E_{\al_k}, E_{\al_i}\in\g_{\al_i}$, where $e_1,\ldots,e_l$ is any basis in $\h$. The matrix of the endomorphism $\ad x$ in this basis has  the block form
$$
\left[
  \begin{array}{cc}
    0 & 0 \\
    0 & D(x) \\
  \end{array}
\right],
$$
where $D(x):=\diag(\al_1(x),\ldots,\al_k(x))$. For generic $x\in\h$ the numbers $\al_1(x),\ldots,\al_k(x)$ are pairwise distinct, hence the matrix of $W$ in the same basis is of the form
$$
\left[
  \begin{array}{cc}
    M & 0 \\
    0 & D \\
  \end{array}
\right],
$$
where $M$ is an $l\times l$ matrix and $D$ is a diagonal $k\times k$ matrix.

On the other hand,  Item 1 of Theorem \ref{t10.1} implies that $\z^{\theta_i}\subset\ker W^{\theta_i}$, i.e. $\z=\z^{\theta_1}\oplus \cdots\oplus\z^{\theta_m}$ is a sum of eigenspaces of $W$ containing a $W$-invariant subspace $\h$. Hence $\h$ also is a direct sum of eigenspaces of $W$.

Vice versa, assume condition 1 holds. Then, obviously, $[\ad x,W]=0$ for any $x\in\h$ and $\h$ is a direct sum of eigenspaces of $W$, hence by Theorem \ref{t10.1} $x\in\z$.

The proof of the last statement also follows from this theorem.
\qed

\medskip

\abz\label{110.15}
\begin{defi}\rm A semisimple bi-Lie structure $(\g,[,],[,]')$ with the central subalgebra $\z$ will be called {\em regular }if there exists a Cartan subalgebra $\h$ of the semisimple Lie algebra $(\g,[,])$ such that $\h \subset\z$. (This terminology is motivated by the fact that the central subalgebra of a regular bi-Lie structure is a regular reductive Lie subalgebra, see Theorem \ref{110.20}.)
\end{defi}

Fix a Cartan subalgebra $\h \subset\g$ and assume that the equivalent conditions of Theorem \ref{110.10} are satisfied for a semisimple bi-Lie structure $(\g,[,],[,]')$ with the principal WNO $W$.

Choose $E_\al\in\g_\al$ for any $\al\in R$ such that $B_\g(E_\al,E_{-\al})=1$, in particular, $[E_\al,E_{-\al}]=H_\al$, where $B_\g(H_\al,H)=\al(H)$ for any $H\in\h$. Let $\la_\al\in\C$ be such that $W(E_\al)=\la_\al E_\al$. We conclude from Corollary \ref{100.25} that for any $\al\in R$ there exists $\pi_\al\in\C$ such that $\pi_\al=\pi_{-\al}$ and $P(E_\al)=\pi_\al E_\al$, here $P$ is the principal primitive of the principal WNO $W$. Moreover, since $[\ad x,W]=0$ for any $x\in\z$, the definition of $P$ implies that $P|_\z=0$, in particular $P|_\h=0$.

If $W^*$ stands for the operator adjoint to $W$ with respect to the Killing form $B_\g$, then  $W^*E_\al=\la_{-\al}E_\al$ for any $\al\in R$ by Item 4 of Theorem \ref{100.20}. Put $\sigma_\al:=(1/2)(\la_\al+\la_{-\al}),\kappa_\al:=(1/2)(\la_\al-\la_{-\al})$.

\abz\label{110.20}
\begin{theo} Retain the assumptions and notations mentioned. Then
\begin{enumerate}\item
For any $\al\in R$  the set of times $\{t\in T|E_\al\in \ker B^t\}$ (see Definition \ref{defiTimes})  coincides  with the set $T_\al=\{t_{1,\al},t_{2,\al}\}$ of the solutions of the quadratic equation
$$
(t-\la_\al)(t-\la_{-\al})-2\pi_\al=0.
$$
Moreover, $T_\al=T_{-\al}$.

\item  For any $\al\in R$   $$
    \sigma_\al=(1/2)(t_{1,\al}+t_{2,\al}),\kappa_\al\in\sqrt{\zeta_\al-2\pi_\al},
    $$
    where $\zeta_\al:=((t_{1,\al}-t_{2,\al})/2)^2$ and $\sqrt{a}$ stands for the set of square roots of $a$.
\item For any $\al\in R$ the following formula is valid:
$$
W^{t_{1,\al}}W^{t_{2,\al}}H_\al=0
$$
(recall $W^t:=W-t\Id_\g$).
Consequently, $H_\al$ is either an eigenvector of $W$ corresponding to the eigenvalue $t_{1,\al}$, or an eigenvector of $W$ corresponding to the eigenvalue $t_{2,\al}$, or a sum of such eigenvectors.
\item The centre $\z^{\theta_{i}}$ (see Section \ref{90})   of the exceptional bracket $[,]^{\theta_i}$,  is a regular with respect to $\h$ (i.e. stabilized by $\ad\h$) subalgebra reductive in $(\g,[,])$, in particular,  $E_\al\in\z^{\theta_{i}}\Longrightarrow E_{-\al}\in\z^{\theta_{i}}$. Moreover,   $\z^{\theta_{i}}=\h\cap\ker W^{\theta_i}+\bigoplus_{\al\in R^{\theta_i}}\g_\al$, where $R^{\theta_i} \subset R$ is the set of roots $\al$ such that  $\la_\al=\la_{-\al}=t_{1,\al}=t_{2,\al}=\theta_i$ and $\al(\z^{\theta_{j}}\cap\h)=0$ for any $\theta_j\in\Theta,\theta_j\not=\theta_i$.

   In particular, the central subalgebra $\z=\z^{\theta_{1}}\oplus \cdots \oplus\z^{\theta_{m}}$ is a regular with respect to $\h$  reductive subalgebra in $(\g,[,])$.
    \item The set $\Theta=\{\theta_1,\ldots,\theta_m\}$ of the times corresponding to the nontrivial centres is equal to the spectrum of the operator $W|_\h$.
\item The operator $W$ preserves the irreducible components of the representation $x\mapsto\ad_\g x$ of the Lie algebra $\z$ in $\z^\bot \subset\g$ and the  restriction of $W$ to any of them is a scalar operator (here $\z^\bot$ is the orthogonal complement in $\g$ to $\z$ with respect to the Killing form).
\item The set $T$ of times of the bi-Lie structure is exhausted by the numbers $\theta_{i}$ and $t_{j,\al}$, i.e. $T=\Theta\cup\bigcup_{\al\in R}T_\al$.
\end{enumerate}

\end{theo}

\noindent {\em Ad  1.} Indeed, since $WE_\al=\la_\al E_\al,W^*E_\al=\la_{-\al}E_\al,PE_\al=\pi_\al E_\al$, then $\hat B^tE_\al=[\la_{-\al}\la_\al-2\pi_\al-t(\la_\al+\la_{-\al})+t^2]E_\al=
[(t-\la_\al)(t-\la_{-\al})-2\pi_\al]E_\al$ (see formula (\ref{killf})).

\medskip

\noindent {\em Ad  2.} The Vi\`{e}te formulae and Item 1 give
$$
t_{1,\al}+t_{2,\al}=\la_\al+\la_{-\al},\ \ t_{1,\al}t_{2\al}=\la_\al\la_{-\al}-2\pi_\al,
$$
which implies the formula for $\sigma_\al$ and the formula $(\la_\al-\la_{-\al})^2=(t_{1,\al}-t_{2\al})^2-8\pi_\al$, hence the formula for $\kappa_\al$.

\medskip

\noindent {\em Ad  3.} Recall the basic equality $T_W(\cdot,\cdot)=[\cdot,\cdot]_P$ relating the WNO $W$ to its primitive $P$. Calculate this expression on the pair $E_\al,E_{-\al}$:
$$
\la_\al\la_{-\al}H_\al-W(\la_\al+\la_{-\al}-W)H_\al=2\pi_\al H_\al
$$
(we used the fact that $PH_\al=0$). Taking into account the Vi\`{e}te formulae above we get
$W(t_{1,\al}+t_{2\al}-W)H_\al=t_{1,\al}t_{2\al}H_\al$, or finally $(W-t_{1,\al})(W-t_{2\al})H_\al=0$.
\medskip

\noindent {\em Ad  4.} Recall that by Theorem \ref{t10.1} $E_\al\in\z$ if and only if $[\ad E_\al,W]=0$. To show that $\z$ is reductive in $\g$ in view of \cite[Prop. 2, \S 3, Ch. VIII]{bourb7-8} it is sufficient to prove that $E_{-\al}\in\z$. Since $WE_{-\al}=\la_{-\al}E_{-\al}$, this last  is equivalent to $[\ad E_{-\al},W]=0$ again by Theorem \ref{t10.1}.

The equality $[\ad E_\al,W]=0$ is equivalent to the following list of conditions:
\begin{enumerate}
\item[a)] $0=[\ad E_\al,W]E_\be=N_{\al,\be}(\la_{\be}-\la_{\al+\be})E_{\al+\be}$ for any $\be\in R$ such that $\al+\be\in R$;
\item[b)] $0=[\ad E_\al,W]E_{-\al}=(\la_{-\al}-W)H_{\al}=-W^{\la_{-\al}}H_{\al}$;
\item[c)] $0=[\ad E_\al,W]H=(\la_\al\al(H)-\al(WH))E_\al=-\al(W^{\la_\al} H)E_\al$ for any $H\in\h$;
\end{enumerate}
here $N_{\al,\be}\not=0$ is given by $[E_\al,E_\be]=N_{\al,\be}E_{\al+\be}$.

On the other hand, the equality $[\ad E_{-\al},W]=0$ is equivalent to a similar list of conditions:
\begin{enumerate}
\item[a')] $0=[\ad E_{-\al},W]E_{\be'}=N_{-\al,\be'}(\la_{\be'}-\la_{-\al+\be'})E_{-\al+\be'}$ for any $\be'\in R$ such that $-\al+\be'\in R$;
\item[b')] $0=[\ad E_{-\al},W]E_{\al}=W^{\la_\al}H_{\al}$;
\item[c')] $0=[\ad E_{-\al},W]H=\al(W^{\la_{-\al}}H)E_\al$ for any $H\in\k$;
\end{enumerate}
It is clear that a') follows from a) (by putting $\be:=-\al+\be'$). Now we will prove that b), c) imply b'), c').

Indeed, by b) we have $WH_\al=\la_{-\al}H_\al$, whence $W^{\la_\al}H_\al=(\la_{-\al}-\la_\al)H_\al$, from which by c) we deduce that $\la_\al=\la_{-\al}$ (since $\al(H_\al)\not=0$). Thus b), c) give $W^{\la_\al}=W^{\la_{-\al}}$ and b'), c') follow and we have proven that $\z$ is reductive in $\g$.

Now let $E_\al\in \z^{\theta_{i}}$, i.e. $E_\al\in\z\cap\ker W^{\theta_i}$, in particular, $\theta_i=\la_\al$. We have shown that  $\la_\al=\la_{-\al}$ and $E_\al,E_{-\al},H_\al\in \ker W^{\la_\al}$.
Since $\pi_\al=0$, Item 1 implies also that $T_\al=\{\la_\al\}$. Condition $\al(\z^{\theta_{j}}\cap\h)=0,\theta_j\not=\theta_i$, follows from c).

\medskip

\noindent {\em Ad  5.} The fact that the spectrum of $W|_\h$ lies in $\Theta$ follows from Item 1 of Theorem \ref{t10.1} and the inclusion $\h \subset\{x\mid [\ad x,W]=0\}$. Conversely, any $\z^{\theta_i}$ is a regular with respect to $\h$ reductive subalgebra, hence $\z^{\theta_i}\cap\h=\ker W^{\theta_i}\cap\h\not=\{0\}$.

\noindent {\em Ad  6.} Note that any irreducible component of the representation $\ad_\g|_\z$ in $\z^\bot$ is spanned by some root spaces $\g_\al$ hence is preserved by $W$. Moreover, due to the equality $[\ad_\g x, W]=0, x\in\z$, (see Theorem \ref{t10.1}) the operator $W$ is intertwining for the $\ad_\g$ representation of $\z$. Now the result follows by the Schur lemma (another argument  is Item 4 of Theorem \ref{t10.1}, saying that the eigenspaces of $W$ are invariant for $\z^{\theta_i}$).

\medskip

\noindent {\em Ad  7.} Taking into account the block structure of the operators $W,P$ and the fact that $P|_\h=0$ we have $\det\hat B^t=\det((W^*)^tW^t-2P)=\det(W|_\h-t\Id_\h)^2\prod_{\al\in R}((t-\la_\al)(t-\la_{-\al})-2\pi_\al)$. It remains to use Item 5.
\qed

\medskip

Identify for a moment $\h$ with $\h^*$ by means of  the Killing form.  Then each vector $H_\al$ will be identified with the root $\al$ itself and the root system $R$ induces a reduced root system on $\h$ which will be denoted also by $R$.
Item 3 of the theorem above shows that $W|_{\h}$ is an $R$-admissible operator in the following sense.

\abz\label{110.255}
\begin{defi}\rm Let $V$ be a vector space over $\R$ and let $R \subset V$ be a reduced root system in $V$. A diagonalizable linear operator $U:V^\C\to V^\C$ (here $V^\C$ stands for the complexification of the vector space $V$) will be called $R$-{\em admissible} if for any $\al\in R \subset V^\C$ either
\begin{enumerate}\item
there exist two eigenvectors $w_{1,\al},w_{2,\al}\in V^\C$ corresponding to different eigenvalues $t_{1,\al},t_{2,\al}$ of the operator $U$ such that
$$
\al=w_{1,\al}+w_{2,\al};
$$
or
\item $\al$ itself is an eigenvector of $U$ corresponding to an eigenvalue $t_\al$.
\end{enumerate}

We say that a root $\al\in R$ is {\em complete} with respect to an $R$-admissible operator $U$ if condition 1 holds.
An $R$-admissible operator $U$ is called {\em complete} if  any $\al\in R$ is complete with respect to $U$.
\end{defi}

\abz\label{remaa}
\begin{rema}\rm Note that, since the eigenvectors corresponding to different eigenvalues are linearly independent, the eigenvectors $w_{1,\al},w_{2,\al}$ if exist are defined uniquely. Put $V_\al:=\langle w_{1,\al},w_{2,\al}\rangle, U_\al:=\{t_{1,\al},t_{2,\al}\}$  if $\al$ is complete and $V_\al:=\langle \al\rangle,U_\al:=\{t_\al\}$ otherwise (here $\langle  \cdot  \rangle$ stands for the linear span).
\end{rema}

 We see that in order to classify semisimple bi-Lie structures $(\g,[,],[,]')$ with condition $\z \supset \h$  one needs, in particular, to classify $R$-admissible operators on $\h$, where $R=R(\g,\h)$ is the root system of the semisimple Lie algebra $(\g,[,])$ transported to $\h$ by means of the Killing form, up to conjugations by the restrictions to $\h$ of the automorphisms of $(\g,[,])$.

 Another piece of data which can be extracted from the principal WNO of a semisimple bi-Lie structure with the help of Theorem \ref{110.20} and which will be important in the classification matters is a prescription of the pair of times $T_\al$ to any root $\al\in R$. This will be formalized in the following definition (the "times selection rules" will be justified below in Theorem \ref{110.30}).

 \abz\label{110.50}
\begin{defi}\rm Let $R$ be a reduced root system on a vector space $V$. A collection $\{T_\al\}_{\al\in R}$ of   unordered  pairs $T_\al=\{t_{1,\al}t_{2,\al}\}$ of complex numbers is called a {\em  diagram of pairs of times} (or simply a {\em pairs diagram}) if
$T_\al=T_{-\al}$ for any $\al\in R$ and
for any triple $\al,\be,\ga\in R$  such that $\al+\be+\ga=0$ (such triples will be called {\em triangles}) the pairs $T_\al,T_\be,T_\ga$ obey the following "times selection rules":
\begin{enumerate}\item either there exist $t_1,t_2,t_3\in\C$ such that
$$
T_\al=\{t_1t_2\},T_\be=\{t_2t_3\},T_\ga=\{t_3t_1\};
$$
\item or there exist $t_1,t_2\in\C, t_1\not=t_2$, such that
$$
T_\al=T_\be=T_\ga=\{t_1t_2\}.
$$
\end{enumerate}

We write $(T_\al)$ if we understand $T_\al$ as a set. A pair $(U,\T)$, where $U:V^\C\to V^\C$  is an $R$-admissible operator and $\T:=\{T_\al\}_{\al\in R}$  is a pairs diagram such that $U_\al \subset (T_\al)$ (see Remark \ref{remaa}) for any $\al\in R$,  will be called {\em admissible} too.

We put $T_\T:=\bigcup_{\al\in R}(T_\al)$ and call $T_\T$ the {\em set of times} of the diagram $\T$.
\end{defi}

Let  $(U,\T)$  be an admissible pair.
Note that,  if $\al \in R$ is complete,  the set of eigenvalues $U_\al$ coincides with the set of times $(T_\al)$. In  general,  a root $\al$ itself can be an eigenvector of $U$ and the corresponding eigenvalue $t_\al$ is one of the elements of $(T_\al)$.  The second element  is called {\em virtual} and can not be read from $U$ (cf. Examples \ref{140.90}, \ref{140.180}).

The next theorem elucidates further properties of regular semisimple bi-Lie structures, in particular, shows that the set of the pairs $T_\al$ from Theorem \ref{110.20} forms a pairs diagram. We will refer to the "times selection rules" of Definition \ref{110.50} as to TSR 1 and TSR 2.

\abz\label{110.30}
\begin{theo}
Retain the assumptions and notations of Theorem \ref{110.20}.  Then
\begin{enumerate}\item Given a triangle $\al,\be,\ga\in R$,
the sets $T_\al,T_\be,T_\ga$ from Item 1 of Theorem \ref{110.20} understood as unordered pairs obey the "times selection rules" of Definition \ref{110.50}.
\item If the triangle is such that  TSR 1 (respectively  TSR 2)  occurs,
the following equality holds:
$$
\kappa_\al+\kappa_\be+\kappa_\ga=0;
$$
respectively:
$$
\kappa_\al+\kappa_\be+\kappa_\ga=\pm(t_1-t_2)/2.
$$
\item The set $\g^{t}_0:=\h\cap\ker W^t+\bigoplus_{\al\in R^t_0}\g_\al$, where $R^t_0 \subset R$ is the set of roots $\al$ with $t_{1,\al}=t_{2,\al}=t$, forms a regular with respect to $\h$ subalgebra reductive in $(\g,[,])$ such that $\g^{t}_0\supset \z^t$, where $\z^t$ is the centre of the corresponding bracket $[,]^{t}$. The subalgebra $\g_0^t$ is nontrivial if and only if so is the subalgebra $\z^t$.

\item (In Items 4--7 assume additionally that $R$ is irreducible, i.e. that $(\g,[,])$ is simple.) For any $\al\in R^t_0$ we have $\ka_\al=0$ and $\pi_\al=0$.
\item For any $t\in \C$ we have $\g_0^t \subset\ker W^t$. In particular, $\g_0^{t_1}\cap\g_0^{t_2}=\{0\}$ for $t_1\not=t_2$.
\item The set $\g_0=\g_0^{\theta_1}\oplus \cdots\oplus\g_0^{\theta_m}$ (here $\{\theta_1,\ldots,\theta_m\}=\Theta$ is the set of times for which the centre $\z^{\theta_i}$ of the corresponding bracket $[,]^{\theta_i}$ is nontrivial, see Theorems \ref{t10.1}, \ref{110.20}) is a subalgebra reductive in $(\g,[,])$ such that $\g_0\supset \z$. The subspace $\g_0^{\theta_i}\oplus\g_0^{\theta_j}$ is a subalgebra in $(\g,[,])$ for any $i,j$. The subalgebra $\g_0$ will be called the {\em basic subalgebra} of a regular semisimple bi-Lie structure.
\item The antisymmetric part $W_{\mathrm{a}}$ of the operator $W|_{\h^\bot}$ preserves the irreducible components of the representation $x\mapsto\ad_\g x$ of the Lie algebra $\g_0$ in $\g_0^\bot \subset\g$ and the  restriction of $W_{\mathrm{a}}$ to any of them is a scalar operator, here the orthogonal complement  is taken with respect to the Killing form.
\end{enumerate}
\end{theo}

\abz\label{remaover}
\begin{rema}\rm
It can be shown that $\g_0 \subset\hat\z$, where $\hat\z$ is the subalgebra from Theorem \ref{overcentral}. We conjecture that in fact $\g_0=\hat\z$.
\end{rema}

\abz\label{rema23}
\begin{rema}\rm
The assumption that $(\g,[,])$ is simple in  Items 4--7 reduces complexity of formulations. It can be substituted by a less restrictive assumption that $(\g,[,])$ is semisimple and the subalgebras $\g_0^t$ do not contain any of its simple components.
\end{rema}

The proof of the theorem will use the following lemma.

\abz\label{110.40}
\begin{lemm}Retain the assumptions of Item 1 of the theorem. Let $T_\al:=\{t_1t_2\},T_\be:=\{t_3t_4\}, \lbr T_\ga:=\{t_5t_6\}$. Then the following equalities hold:
\begin{eqnarray*}
(t_5+t_6)(t_1+t_2-t_3-t_4)+(t_3^2+t_4^2-t_1^2-t_2^2)&=&0\\
(t_1+t_2)(t_3+t_4-t_5-t_6)+(t_5^2+t_6^2-t_3^2-t_4^2)&=&0\\
(t_3+t_4)(t_5+t_6-t_1-t_2)+(t_1^2+t_2^2-t_5^2-t_6^2)&=&0.
\end{eqnarray*}

\end{lemm}

\noindent The proof will be deduced from the basic equality $T_W(\cdot,\cdot)=[\cdot,\cdot]_P$. Taking $E_\al,E_\be$ as the arguments, we get $\la_\al\la_\be-\la_{-\ga}(\la_\al+\la_\be-\la_{-\ga})=\pi_\al+\pi_\be-\pi_{-\ga}$, hence $(\la_{-\ga}-\la_\al)(\la_{-\ga}-\la_\be)=\pi_\al+\pi_\be-\pi_{-\ga}$.  Substituting here $-\al,-\be,-\ga$ instead of $\al,\be,\ga$ respectively and recalling that
$\pi_\al=\pi_{-\al}$ for any $\al\in R$ we obtain
$$
(\la_{-\ga}-\la_\al)(\la_{-\ga}-\la_\be)=(\la_{\ga}-\la_{-\al})(\la_{\ga}-\la_{-\be})
$$
Recalling that $\la_{\pm\al}=\sigma_\al\pm\kappa_\al$ (and $\sigma_{-\al}=\sigma_\al,\kappa_{-\al}=-\kappa_\al$) we get
$$
((\sigma_\ga-\sigma_\al)-(\kappa_\ga+\kappa_\al))
((\sigma_\ga-\sigma_\be)-
(\kappa_\ga+\kappa_\be))=
((\sigma_\ga-\sigma_\al)+(\kappa_\ga+\kappa_\al))
((\sigma_\ga-\sigma_\be)+
(\kappa_\ga+\kappa_\be)).
$$
Put $A:=\sigma_\ga-\sigma_\al,B:=\kappa_\ga+\kappa_\al,
C:=\sigma_\ga-\sigma_\be, D:=\kappa_\ga+\kappa_\be$. Then the above equality can be rewritten as $(A-B)(C-D)=(A+B)(C+D)$, whence $AD+BC=0$ and $(A-B)(C-D)=AC+BD$. Thus
\begin{eqnarray}\equ\label{eqsh}
(\sigma_\ga-\sigma_\al)(\kappa_\ga+\kappa_\be)+
(\kappa_\ga+\kappa_\al)(\sigma_\ga-\sigma_\be)&=&0,\\
\equ\label{eqlong}
(\sigma_\ga-\sigma_\al)(\sigma_\ga-\sigma_\be)+
(\kappa_\ga+\kappa_\al)(\kappa_\ga+\kappa_\be)&=&
\pi_\al+\pi_\be-\pi_\ga.
\end{eqnarray}
Now we note that the formula from Item 2 of Theorem \ref{110.20} implies the formula $\pi_\al=(1/2)(\zeta_\al-\kappa_\al^2)$ in account of  which  equality (\ref{eqlong}) gives
$$
(\sigma_\ga-\sigma_\al)(\sigma_\ga-\sigma_\be)+
(\kappa_\ga+\kappa_\al)(\kappa_\ga+\kappa_\be)=\frac{1}{2 }(\zeta_\al+\zeta_\be-\zeta_\ga
-(\kappa_\al^2+\kappa_\be^2-\kappa_\ga^2)),
$$
or, equivalently,
\begin{equation}\equ\label{eqI}
(\sigma_\ga-\sigma_\al)(\sigma_\ga-\sigma_\be)+\frac{1}{2}
(\kappa_\al+\kappa_\be+\kappa_\ga)^2=\frac{1}{2 }(\zeta_\al+\zeta_\be-\zeta_\ga).
\end{equation}
Cyclic permutations of $\al,\be,\ga$ give
\begin{eqnarray}
\equ\label{eqII}(\sigma_\be-\sigma_{\ga})(\sigma_\be-\sigma_\al)+ \frac{1}{2 }\kappa=\frac{1}{2 }(\zeta_{\ga}+\zeta_\al-\zeta_\be)\\
\equ\label{eqIII}(\sigma_\al-\sigma_\be)(\sigma_\al-\sigma_{\ga})+ \frac{1}{2 }\kappa=\frac{1}{2 }(\zeta_{\be}+\zeta_\ga-\zeta_\al),
\end{eqnarray}
where we put $\kappa:=(\kappa_\al+\kappa_\be+\kappa_\ga)^2$.

Adding respectively equalities (\ref{eqI}) and (\ref{eqII}), (\ref{eqI}) and (\ref{eqIII}), (\ref{eqII}) and (\ref{eqIII}) we get the following equalities:
\begin{eqnarray}\equ
(\sigma_{\be}-\sigma_{\ga})^2+\kappa&=&\zeta_{\al}\\ \nonumber
(\sigma_{\ga}-\sigma_{\al})^2+\kappa&=&\zeta_{\be}\\ \nonumber
(\sigma_\al-\sigma_{\be})^2+\kappa&=&\zeta_{\ga}.\label{kappa}
\end{eqnarray}
 From this we get $(\sigma_{\be}-\sigma_{\ga})^2-\zeta_\al=(\sigma_{\ga}-\sigma_{\al})^2-\zeta_\be$, or, equivalently (using formulae from Item 2 of Theorem \ref{110.20}),
$$
(t_5+t_6-(t_3+t_4))^2-(t_1-t_2)^2=(t_5+t_6-(t_1+t_2))^2-(t_3-t_4)^2.
$$
Elementary calculations show that this is equivalent to the first equality of the lemma. The remaining equalities are proved analogously. \qed

\medskip

\noindent{\em Proof of Theorem \ref{110.30}. Ad 1.} We have to consider several cases which exhaust all possible ones. Recall that $W|_{\h}$ is $R$-admissible in the sense of Definition \ref{110.255} and use notation from Remark \ref{remaa}. Recall that we identify $\h$ with $\h^*$ (in particular,  $H_\al$ with $\al$) by means of the Killing form $B_\g$.

\medskip

\noindent{\em Case a)}: $\al,\be,\ga$ are complete. Assume first that $V_\al\cap V_\be=\{0\}$. Then, necessarily, $T_\al=T_\be=\{t_1t_2\}$, where $t_1\not=t_2$ (otherwise the vectors $w_{1,\al},w_{2,\al},w_{1,\be},w_{2,\be}$ would be linearly independent and $\ga$ would depend on more than two eigenvectors), and $T_\ga=\{t_1t_2\}$ (TSR 2). The cases $V_\al\cap V_\ga=\{0\},V_\ga\cap V_\be=\{0\}$ are analogous.

Now assume that $V_\al\cap V_\be\not=\{0\},V_\al\cap V_\ga\not=\{0\},V_\ga\cap V_\be\not=\{0\}$. Then either $V_\al=V_\be=V_\ga=\langle w_1,w_2\rangle$, where $w_i$ is an eigenvector corresponding to an
eigenvalue $t_i,i=1,2$, and $T_\al=T_\be=T_\ga=\{t_1t_2\}, t_1\not=t_2$ (TSR 2), or $\al=w_1-w_2,\be=w_2-w_3,\ga=w_3-w_1$, where $w_i$ is an eigenvector corresponding to an
eigenvalue $t_i,i=1,2,3$, and $T_\al=\{t_1t_2\},T_\be=\{t_2t_3\},T_\ga=\{t_3t_1\}$ with $t_1\not=t_2,t_2\not=t_3,t_3\not=t_1$ (TSR 1).

\noindent{\em Case b)}: one of the roots, say $\ga$, is not complete, but the other two are complete. Then $V_\al=V_\be=\langle  w_1,w_2  \rangle$, where $w_i$ is an eigenvector corresponding to an
eigenvalue $t_i,i=1,2, t_1\not=t_2$, and $\al=w_1+w_2,\be=-w_1+aw_2$ for some $a\not=-1$, i.e. $V_\ga=\langle  w_2  \rangle, T_\al=T_\be=\{t_1t_2\}, T_\ga=\{t_2(t_6)\}$ for some $t_6\in\C$ (here
the brackets $()$ mean that the corresponding element of the unordered pair is virtual, see discussion after Definition \ref{110.50}). A priori $t_6$ can be arbitrary, but in fact  is not so. Indeed, the substitution $t_3=t_1,t_4=t_2,t_5=t_2$ reduces the equalities of  Lemma \ref{110.40}  to one equality $(t_6-t_1)(t_2-t_6)=0$. Hence $t_6=t_1$ (TSR 2) or $t_6=t_2$ (TSR 1).

\noindent{\em Case c)}: two of the roots, say $\al,\be$, are not complete, but the remaining one is compete. Then $V_\al=\langle  \al  \rangle,V_\be=\langle  \be  \rangle, V_\ga=\langle  \al,\be  \rangle$, the roots $\al,\be$ are the eigenvectors corresponding to eigenvalues $t_1,t_3$ respectively with $t_1\not=t_3$. Moreover $T_\al=\{t_1(t_2)\},T_\be=\{t_3(t_4)\},T_\ga=\{t_1t_3\}$ for some $t_2,t_4\in \C$.
The second and third of the equalities from the lemma give $(t_4-t_1)(t_2-t_4)=0$ and $(t_3-t_2)(t_4-t_2)=0$ respectively.
Thus either $t_2=t_4$ (TSR 1), or $t_2\not=t_4$, but $t_4=t_1$ and $t_3=t_2$ (TSR 2).

\noindent{\em Case d)}: all the three roots $\al,\be,\ga$ are not complete. Then $V_\al=\langle  \al  \rangle,V_\be=\langle  \be  \rangle, V_\ga=\langle  \ga  \rangle$ and  the roots $\al,\be,\ga$ are the eigenvectors corresponding to same eigenvalue $t_1$. Then
$T_\al=\{t_1(t_2)\},T_\be=\{t_1(t_4)\},T_\ga=\{t_1(t_6)\}$. The first two equalities of the lemma give:
\begin{eqnarray*}
(t_4-t_2)(t_4+t_2-t_1-t_6)=0\\
(t_6-t_4)(t_6+t_4-t_1-t_2)=0.
\end{eqnarray*}
Hence, either  $t_4=t_2$, which implies $(t_6-t_2)(t_6-t_1)=0$ and $t_6$ equal to $t_1$ (TSR 1) or $t_2$ (TSR 2 or 1 with all times equal), or  $t_4=-t_2+t_1+t_6$, which implies $(t_6-t_4)(t_6-t_2)=0$. In the last case we have either $t_6=t_4$ which implies $t_1=t_2$ (TSR 1), or $t_6=t_2$ which implies $t_4=t_1$ (TSR 1).

This finishes the list of possible cases and the proof of the first item  of the theorem.

\medskip

\noindent{\em  Ad 2.} Item 2 is proved by direct check  using formulae (\ref{kappa}) and the definition of $\sigma_\al, \zeta_\al$ (given before Theorem \ref{110.20}).

\medskip

\noindent{\em  Ad 3.} The fact that $R^t_0$ is a closed set of roots follows from the times selection rules (TSR 1). The symmetry of $R^t_0$ (which will imply the reductivity), is due to the property $T_\al=T_{-\al}$. The fact that $\g^t_0$ is a subalgebra follows from the condition $H_\al=[E_\al,E_{-\al}]\in\ker W^t,\al\in R_0^t$, which is due to Item 3 of Theorem \ref{110.20} and diagonalizability of $W$.

The inclusion $\g^t_0\supset\z^t$ follows from Item 4 of Theorem \ref{110.20}. Thus the nontriviality of $\z^t$ impies that of $\g_0^t$. Finally, if $\g_0^t$ is nontrivial, then $\h\cap\ker W^t\not=\{0\}$ and $\z^t$ is nontrivial by Item 5 of Theorem \ref{110.20}.

\medskip

\noindent{\em  Ad 4.} Assume $R_0^t\not=\emptyset$. We will first prove that
\begin{equation}\equ\label{ass1}
 \forall \al \in R^t_0\ \exists\be\in R\setminus R^t_0:\al+\be\in R.
\end{equation}
Let $R_\mathrm{max}$ be the maximal symmetric closed proper root set containing $R^t_0$. We will show that
\begin{equation}\equ\label{ass2}
\forall\al \in R_\mathrm{max}\ \exists\be,\ga\in R\setminus R_\mathrm{max}:\al=\be+\ga,
\end{equation} which will prove assertion (\ref{ass1}).

All maximal symmetric closed  root sets, which are the root sets of maximal reductive regular subalgebras of simple Lie algebras, are known \cite[Ch. 6]{vinbergOnishIII}. It appears that each of these subalge- \linebreak bras is a fixed point subalgebra $\g_\mathrm{fix}$ of an inner automorphism of order 2, 3 or 5. Obviously, the assertion (\ref{ass2}) is equivalent to the following one: the subspace $[\m,\m]\cap\g_\mathrm{fix}$ coincides with $\g_\mathrm{fix}$, here $\m=(\g_\mathrm{fix})^\bot$ is the orthogonal complement to $\g_\mathrm{fix}$ with respect to the Killing form.

The last assertion follows from the fact that $[\m,\m]+\m$ is an ideal in the initial simple Lie algebra, which can be proved directly using the commutation relations of the corresponding $\Z_2,\Z_3$ or $\Z_5$ gradings.

Now let $\al \in R^t_0,\be\in R\setminus R^t_0$ be such that $-\ga:=\al+\be\in R$. Then by the times selection rules there exists $t'\in \C,t'\not= t$, such that  $T_\be=T_\ga=\{tt'\}$ (recall $T_\al=\{tt\}$). By formula (\ref{eqsh}) we have $((t'-t)/2)(\ka_\ga+\ka_\be)=0$, whence $\ka_\ga=-\ka_\be$. Using  Item 2, we conclude that $\ka_\al=0$.

The equality $\pi_\al=0$ now follows from the formula $\ka_\al^2=\zeta_\al-2\pi_\al$ (see Item 2 of Theorem \ref{110.20}).

\medskip

\noindent{\em  Ad 5.} The previous item and Item 1 of Theorem \ref{110.20} imply that $\la_\al=\la_{-\al}=t$ for any $\al\in R_0^t$, hence $E_\al,E_{-\al}\in\ker W^t$.

\medskip

\noindent{\em  Ad 6.} The proof of the fact that $\g_0^{\theta_i}\oplus\g_0^{\theta_j}$ is a subalgebra in $(\g,[,])$ follows from the definition of $\g_0^t$. The rest of the proof follows from Items 3 and 5.

\medskip

\noindent{\em  Ad 7.} Any irreducible component of the representation $\ad_\g|_{\g_0}$ in $\g_0^\bot$ is spanned by some root spaces $\g_\al$, hence is preserved by $W_{\mathrm{a}}$. If $E_\al\in\g_0,E_\be\in \g_0^\bot$, we have $[\ad E_\al,W_{\mathrm{a}}]E_\be=N_{\al,\be}(\ka_\be-\ka_{\al+\be})E_{\al+\be}=0$ (cf. the proof of Item 4 of Theorem \ref{110.20}), where the last equality is due to the first equality of Item 2 and the fact that $\ka_\al=0$ (Item 4). Moreover, $[\ad H,W_{\mathrm{a}}]E_\be=\be(H)\ka_\be E_\be-\ka_\be\be(H)E_\be=0$ for any $H\in\h$. Thus the operator $W_{\mathrm{a}}$ is intertwinning for the corresponding representation and the Schur lemma can be applied.
\qed

\medskip

Items 2 and 4 of Theorem \ref{110.30} imply the following result in case when $(\g,[,])$ is simple (cf. Remark \ref{rema23}).

\abz\label{110.501}
\begin{coro} Assume $(\g,[,])$ is simple.
Then the antisymmetric part of the operator $W|_{\h^\bot}$  obeys the $((t_1-t_2)/2)$-triangle rule subject to the root set $R_0 \subset R$ in the sense of the following definition, here $R_0:=R_0^{\theta_1}\cup\cdots\cup R_0^{\theta_m}$ (see Items 3, 6 of Theorem \ref{110.30}).
\end{coro}

\abz\label{130.30}
\begin{defi}\rm  Let $R_0 \subset R$ be a closed symmetric root set (see Section \ref{130}) and  let $S\in\End(\h^\bot)$ be an  antisymmetric  operator preserving the root spaces $\g_\al$ of the corresponding root decomposition, $S|_{\g_\al}=\ka_\al\Id_{\g_\al},\ka_\al\in\C, \ka_\al=-\ka_{-\al}$ (the last equality follows from Item 4 of Theorem \ref{100.20}). Given a constant $a\in\C$, we say that $S$ obeys the {\em $a$-triangle rule subject to the root set $R_0$}, if $\ka_\al=0$ for any $\al\in R_0$ and  for any triangle $\al,\be,\ga\in R$ (see Definition \ref{110.50})
\begin{enumerate}\item $\ka_\al+\ka_\be+\ka_\ga=0$  whenever two of three roots belong to $R\setminus R_0$ and the third one to $R_0$ or all the three roots belong to $R_0$;
\item $\ka_\al+\ka_\be+\ka_\ga=\pm a$ whenever $\al,\be,\ga\in R\setminus R_0$.
\end{enumerate}
  We say that  $0,+$ or $-$ is a {\em label} of the corresponding triangle. The whole family of labels is called {\em the system of labels} corresponding to the operator $S$ and is denoted by $\L_S$.
\end{defi}

In fact due to Items 5--7 of Theorem \ref{110.30} we can say more about the antisymmetric part $W_{\mathrm{a}}$ of the operator $W|_{\h^\bot}$.

\abz\label{110.60}
\begin{defi}\rm
Let a $\Ga$-grading $\g=\bigoplus_{i\in\Ga}\g_i$ of a Lie algebra $\g$, where $\Ga$ is an abelian group, be given. Put $\overline\Ga:=\{i\in\Ga\setminus\{0\}\mid \g_i\not=\{0\}\}$ and call the elements of this set {\em quasiroots}. We call a {\em triangle} any triple $i,j,k\in\overline\Ga$ such that $i+j+k=0$. Given an antisymmetric operator $S\in\End(\g)$ which is scalar on each $\g_i$, i.e. such that $S|_{\g_i}=\ka_i\Id_{\g_i}$ for some $\ka_i\in\C$ (necessarily satisfying equalities $\ka_i=-\ka_{-i}$, cf. Item 4 of Theorem \ref{100.20}), and a constant $a\in\C$, we say that $S$ obeys the {\em $a$-triangle rule subject to the grading}, if
$$
\ka_i+\ka_j+\ka_k=\pm a
$$
 for any triangle $i,j,k$.
  We say that  $+$ or $-$ is a {\em label} of the corresponding triangle. The whole family of labels is called {\em the system of labels} corresponding to the operator $S$ and is denoted by $\L_S$.
\end{defi}

\abz\label{110.65}
\begin{rema}\rm
Note that if an endomorphism $S\in\End(\h^\bot)$ satisfies the conditions of Definition \ref{130.30}, the extended  by zero on $\h$ endomorphism will obey the $a$-triangle rule subject to the root grading in the sense of Definition \ref{110.60}.

On the other hand, let $\g_0 \subset\g$ be a subalgebra reductive in $\g$ of maximal rank, let $\h \subset\g_0$ be some Cartan subalgebra  and let $\g=\bigoplus_{i\in\Ga}\g_i$ be the irreducible toral grading (see Appendix \ref{ApToral}). Then an endomorphism $S\in\End(\g)$ satisfies the conditions of Definition \ref{110.60} if and only if its restriction to $\h^\bot$  obeys the $a$-triangle rule subject to the root set $R_0$ in the sense of Definition \ref{110.50}, where $R_0 \subset R(\g,\h)$ is the root set corresponding to the subalgebra $\g_0$.
\end{rema}

Items 2, 4--7 of Theorem \ref{110.30} imply the following result.

\abz\label{110.70}
\begin{coro} Let $(\g,[,])$ be simple.
The antisymmetric part of the operator $W|_{\h^\bot}$ extended by zero on $\h$ obeys the $((t_1-t_2)/2)$-triangle rule subject to the irreducible toral $\Ga(\g_0)$-grading on $(\g,[,])$ (see Definition \ref{ApToral.20}) corresponding to the reductive Lie subalgebra $\g_0 \subset\g$ (the basic subalgebra).
\end{coro}

\section{Two classes of pairs diagrams}
\label{s120}

The aim of this section is to prove the following theorem.

\abz\label{120.10}
\begin{theo} Let $R$ be a reduced irreducible root system in a vector space $V$ over $\R$ and let $\T:=\{T_\al\}_{\al\in R}$ be a pairs diagram. Assume that there exist $\al,\be,\ga\in R$  such that $\al+\be+\ga=0$ and
$$
T_\al=T_\be=T_\ga=\{t_1t_2\}
$$
for some $t_1,t_2\in \C,t_1\not=t_2$. Then the set of times $T_\T$ is equal to $\{t_1,t_2\}$ (see Definition \ref{110.50}).
\end{theo}

\abz\label{120.15}
\begin{defi}\rm Let $R$ be a reduced root system. We say that a pairs diagram $\{T_\al\}_{\al\in R}$
 is of  {\em Class II}, if there  exist $\al,\be,\ga\in R$ satisfying the hypotheses of the theorem, and of {\em Class I}, if  such roots do not exist.
\end{defi}

We start from auxiliary results.

Let $R$ be a reduced root system. We define on $R$ an equivalence relation by putting $\al\sim\al$ and $\al\sim -\al$ for any $\al\in R$. Write $\widetilde{R}:=R/\sim$. Obviously, any pairs diagram is in fact indexed by the set $\widetilde{R}$.
An unordered triple $\{a,b,c\} \subset \widetilde{R}$ is called a {\em triangle} if there exist $\al\in a,\be\in b,\ga\in c$ such that $\al+\be+\ga=0$. Given a pairs diagram $\{T_a\}_{a\in \widetilde{R}}$, a triangle $\{a,b,c\}$ is said to be $\{t_1t_2\}$-{\em triangle}, if
$$
T_a=T_b=T_c=\{t_1t_2\}
$$
for some $t_1,t_2\in \C,t_1\not=t_2$.

\abz\label{120.20}
\begin{lemm} Let $R$ be a reduced irreducible root system in a vector space $V$ and let $\{T_a\}_{a\in \widetilde{R}}$ be a pairs diagram. Let
$a,b,c,a',b',c'\in \widetilde{R}$ be pairwise distinct elements such that $\{a,b,c'\},\{a,c,b'\},\lbr\{b,c,a'\},\{a',b',c'\}$ are triangles (see the figure below).
Then, if $\{a,b,c'\}$ is a
$\{t_1t_2\}$-triangle, the pairs of times $T_{a'},T_{b'},T_c$ are equal to one of the pairs
$$
\{t_1t_2\},\{t_1t_1\},\{t_2t_2\}
$$
and moreover at least one of the triangles $\{a,c,b'\},\{b,c,a'\},\{a',b',c'\}$ is a $\{t_1t_2\}$-triangle.
\end{lemm}

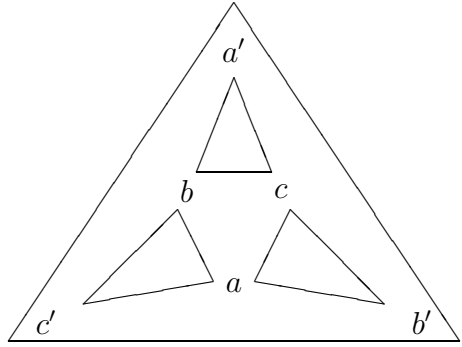
\begin{center}
\setlength{\unitlength}{0,25mm}
\begin{picture}(300,200)
\put(20,20){\line(2,3){120}}
\put(140,200){\line(2,-3){120}}
\put(20,20){\line(1,0){240}}
\put(60,40){\line(1,1){50}}
\put(60,40){\line(6,1){69}}
\put(110,90){\line(1,-2){19}}
\put(170,90){\line(-1,-2){19}}
\put(170,90){\line(1,-1){50}}
\put(220,40){\line(-6,1){69}}
\put(120,110){\line(2,5){20}}
\put(120,110){\line(2,0){40}}
\put(140,160){\line(2,-5){20}}
\put(40,30){\makebox(0,0){$c'$}}
\put(240,30){\makebox(0,0){$b'$}}
\put(140,50){\makebox(0,0){$a$}}
\put(115,100){\makebox(0,0){$b$}}
\put(165,100){\makebox(0,0){$c$}}
\put(140,175){\makebox(0,0){$a'$}}
\end{picture}
\end{center}

\noindent The times selection rules (TSR for short) applied to the triangle $a,c,b'$ imply, that each of the pairs $T_{c},T_{b'}$ should contain at least one of the times $t_1,t_2$. So, a priori, we have the following possibilities:
\begin{enumerate}\item
 $T_a=\{t_1t_2\},T_{c}=\{t_1t_2\},T_{b'}=\{t_1t_2\}$;
\item $T_a=\{t_1t_2\},T_{c}=\{t_2t_1\},T_{b'}=\{t_1t_1\}$;
\item $T_a=\{t_2t_1\},T_{c}=\{t_1t_1\},T_{b'}=\{t_1t_2\}$;
\item $T_a=\{t_2t_1\},T_{c}=\{t_1t_2\},T_{b'}=\{t_2t_2\}$;
\item $T_a=\{t_1t_2\},T_{c}=\{t_2t_2\},T_{b'}=\{t_2t_1\}$;
\item $T_a=\{t_1t_2\},T_{c}=\{t_2t_3\},T_{b'}=\{t_3t_1\}$ for some $t_3\not=t_1,t_2$;
\item $T_a=\{t_2t_1\},T_{c}=\{t_1t_3\},T_{b'}=\{t_3t_2\}$ for some $t_3\not=t_1,t_2$.
 \end{enumerate}
However the last two do not occur. Indeed, assume possibility 6 occurs. Then the TSR for the triangle $b,c,a'$ would imply $T_{a'}=\{t_1t_3\}$ which contradicts to the TSR for the triangle $c',b',a'$. Analogous considerations show impossibility of 7.

Thus we have proved that pairs of times $T_{c},T_{b'}$ are equal to one of the pairs $\{t_1t_2\},\{t_1t_1\},\{t_2t_2\}$.
By the symmetry of the triangles $a,c,b'$ and $b,c,a'$ one comes to the same conclusion about $T_{a'}$.

To finish the proof note that the TSR imply that the pairs $\{t_1t_1\}$ or $\{t_2t_2\}$ can appear only once among $T_c,T_{b'},T_{a'}$ and the rest should be equal to $\{t_1t_2\}$.
\qed

\abz\label{120.30}
\begin{theo} Theorem \ref{120.10} holds for the root system $R:=\a_n$.
\end{theo}

\noindent We will give a "geometrical" proof of this theorem. The elements of the set $\widetilde{R}$ will be represented as points of the $xy$-coordinate plane with integer coordinates $(l,m)$ with $l,m\ge 0,l+m\le n$.  The points $(n,0),(n-1,1),\ldots,(0,n)$ represent the elements $[\al_1],\ldots,[\al_n]$, where $B:=\{\al_1,\ldots,\al_n\}$ is a basis of the root system $\a_n$. Recall \cite[Ch. VI,\S 1,Cor. 3 of Prop. 19]{bourb7-8} that, $B$ is identified with the set of vertices of the graph of corresponding Coxeter system (which  will be denoted by $\Ga(R)$ and called the {\em Coxeter graph} of $R$) and that   for any subset $Y \subset B$ which is connected as a subset in $\Ga(R)$  we have $\sum_{\be\in Y}\be\in R$. On the other hand, for the root system $\a_n$ all the positive roots are obtained this way. From this we can get the following lemma.

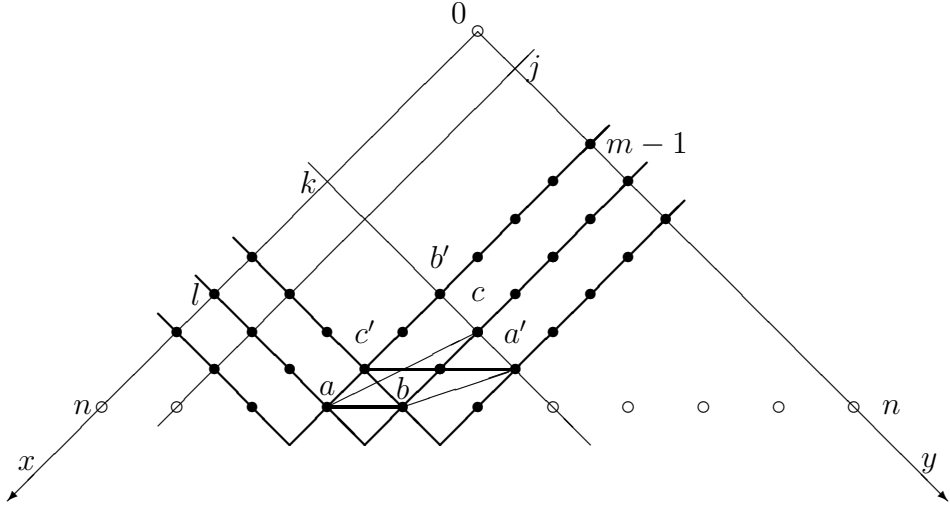
\abz\label{120.40}
\begin{lemm} Let $a=(l_0,m_0)\in \widetilde{\a}_n$. Then the elements $b\in \widetilde{\a}_n$ such that there exist a triangle $a,b,c$ for some $c\in \widetilde{\a}_n$ lie on the "X-shape" $X(a)$, i.e. the union of lines $x=l_0,y=m_0,x=n-m_0+1,y=n-l_0+1$ intersected with $\widetilde{\a}_n$ (see the figure below) with the crossing at $a$.

Given $a\in \widetilde{\a}_n$ and $b\in X(a)$, the element $c\in \widetilde{\a}_n$ such that $a,b,c$ is a triangle is uniquely defined by $c=X(a)\cap X(b)$.
\end{lemm}

\begin{center}
\setlength{\unitlength}{0,25cm}
\begin{picture}(50,28)
\put(25,25){\vector(-1,-1){25}}
\put(25,25){\vector(1,-1){25}}
\put(25,25){\circle*{0,5}}
\put(23,23){\circle*{0,5}}
\put(21,21){\circle*{0,5}}
\put(27,23){\circle*{0,5}}
\put(29,21){\circle*{0,5}}
\put(25,21){\circle*{0,5}}

\put(5,5){\circle*{0,5}}
\put(9,5){\circle*{0,5}}
\put(13,5){\circle*{0,5}}
\put(17,5){\circle*{0,5}}
\put(21,5){\circle*{0,5}}
\put(29,5){\circle*{0,5}}
\put(33,5){\circle*{0,5}}
\put(37,5){\circle*{0,5}}
\put(41,5){\circle*{0,5}}
\put(45,5){\circle*{0,5}}

\put(7,7){\circle*{0,5}}
\put(11,7){\circle*{0,5}}
\put(15,7){\circle*{0,5}}
\put(9,9){\circle*{0,5}}
\put(13,9){\circle*{0,5}}
\put(11,11){\circle*{0,5}}

\put(23,19){\circle*{0,5}}
\put(25,17){\circle*{0,5}}
\put(27,15){\circle*{0,5}}
\put(29,13){\circle*{0,5}}
\put(31,11){\circle*{0,5}}
\put(33,9){\circle*{0,5}}
\put(35,7){\circle*{0,5}}
\put(43,7){\circle*{0,5}}

\put(23,7){\circle*{0,5}}
\put(25,9){\circle*{0,5}}
\put(27,11){\circle*{0,5}}
\put(29,13){\circle*{0,5}}
\put(31,15){\circle*{0,5}}
\put(33,17){\circle*{0,5}}

\thicklines
\put(20,22){\line(1,-1){19}}
\put(39,3){\line(1,1){5}}
\put(10,12){\line(1,-1){9}}
\put(19,3){\line(1,1){15}}

\thinlines
\put(16,18){\line(1,-1){15}}
\put(31,3){\line(1,1){9}}
\put(25,9){\line(5,-1){10}}

\put(1,2){\makebox(0,0){$x$}}
\put(4,5){\makebox(0,0){$n$}}
\put(19,21){\makebox(0,0){$l_0$}}
\put(6,11){\makebox(0,0){$n-m_0+1$}}
\put(24,26){\makebox(0,0){$0$}}
\put(35,17){\makebox(0,0){$m_0$}}
\put(48,7){\makebox(0,0){$n-l_0+1$}}
\put(47,5){\makebox(0,0){$n$}}
\put(49,2){\makebox(0,0){$y$}}
\put(29,14){\makebox(0,0){$a$}}
\put(25,10){\makebox(0,0){$b$}}
\put(35,8){\makebox(0,0){$c$}}
\end{picture}
\end{center}

\abz\label{120.50}
\begin{lemm}Let $R\not=\g_2,\T$ and $\al,\be,\ga\in R$ satisfy hypotheses of Theorem \ref{120.10}. Then there exists a basis $\al_1,\ldots,\al_n$ of the root system $R$ and two its elements $\al_i,\al_j$ such that $[\al_i],[\al_j]$ belong to the  triangle $\{[\al_i],[\al_j],[\al_i+\al_j]\}$ (in particular, $\al_i,\al_j$ form a connected subgraph of the Coxeter graph of $R$) and $T_{[\al_i]}=T_{[\al_j]}=T_{[\al_i+\al_j]}=\{t_1t_2\}$.
\end{lemm}

\noindent Let $V' \subset V$ be a (2-dimensional) subspace generated by $\al,\be,\ga$ an let $R':=R\cap V'$. Then $R'$ is a root system in $V'$ by \cite[Cor. of Prop. 4, \S 1, Ch. VI]{bourb7-8} and is equal to one of the root systems $\a_2,\b_2$. Indeed, $R'$ is an irreducible reduced root system of rank 2; the list of all such  systems is $\a_2,\b_2,\g_2$ (cf. \cite[Ch. X, Exc. B1]{helgason}). However the case $\g_2$ is excluded by the assumption $R\not=\g_2$ as we will show below.

Let $B'$ be a basis of $R'$. Then by \cite[Prop. 24, \S 1, Ch. VI]{bourb7-8} there exists a basis $B$ of $R$ such that $B' \subset B$. In particular, this shows that $R'\not=\g_2$, since neither of the Coxeter graphs of reduced irreducible root systems $R$ contains the graph of $\g_2$ except $R=\g_2$ itself.

Now, the direct inspection of the systems $\a_2,\b_2$ shows that among the roots $\al,\be,\ga,-\al,-\be,-\ga$ there exist two roots $\al',\be'$ forming a basis $B'$ of $R'$ and the proposition cited gives the result. \qed

\medskip

Now we are able to prove Theorem \ref{120.30}. Let $[\al_i],[\al_j]$ be as in Lemma \ref{120.50}. Then they have to be neighbours in the first row of the coordinate representation of $\widetilde{R}$. Denote them $a,b$ and let $a=(l,m-1),b=(l-1,m),l+m-1=n$. The element   $c':=(l-1,m-1)$ correspond to the third element $[\al_i+\al_j]$ of the $\{t_1t_2\}$-triangle, which will be called {\em basic} for a moment.

\begin{center}
\setlength{\unitlength}{0,25cm}
\begin{picture}(50,28)
\put(25,25){\vector(-1,-1){25}}
\put(25,25){\vector(1,-1){25}}
\put(25,25){\circle{0,5}}

\put(5,5){\circle{0,5}}
\put(9,5){\circle{0,5}}
\put(13,5){\circle*{0,5}}
\put(17,5){\circle*{0,5}}
\put(21,5){\circle*{0,5}}
\put(25,5){\circle*{0,5}}
\put(29,5){\circle{0,5}}
\put(33,5){\circle{0,5}}
\put(37,5){\circle{0,5}}
\put(41,5){\circle{0,5}}
\put(45,5){\circle{0,5}}

\put(11,7){\circle*{0,5}}
\put(15,7){\circle*{0,5}}
\put(19,7){\circle*{0,5}}
\put(23,7){\circle*{0,5}}
\put(27,7){\circle*{0,5}}
\put(9,9){\circle*{0,5}}
\put(13,9){\circle*{0,5}}
\put(17,9){\circle*{0,5}}
\put(21,9){\circle*{0,5}}
\put(25,9){\circle*{0,5}}
\put(29,9){\circle*{0,5}}

\put(25,9){\circle*{0,5}}

\put(11,11){\circle*{0,5}}
\put(15,11){\circle*{0,5}}
\put(23,11){\circle*{0,5}}

\put(31,11){\circle*{0,5}}

\put(13,13){\circle*{0,5}}
\put(25,13){\circle*{0,5}}
\put(29,13){\circle*{0,5}}
\put(33,13){\circle*{0,5}}
\put(27,15){\circle*{0,5}}
\put(31,15){\circle*{0,5}}
\put(35,15){\circle*{0,5}}
\put(29,17){\circle*{0,5}}
\put(31,19){\circle*{0,5}}


\put(23,7){\circle*{0,5}}

\put(27,11){\circle*{0,5}}
\put(29,13){\circle*{0,5}}
\put(33,17){\circle*{0,5}}

\thicklines
\put(10,12){\line(1,-1){9}}
\put(19,3){\line(1,1){15}}


\put(12,14){\line(1,-1){11}}
\put(23,3){\line(1,1){13}}
\put(8,10){\line(1,-1){7}}
\put(15,3){\line(1,1){17}}

\put(17,5){\line(1,0){4}}

\put(1,2){\makebox(0,0){$x$}}
\put(4,5){\makebox(0,0){$n$}}
\put(10,11){\makebox(0,0){$l$}}
\put(24,26){\makebox(0,0){$0$}}
\put(34,19){\makebox(0,0){$m-1$}}
\put(47,5){\makebox(0,0){$n$}}
\put(49,2){\makebox(0,0){$y$}}
\put(16,17){\makebox(0,0){$k$}}

\put(17,6){\makebox(0,0){$a$}}
\put(21,6){\makebox(0,0){$b$}}
\put(19,9){\makebox(0,0){$c'$}}
\put(23,13){\makebox(0,0){$b'$}}
\put(25,11){\makebox(0,0){$c$}}
\put(27,9){\makebox(0,0){$a'$}}
\put(28,23){\makebox(0,0){$j$}}

\thinlines
\put(17,5){\line(2,1){8}}
\put(19,7){\line(1,0){8}}
\put(16,18){\line(1,-1){15}}
\put(21,5){\line(3,1){6}}
\put(8,4){\line(1,1){20}}
\end{picture}
\end{center}

We claim that for any  $e\in X(a)\cup X(b)\cup X(c')$ (the union of the "X-shapes" with crossings at $a,b,c'$, see Lemma \ref{120.40}) the pair of times $T_{e}$ is equal to one of the pairs
$\{t_1t_2\},\{t_1t_1\},\{t_2t_2\}$.
Indeed, let $e$ lie on the line $x=k$ with $0\le k\le l-2$, i.e. $e$ is one of the points $b':=(k,m-1),c:=(k,m),a':=(k,m+1)$. Then the elements $a,b,c,a',b',c'$ satisfy hypotheses of Lemma \ref{120.20} and our claim follows (analogous arguments work if $e$ lies on the line $y=j$ with $0\le j\le m-2$).

Moreover, among the triangles $\{a,c,b'\},\{b,c,a'\},\{a',b',c'\}$ at least one is a $\{t_1t_2\}$-triangle. In any of these cases we can repeat the above considerations taking this triangle as a basic one. Irrespectively to which of these cases is considered we will obtain new points $e$ such that  $T_{e}$ is equal to one of the pairs
$\{t_1t_2\},\{t_1t_1\},\{t_2t_2\}$. These points will lie on  on $X(a')\cap X(b')\cap X(c)$, i.e. on the intersection  of the lines $x=k,y=n-k+1$ with $\widetilde{R}$, see the figure below (which corresponds to the new basic triangle $\{a,c,b'\}$).

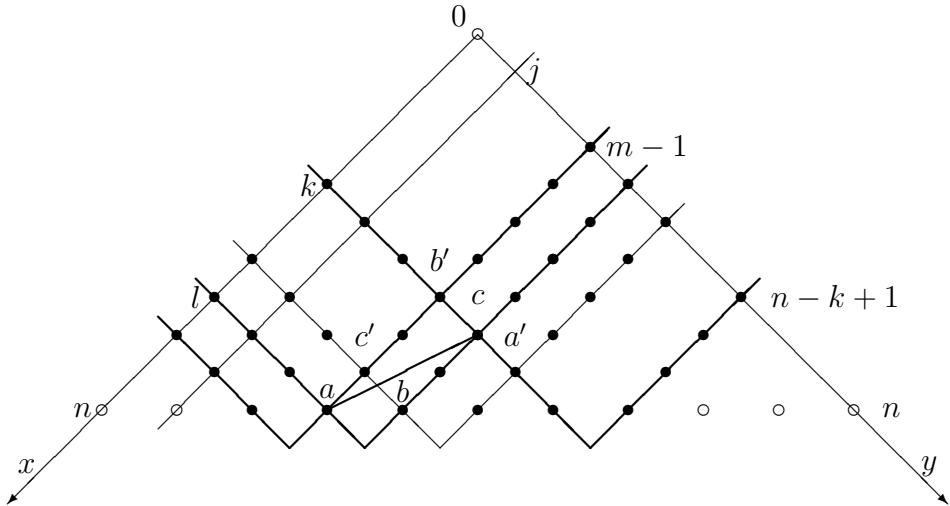
\begin{center}
\setlength{\unitlength}{0,25cm}
\begin{picture}(50,28)
\put(25,25){\vector(-1,-1){25}}
\put(25,25){\vector(1,-1){25}}
\put(25,25){\circle{0,5}}

\put(5,5){\circle{0,5}}
\put(9,5){\circle{0,5}}
\put(13,5){\circle*{0,5}}
\put(17,5){\circle*{0,5}}
\put(21,5){\circle*{0,5}}
\put(25,5){\circle*{0,5}}
\put(29,5){\circle*{0,5}}
\put(33,5){\circle*{0,5}}
\put(37,5){\circle{0,5}}
\put(41,5){\circle{0,5}}
\put(45,5){\circle{0,5}}

\put(11,7){\circle*{0,5}}
\put(15,7){\circle*{0,5}}
\put(19,7){\circle*{0,5}}
\put(23,7){\circle*{0,5}}
\put(27,7){\circle*{0,5}}
\put(35,7){\circle*{0,5}}
\put(9,9){\circle*{0,5}}
\put(13,9){\circle*{0,5}}
\put(17,9){\circle*{0,5}}
\put(21,9){\circle*{0,5}}
\put(25,9){\circle*{0,5}}
\put(29,9){\circle*{0,5}}
\put(37,9){\circle*{0,5}}

\put(25,9){\circle*{0,5}}

\put(11,11){\circle*{0,5}}
\put(15,11){\circle*{0,5}}
\put(23,11){\circle*{0,5}}
\put(39,11){\circle*{0,5}}

\put(31,11){\circle*{0,5}}

\put(13,13){\circle*{0,5}}
\put(21,13){\circle*{0,5}}
\put(25,13){\circle*{0,5}}
\put(29,13){\circle*{0,5}}
\put(33,13){\circle*{0,5}}
\put(19,15){\circle*{0,5}}
\put(27,15){\circle*{0,5}}
\put(31,15){\circle*{0,5}}
\put(35,15){\circle*{0,5}}
\put(17,17){\circle*{0,5}}
\put(29,17){\circle*{0,5}}
\put(31,19){\circle*{0,5}}


\put(23,7){\circle*{0,5}}

\put(27,11){\circle*{0,5}}
\put(29,13){\circle*{0,5}}
\put(33,17){\circle*{0,5}}


\put(12,14){\line(1,-1){11}}
\put(23,3){\line(1,1){13}}
\thicklines
\put(19,3){\line(1,1){15}}
\put(10,12){\line(1,-1){9}}
\put(31,3){\line(1,1){9}}
\put(16,18){\line(1,-1){15}}

\put(8,10){\line(1,-1){7}}
\put(15,3){\line(1,1){17}}


\put(1,2){\makebox(0,0){$x$}}
\put(4,5){\makebox(0,0){$n$}}
\put(10,11){\makebox(0,0){$l$}}
\put(24,26){\makebox(0,0){$0$}}
\put(34,19){\makebox(0,0){$m-1$}}
\put(47,5){\makebox(0,0){$n$}}
\put(49,2){\makebox(0,0){$y$}}
\put(16,17){\makebox(0,0){$k$}}

\put(17,6){\makebox(0,0){$a$}}
\put(21,6){\makebox(0,0){$b$}}
\put(19,9){\makebox(0,0){$c'$}}
\put(23,13){\makebox(0,0){$b'$}}
\put(25,11){\makebox(0,0){$c$}}
\put(27,9){\makebox(0,0){$a'$}}
\put(28,23){\makebox(0,0){$j$}}
\put(44,11){\makebox(0,0){$n-k+1$}}

\put(17,5){\line(2,1){8}}
\thinlines


\put(8,4){\line(1,1){20}}
\end{picture}
\end{center}
Varying $k$ from $0$ to $l-2$ and $j$ from $0$ to $m-2$ we will cover the whole set $\widetilde{R}$. \qed

\medskip

\noindent {\em Proof of Theorem \ref{120.10}}
Let $R\not=\g_2$  be a reduced irreducible root system with the Coxeter graph $\Ga(R)$ (see the beginning of the proof of Theorem \ref{120.30}).          We say that a connected subgraph $\Ga' \subset \Ga(R)$ is a {\em chain} if it does not have ramifying points (i.e. vertices connected with al least three other vertices, see \cite[Ch. IV, App.]{bourb7-8}). In particular, the graphs of $\a_n,\b_n,\c_n,\f_4$ are chains in themselves.

Let $\al_i,\al_j$ be as in Lemma \ref{120.50}.
We can consider any chain $\Ga'\subset \Ga(R)$ containing $\al_i,\al_j$ and proceed as in the proof of Theorem \ref{120.30}
to prove that $T_e$ is equal to one of the pairs
$\{t_1t_2\},\{t_1t_1\},\{t_2t_2\}$ for any $e$ of the form $[\sum_{\be\in Y}\be]$, where $Y$ is a subset of vertices of $\Ga'$ connected as a subgraph. Taking all maximal chains $\Ga'$ containing $\al_i,\al_j$ we will prove in this way that $T_{[\al]}$ is equal to one of the pairs
$\{t_1t_2\},\{t_1t_1\},\{t_2t_2\}$  for all the elements $\al$ of the basis $B$ (and also for all positive roots which are the combinations of the simple roots with coefficients not greater than $1$).

To prove this for all positive roots we will use the fact that any such root can be written as $\be=\be_1+\cdots+\be_k$, where all $\be_j$  (not necessarily distinct) belong to $B$ and each partial sum $\be_1+\cdots+\be_j\in R$ (see \cite[Lem. 3.10, Ch. X]{helgason}). Proceed by induction. We already know that $T_{[\be_j]}$ is equal to one of the pairs
$\{t_1t_2\},\{t_1t_1\},\{t_2t_2\}$. Assume this is true for $T_{[\be_1+\cdots+\be_{j-1}]}$. Then the times selection rules applied to the triangle $\{[\be_1+\cdots+\be_{j-1}],[\be_j],[\be_1+\cdots+\be_{j}]\}$ imply that also $T_{[\be_1+\cdots+\be_{j}]}$ is equal to one of the pairs
$\{t_1t_2\},\{t_1t_1\},\{t_2t_2\}$.

To complete the proof consider the case $R=\g_2$ and proceed by direct inspection. \qed

\section{Bi-Lie structures of Class I}
\label{140}

Let $(\g,[,],[,]')$ be a semisimple bi-Lie structure with a simple Lie algebra $(\g,[,])$ and let  $W:\g\to\g$ be the corresponding  principal WNO. Assume that there exists a Cartan subalgebra $\h$ contained in the central subalgebra $\z$. Then the results of Section \ref{110} show that $W$ defines the time diagram $\T=\T_W$ and that $(W|_{\h},\T)$ is an admissible pair (see Definition \ref{110.50}).

\abz\label{130.05}
\begin{defi}\rm
We say that a regular semisimple bi-Lie structure $(\g,[,],[,]')$ with the principal WNO $W$ satisfying Condition 1 of Theorem \ref{110.10}  is of {\em Class I or II} depending on the class of the pairs diagram $\T_W$ (see Definition \ref{120.15}).
\end{defi}

An algebraic structure which will be introduced below is essentially what we obtain from the structure of the pair groupoid (see eg. \cite{weinstGr}) if we substitute ordered pairs by unordered ones.

\abz\label{140.10}
\begin{defi}\rm Let $X$ be a finite set and let  $X^{\hat{2}}:=\{\{xy\}\mid x,y\in X\}$ be the set of unordered pairs of elements from $X$. Introduce a partial binary operation  $X^{\hat{2}}\times X^{\hat{2}}\to X^{\hat{2}}, (\{xy\},\{zu\})\mapsto \{xy\}\circ\{zu\}$ by the following rules:
\begin{enumerate}\item $\{xy\}\circ\{zu\}$ is defined if and only if the sets $\{x,y\},\{z,u\}$ have a common element; \item $\{xy\}\circ\{yu\}:=\{xu\}$ if $x\not=u$;
\item $\{xy\}\circ\{xy\}$ is either equal to $\{xx\}$ or to $\{yy\}$.\end{enumerate}
The pair $(X^{\hat{2}},\circ)$ will be called a {\em pairoid  of Class I with the base $X$}.
\end{defi}

\abz\label{140.15}
\begin{rema}\rm Note that unlike the pair groupoid structure on the $X\times X$ a pairoid structure is not defined uniquely due to uncertainty of Condition 3.
\end{rema}

Let $(\g,[,])$ be a semisimple Lie algebra, $\h \subset\g$ a fixed Cartan subalgebra. Let $R=R(\g,\h) \subset \h_\R^*$ be the corresponding reduced irreducible root system and $\g=\h+\sum_{\al\in R}\g_\al$ the root grading.

\abz\label{140.20}
\begin{defi}\rm Let $X$ be a finite set.
A decomposition $\g=\bigoplus_{\{xy\}\in X^{\hat{2}}}\g_{\{xy\}}$ such that $\g_{\{xy\}}\not=\{0\}$ for any pair $\{xy\}\in X^{\hat{2}},x\not=y$, is a {\em pairoid quasigrading on $\g$ of Class I with the base $X$} if
\begin{enumerate}\item
$[\g_{\{xy\}},\g_{\{yu\}}]\subset\g_{\{xu\}}$ whenever $x\not=u$;
\item $[\g_{\{xy\}},\g_{\{xy\}}]\subset\g_{\{xx\}}\oplus\g_{\{yy\}}$;
\item $[\g_{\{xy\}},\g_{\{zu\}}]=\{0\}$ whenever $\{xy\}\circ\{zu\}$ is not defined.
\end{enumerate}
A pairoid quasigrading is {\em toral} (with respect to $\h$) if: each $\g_{\{xy\}}$ is the sum of the spaces $\g_\al$ and some subspaces of $\h$; $\h \subset \bigoplus_{x\in X}\g_{xx}$(cf. Definition \ref{ApToral.05}). A toral pairoid quasigrading is {\em symmetric} if $\g_{-\al}\subset\g_{\{xy\}}$ whenever $\g_{\al}\subset\g_{\{xy\}}$.

The elements of the set $X$ are called the {\em times} of the quasigrading. A time $t\in X$ is {\em virtual} if $\g_{\{tt\}}=\{0\}$. Given an unordered pair $\{xy\}\in X^{\hat{2}}$, the element $x$ is called a {\em virtual element}  of the pair if $[\g_{\{xy\}},\g_{\{xy\}}] \subset \g_{\{yy\}}$. In particular, a virtual time is a virtual element of any pair containing it (cf. Examples \ref{140.90}, \ref{140.180}). If we want to distinguish virtual elements we  denote them in brackets (cf. the proof of Theorem \ref{110.30} and examples below).

Two pairoid quasigradings of Class I, $\g=\bigoplus_{\{xy\}\in X^{\hat{2}}}\g_{\{xy\}}$ and $\g=\bigoplus_{\{x'y'\}\in (X')^{\hat{2}}}\g'_{\{x'y'\}}$, with the bases $X$ and $X'$ are {\em equivalent} if there exists $\phi\in\Aut(\g)$ and a bijection $\chi:X\to X'$ such that $\phi(\g_{\{ab\}})=\g'_{\{\chi(a)\chi(b)\}}$ for any $a,b\in X$ (here  $\Aut(\g)$ denotes the set of automorphisms of the Lie algebra $(\g,[,])$).
\end{defi}

For the sake of simplicity we will omit the braces in the notations of unordered pairs indexing  grading subspaces.

\abz\label{140.30}
\begin{exa}
\rm Let $X:=\{x,y\}$, where $x\not=y$, and let $\g=\g_0\oplus\g_1$ be a $\Z_2$-grading on $\g$. Putting $\g_{xx}:=\g_0,\g_{yy}:=\{0\},\g_{x(y)}:=\g_1$ we get a pairoid quasigrading on $\g$ with the base $X$ and a virtual time $y$. An equivalent quasigrading can be given by $\g_{xx}:=\{0\},\g_{yy}:=\g_0,\g_{(x)y}:=\g_1$. An example of quasigrading nonequivalent to that from previous examples is obtained when $\g_0$ is a sum of two nontrivial subalgebras, $\g_0=\g_0^1\oplus\g_0^2$: $\g_{xx}:=\g_0^1,\g_{yy}:=\g_0^2,\g_{xy}:=\g_1$. These examples exhaust up to equivalence all nontrivial pairoid quaisgradings of Class I with the base $X$.

\end{exa}

Let us mention some properties of pairoid quasigrading of Class I.
 Obviously, the notion of pairoid quasigrading of Class I on $\g$ with the base $X$ is independent of the choice of pairoid structure on $X^{\hat{2}}$ (see Remark \ref{140.15}). Moreover,  the subspaces $\g_{xx}$ are subalgebras and if the pairoid quasigrading is toral symmetric these subalgebras  are reductive in $\g$ (by \cite[\S 3, Ch. VII]{bourb7-8}, see also the beginning of Section \ref{130}). The following lemma is a direct consequence of  the definitions.

\abz\label{140.40}
\begin{lemm} Let $X:=\{t_1,\ldots,t_n\}$ be a set of pairwise distinct complex numbers.

Let $(U,\T), \T=\{T_\al\}_{\al\in R}$, be an admissible pair (see Definition \ref{110.50}) such that $T_\T=X,U\in\End(\h)$ and the pairs diagram is of Class I (we transport the root system $R$ to $\h$ by means of the Killing form, see the discussion before Definition \ref{110.255}). Write   $\h_{t_i}$   for the eigenspace of the operator $U$ corresponding to the eigenvalue $t_i$.

Put $\g_{t_it_j}:=\bigoplus_{\al,T_\al=\{t_it_j\}}\g_\al$ if $i\not=j$ and $\g_{t_it_i}:=\h_{t_i}\oplus\bigoplus_{\al,T_\al=\{t_it_i\}}\g_\al$. Then $\g=\bigoplus_{t_i,t_j\in X}\g_{t_it_j}$ is a toral symmetric pairoid quasigrading of Class I with the base $X$.

Vice versa, given a toral symmetric pairoid quasigrading $\g=\bigoplus_{t_i,t_j\in X}\g_{t_it_j}$ of Class I with the base $X$, one obtains an admissible pair $(U,\T)$ with $T_\T=X$ and the pairs diagram of Class I by putting $T_\al:=\{t_it_j\}$ for all  $\al\in R$ such that $\g_\al \subset\g_{t_it_j}$ and $U|_{\g_{t_it_i}\cap\h}=t_i\Id_{\g_{t_it_i}\cap\h}$.

The built correspondence is one-to-one.
\end{lemm}
In the following theorem we will build a WNO from a given toral symmetric pairoid quasigrading of Class I, or, in view of the lemma above, from an admissible pair.
This theorem shows that the necessary conditions for an operator to be a WNO obtained in Section \ref{110} are in fact sufficient in the case of Class I bi-Lie structures.

\abz\label{140.50}
\begin{theo} Let $(\g,[,])$ be a semisimple Lie algebra.
\begin{enumerate}\item
Let $X:=\{t_1,\ldots,t_n\}$ be a set of pairwise distinct complex numbers and let $\g=\bigoplus_{t_i,t_j\in X}\g_{t_it_j}$ be a toral symmetric pairoid quasigrading of class I with the base $X$. Define an operator $W\in\End(\g)$ by
$$
W|_{\g_{t_it_j}}:=\frac{t_i+t_j}{2}\Id_{\g_{t_it_j}}.
$$
Then  the triple $(\g,[,],[,]_W)$ is a semisimple bi-Lie structure of Class I such that
\begin{enumerate}
\item $X$ is its set of times, the subalgebras $\g_{t_it_i}$ coincide with the subalgebras $\g_0^{t_i}$ from Theorem \ref{110.30}, in particular, the  basic subalgebra   is equal to  $\g_0:=\bigoplus_{i=1}^n\g_{t_it_i}$;
\item its central subalgebra $\z$ contains $\h$ and,
moreover, $\z=\bigoplus_{i=1}^n\z^{t_i}$ ($\z^{t_i}$ being the centre of the exceptional bracket $(\g,[,]'-t_i[,])$),  $\z^{t_i}= \g_{t_it_i}\cap\h+\bigoplus_{\al\in\widetilde R^i}\g_\al$, where  $\widetilde R^i:=\{\al\in R^i\mid \al(\bigoplus_{j\not=i}\g_{t_jt_j}\cap\h)=0\}$, here $R^i$ is a closed symmetric root set corresponding to the subalgebra $\g_{t_it_i}$ (see the beginning of Section \ref{130});
\item the operator $W|_{\h^\bot}$ is symmetric, in particular, $W$ is the principal WNO of the constructed bi-Lie structure. \end{enumerate}
\item Any semisimple bi-Lie structure $(\g,[,],[,]')$ of Class I is of the form above.
\item If $W,W'$ are the operators built by two toral symmetric pairoid quasigradings of Class I with the bases $X,X'$, then the corresponding bi-Lie structures are strongly isomorphic if and only if the quasigradings are equivalent, $X=X'$, and $\chi(x)=x$ for any $x\in X$, here $\chi$ is the bijection from  Definition \ref{140.20}. The bi-Lie structures are isomorphic if and only if the quasigradings are equivalent and there exist $\la,\la'\in\C$ such that $\chi(x)=\la x+\la'$ for any $x\in X$.\end{enumerate}
\end{theo}

\noindent {\em Ad 1.} We will first prove that the operator $W$ satisfies the basic equality $T_W(x,y)=[x,y]_P$, where $P\in\End(\g)$ is a symmetric operator preserving the root grading, $P\h \subset\h,P|_{\g_\al}=\pi_\al$, where $\pi_\al$ are some complex numbers, $\pi_\al=\pi_{-\al}$.

Let $\{t_{1,\al}t_{2,\al}\}=\{t_it_j\}$ if $\g_\al \subset \g_{t_it_j}$. Put
 $\la_\al:=(t_{1,\al}+t_{2,\al})/2, P|_\h=0, \pi_\al:=(t_{1,\al}-t_{2,\al})^2/8$ (the last choice is suggested by Item 2 of Theorem \ref{110.20} since $\la_\al=\la_{-\al}$ and $\ka_\al=(\la_\al-\la_{-\al})/2=0$).  Then, obviously,
\begin{equation}\equ\label{formulaa}
\la_\al+\la_{-\al}=t_{1,\al}+t_{2,\al}, \la_\al \la_{-\al}=t_{1,\al}t_{2,\al}+2\pi_\al.
\end{equation}
We claim that
\begin{equation}\equ\label{formulaa1}
(\la_\al-\la_\ga)(\la_\be-\la_\ga)=\pi_\al+\pi_\be-\pi_\ga
\end{equation}
whenever $\al+\be=\ga$ for some $\al,\be,\ga\in R$. To prove this  take into account the commutation relations of pairoid quasigrading of Class I and observe that the equality $\al+\be=\ga$ implies the  inclusions
$$
\g_\al \subset\g_{t_it_j},\g_\be \subset\g_{t_jt_k},\g_\ga \subset\g_{t_kt_i},
$$
where some of the indices $i,j,k$ can be equal. Now we have $\la_\al=(t_i+t_j)/2,\la_\be=(t_j+t_k)/2,\la_\ga=(t_k+t_i)/2,
\pi_\al=(t_i-t_j)^2/8,\pi_\be=(t_j-t_k)^2/8,\pi_\ga=(t_k-t_i)^2/8$ and we can make a direct inspection.

To prove the main equality we  have to consider several cases.
Let $x,y\in\h$. Then, obviously, $T_W(x,y)=0=[x,y]_P$.

 If $x\in\h, y\in\g_\al$ (recall $[x,y]=\al(x)y$),  we have $T_W(x,y)=[Wx,\la_\al y]-W([Wx,y]+[x,\la_\al y]-W(\al(x)y))=\la_\al\al(Wx)y-W(\al(Wx)y+\la_\al\al(x)y-\la_\al\al(x)y)=0$. On the other hand, $[x,y]_P=[x,\pi_\al y]-P[x,y]=\pi_\al\al(x)y-P\al(x)y=0$.

 Now, let $x\in\g_\al,y\in\g_{-\al}$ be such that $B(x,y)=1$. Then (recall $[x,y]=H_\al$) $T_W(x,y)=(\la_\al \la_{-\al}\Id-(\la_\al +\la_{-\al})W-W^2)H_\al=(W-\la_\al\Id)(W-\la_{-\al}\Id)H_\al$. By formula (\ref{formulaa}) the last expression is equal to  $(W-t_{1,\al}\Id)(W-t_{2,\al}\Id)H_\al+2\pi_\al H_\al=2\pi_\al H_\al$  (here we used the fact that  $H_\al$ is a sum of eigenvectors of $W$ corresponding to the eigenvalues $t_{1,\al},t_{2,\al}$ or is an eigenvector corresponding to one of them). On the other hand, $[x,y]_P=[\pi_\al x,y]+[x,\pi_\al y]-PH_\al=2\pi_\al H_\al$.

 If $x\in\g_\al,y\in\g_\be$ with $\al+\be\not\in R\cup\{0\}$, then $T_W(x,y)=0=[x,y]_P$.
\noindent Finally, assume $x\in\g_\al,y\in\g_\be$ with $\al+\be=\ga\in R$. Then (cf. the proof of Lemma \ref{110.40}) $T_W(x,y)=(\la_\al\la_\be-\la_\ga(\la_\al+\la_\be-\la_\ga))[x,y]=
    (\la_\al-\la_\ga)(\la_\be-\la_\ga)[x,y]=(\pi_\al+\pi_\be-\pi_\ga)[x,y]=[x,y]_P$, where we used formula (\ref{formulaa1}).

We have proven that $(\g,[,],[,]_W)$ is a (semisimple) bi-Lie structure which by construction $(\g,[,], \lbr [,]_W)$ is of Class I.

Properties (a), (b) follow  from Theorems \ref{110.20}, \ref{110.30}. Property (c) is a consequence of Item 4 of Theorem \ref{100.20} and the symmetric property of the quasigrading.

\medskip

\noindent  {\em Ad 2.} Item 2 follows from Theorems \ref{110.20}, \ref{110.30}.

\medskip

\noindent  {\em Ad 3.} Item 3 is a consequence of Theorem \ref{40.nonum}.
\qed

\medskip

Now we will show that any pairoid quasigrading of class I induces a special type of
$\Z_2^{m-1}$-grading (here $\Z_p:=\Z/p\Z$) on $(\g,[,])$.  Recall that a {\em basis} of a finite  abelian group $G$ is  a set of elements $e_1,\ldots, e_n\in G$ such that any element $g\in G$ has a unique decomposition of the form $g=k_1e_1+\cdots+k_ne_n$ (we use the additive notation) with $0\le k_i<o_i$, where $o_i$ is the order of the element $e_i$. Clearly, given a basis as above we can define an isomorphism $\phi:G\to\Z_{o_1}\times\cdots\times\Z_{o_n}$ by the formula $g\mapsto(k_1,\ldots,k_n)$. If  $G=G_{m-1}$ is a group isomorphic to $\Z_2^{m-1}$, the cardinality of all the bases of $G$ is the same.

Let $\g=\bigoplus_{i\in G_{m-1}}\g_i$ be a $G_{m-1}$-grading. An element $i\in G_{m-1}$ will be called a {\em quasiroot}  if $\g_{i}\not=\{0\}$.

\abz\label{140.60}
\begin{defi}\rm Given $k,l\in\{1,\ldots,m\},k <l$, put
$I_{kl}:=(i_1,\ldots,i_{m-1})$, where $i_1=\cdots=i_{k-1}=i_l= \cdots =i_{m-1}=0,i_k= \cdots= i_{l-1}=1$.

 We say that a toral (see Appendix \ref{ApToral}) $G_{m-1}$-grading $\g=\bigoplus_{i\in G_{m-1}}\g_i$ is {\em admissible} if there \linebreak exists a basis of $G_{m-1}$ such that    the element $(0,\ldots,0)$ and all the elements $I_{kl},k\in\{1,\ldots,m-1\}, l\in\{2,\ldots,m\},k<l$, are the only quasiroots of the induced $\Z_2^{m-1}$-grading $$\g=\bigoplus_{(i_1,\ldots,i_{m-1})\in\Z_2^{m-1}}\g_{(i_1,\ldots,i_{m-1})}.$$
\end{defi}

 The following lemma is a direct consequence of the definition.

\abz\label{140.70}
\begin{lemm}
Let $X:=\{t_1,\ldots,t_n\}$ be a set of pairwise distinct complex numbers and let $\g=\bigoplus_{t_i,t_j\in X}\g_{\{t_it_j\}}$ be a toral symmetric pairoid quasigrading of class I with the base $X$. Then the formula
 $$
 \g_{(i_1,\ldots,i_{m-1})}:=\left\{ \begin{array}{ll}
 \bigoplus_{i=1}^n\g_{\{t_it_i\}} & \mbox{if}\ (i_1,\ldots,i_{m-1})=(0,\ldots,0);\\
 \g_{\{t_kt_l\}} & \mbox{if}\ (i_1,\ldots,i_{m-1})=I_{kl};\\
 \{0\}& \mbox{in\ other\ cases}.\\
 \end{array}\right.
$$
gives an  admissible $\Z^{m-1}_2$-grading on $\g$.
\end{lemm}

Note that in fact the grading above is determined by the corresponding pairs diagram $\T$ (without participation of the $R$-admissible operator $U$, cf. Lemma \ref{140.40}).

Below we will construct a series of examples of pairoid quasigradings of Class I. In order to do this one should start from constructing admissible $\Z_2^{m-1}$-gradings. By definition any $\Z_2^{m-1}$-grading with $m=2,3$ is so, hence in principle this is a nontrivial problem only for $m>3$. However
 the first of the following examples is in a sense a counterexample since it shows that the correspondence from the preceding lemma is not one-to-one, i.e. finding an admissible toral $\Z_2^{m-1}$-grading could be insufficient for building a pairoid quasigrading (even a pairs diagram).

 It is clear that any toral symmetric $\Z_2^{m-1}$-grading is defined by $m-1$ commuting inner automorphisms of order 2. In the examples below we will use the model diffeomorphisms of different types  (see \cite[Ch. X]{helgason} and discussion after Remark \ref{130.65}). We will also use the notations from Tables I-IV of \cite{bourb4-6}.

\abz\label{140.80}
\begin{exa}\rm Let $\g=\d_4$ and let $\al_1,\al_2,\al_3,\al_4$ be the basis of the root system. Here  $\al_1,\al_3,\al_4$ have labels $1$ and $\al_2$ has label $2$. Consider $\Z_2^2$-grading defined by  the  model automorphisms of type $(0,1,0,0,1;1),(0,0,0,1,1;1)$. Then  we have $\g_{\al_1} \subset\g_{(1,0)},\g_{\al_2} \subset\g_{(0,0)},\g_{\al_3}\subset\g_{(0,1)},\g_{\al_4}\subset\g_{(1,1)}$, or, by the hypothetical inverse correspondence to that of Lemma \ref{140.70}, $\g_{\al_1}\subset\g_{t_1t_2},\g_{\al_3}\subset\g_{t_2t_3},\g_{\al_4} \subset\g_{t_1t_3}$. However, there  is a problem: in which of three components $\g_{t_1t_1},\g_{t_2t_2}$ or $\g_{t_3t_3}$ should lie $\g_{\al_2}$ in order that the axioms of the pairoid quasigrading be satisfied? This problem is related to the ramifying point of the Dynkin diagram.
\end{exa}

The next example shows that nonadmissible toral $\Z_2^{m-1}$-gradings with $m>3$ exist.

\abz\label{140.90}
\begin{exa}\rm
Let $\g=\d_4$. Consider the  model  automorphisms of types $(1,1,0,0,0;1),(1,0,0, \lbr 1,0;1),\lbr (1,0,0,0,1;1)$. They  define a $\Z_2^3$-grading on $\g$ which is not admissible. Indeed, it has eight quasiroots: $\g_{\al_1} \subset\g_{(100)},\g_{\al_2} \subset\g_{(000)},\g_{\al_3} \subset\g_{(010)},\g_{\al_4} \subset\g_{(001)},\g_{\al_1+\al_2+\al_3}\subset\g_{(110)},
\g_{\al_2+\al_3+\al_4}\subset\g_{(011)},\g_{\al_1+\al_2+\al_4}\subset\g_{(101)},
\g_{\al_1+\al_2+\al_3+\al_4}\subset\g_{(111)}$. However, an admissible $\Z_2^3$-grading has only seven quasiroots. Again the ramifying point plays crucial role.
\end{exa}

Below  we give a series of examples of toral symmetric pairoid quasigrading of Class I. The first of these examples is universal, i.e. it appears in any semisimple Lie algebra.

\abz\label{140.95}
\begin{exa}\rm Consider pairoid quasigradings from Example \ref{140.30} such that the corresponding $\Z_2$-grading is induced by an inner automorphism. Then the quasigradings are toral and symmetric and by Theorem \ref{140.50} we get a number of regular semisimple bi-Lie structures of Class I with two times.
\end{exa}

Now we will construct a series of (in a sense canonical, see Conjecture \ref{140.170}) toral symmetric pairoid quasigrading of Class I on the classical Lie algebras. First three series were known in the literature. The forth an fifth one, related to the $\a_n$-series, are new.

We will describe the quasigradings  by means of admissible pairs $(U,\T)$ (cf. Lemma \ref{140.40}).

\abz\label{140.90}
\begin{exa}\rm Let $\g=\b_n=\so(2n+1,\C)$, the roots are $\pm\ep_i,\pm\ep_i\pm\ep_j,1\le i<j\le n$, where $\ep_i$ are the elements of an orthonormal basis in $\h_\R^*$. Put $U(H_{\ep_i}):=t_iH_{\ep_i},T_{\pm\ep_i\pm\ep_j}:=\{t_it_j\}$ and $T_{\pm\ep_i}:=\{t_i(t_{n+1})\}$ (recall that the brackets mean that the corresponding element is virtual, cf. Definition \ref{110.50}, discussion after it, and Definition \ref{140.20}).
The set of times is $\{t_1,\ldots,t_{n+1}\}$ and the time $t_{n+1}$ is virtual.
The labels of the standard basis $\al_1,\ldots,\al_n$ are $1,2,\ldots,2$.
The corresponding admissible $\Z_2^n$-grading is defined by the model automorphisms of types $(1,1,0,0,\ldots,0;1),(0,0,1,0,\ldots,0;1),\lbr (0,0,0,1,\ldots,0;1),\lbr \ldots,(0,0,0,0,\ldots,1;1)$. This example corresponds to Example \ref{0.30} with the diagonal matrix
$A=\diag(t_1,t_1,t_2,t_2,\ldots,t_n,t_n,t_{n+1})$.
\end{exa}

\abz\label{140.100}
\begin{exa}\rm Let $\g=\d_n=\so(2n,\C)$, the roots are $\pm\ep_i\pm\ep_j,1\le i<j\le n$. Put $U(H_{\ep_i}):=t_iH_{\ep_i},T_{\pm\ep_i\pm\ep_j}:=\{t_it_j\}$.
The set of times is $\{t_1,\ldots,t_{n}\}$ and there no virtual times (and elements).
The labels of the standard basis $\al_1,\ldots,\al_n$  are $1,2,\ldots,2,1,1$. The corresponding admissible $\Z_2^{n-1}$-grading is defined by the model automorphisms of types $(1,1,0,0,\ldots,0,0,0;1),\lbr (0,0,1,0,\ldots,0,0,0;1),\lbr (0,0,0, \lbr 1,\ldots,0,0,0;1),\lbr \ldots,(0,0,0,0,\ldots,1,0,0;1),(0,0,0,0,\ldots,0,1,1;1)$. This example corresponds to Example \ref{0.30} with the diagonal matrix
$A=\diag(t_1,t_1,t_2,t_2,\ldots,t_n,t_n)$.
\end{exa}

\abz\label{140.110}
\begin{exa}\rm Let $\g=\d_n=\sp(n,\C)$, then the roots are $\pm 2\ep_i,\pm\ep_i\pm\ep_j,1\le i<j\le n$. Put $U(H_{\ep_i}):=t_iH_{\ep_i},T_{\pm\ep_i\pm\ep_j}:=\{t_it_j\}$ and $T_{\pm 2\ep_i}=\{t_it_i\}$.
The set of times is $\{t_1,\ldots,t_{n}\}$ and there no virtual times (and elements).
The labels of the standard basis  $\al_1,\ldots,\al_n$ are $2,2,\ldots,2,1$. The corresponding admissible $\Z_2^{n-1}$-grading is defined by the model automorphisms of types  $(0,1,0,\ldots,0,0;1),\lbr (0,0,1,\ldots,0,0;1),\lbr \ldots,(1,0,0,\ldots,1,0;1)$. This example corresponds to Example \ref{0.35} with the diagonal matrix
$A=\diag(t_1,t_2,\ldots,t_n,t_1,t_2,\ldots,t_n)$.
\end{exa}

\abz\label{140.120}
\begin{exa}\rm Let $\g=\a_n=\sl(n+1,\C)$, then the roots are $\pm(\ep_i-\ep_j),1\le i<j\le n+1$, where $\ep_i$ are the elements of an orthonormal basis in a $(n+1)$-dimensional euclidean space in which $\h_\R^*$ is embedded as the hyperplane orthogonal to the vector $(1,1, \ldots,1)$. Put $U(H_{\ep_i-\ep_{n+1}}):=t_iH_{\ep_i-\ep_{n+1}}, i=1,\ldots,n,T_{\pm(\ep_i-\ep_j)}:=\{t_it_j\}$, if $i<j<n+1$ and $T_{\pm(\ep_i-\ep_{n+1})}=\{t_i(t_{n+1})\}$.
The set of times is $\{t_1,\ldots,t_{n+1}\}$ and the time $t_{n+1}$ is virtual.

The labels of the standard basis  $\al_1,\ldots,\al_n$ are $1,1,\ldots,1$. The corresponding admissible $\Z_2^{n}$-grading is defined by the model automorphisms of types  $(1,1,0,\ldots,0;1),\lbr (1,0,1, \linebreak \ldots,0;1),\lbr \ldots,(1,0,0,\ldots,1;1)$. The WNO constructed via Theorem \ref{140.50} by these data has the form $WX=(1/2)(L_A+R_A)X-\Tr((1/2)(L_A+R_A)X)B$, where $X\in\sl(n+1),A=\diag(t_1,t_2,\ldots,t_{n+1}),B=\diag(0,0,\ldots,0,1)$ (cf. Example \ref{10.60}).
 Explicitly, if $X=||x_{ij}||\in\sl(n+1,\C)$, then $WX=||y_{ij}||$, where $y_{ij}=x_{ij}(t_i+t_j)/2$ for $i\not=j$, $y_{ii}=x_{ii}t_i$ for $i=1,\ldots,n$, and $y_{n+1,n+1}=-\sum_{j=1}^{n}x_{jj}t_j)$.

The proof of Theorem \ref{140.50} suggests that the WNO $W$ will have the principal primitive equal to that of the operator $Z:=(1/2)(L_A+R_A)$ (considered as an element of $\End(\gl(n+1))$). In other words, the following fact is true: $T_W=T_Z$ (here also $W$ is considered as an operator on $\gl(n+1)$). Its direct proof surprisingly is not evident, and we present it below.

 We have to prove that $T_{Z-V}=T_Z$, where  $VX:=(Z-W)X=\Tr(ZX)B$. We will use the following simple observations (which are true for any $X,Y\in\gl(n+1)$): 1) $Z[B,X]=[B,ZX]$  (since the matrices $A$ and $B$ commute); 2) $Z[VX,Y]=[VX,ZY]$ (as a consequence of 1)); 3) $V^2X=t_{n+1}VX$ and $ZVX=t_{n+1}VX$ (by definitions); 4) $\Tr(Z[X,Y]_Z)=0$ (this follows from the main identity $T_Z(X,Y)=[X,Y]_P$, where $P=\ad (A/2)\circ L_A-(1/8)\ad^2(A)$, see Example \ref{10.60}, since $Z[X,Y]_Z=[ZX,ZY]-[X,Y]_P=[ZX,ZY]-([PX,Y]+[X,PY]-[A/2,A[X,Y]]+(1/8)[A,[A,[X,Y]]])$ is the combination of commutators); 5) $[VX,VY]=0$ (obvious); 6) $\Tr(Z[B,X])=0$ (the matrices $[B,X]$ and $Z[B,X]$ has zeroes on the diagonal); 7) $\Tr(Z[VX,Y])=0$ (as a consequence of 6)).

 Now we have $T_{Z-V}(X,Y)=T_Z(X,Y)-[ZX,VY]-[VX,ZY]+[VX,VY]+Z[X,Y]_V+V[X,Y]_Z-V[X,Y]_V=
 T_Z(X,Y)-[ZX,VY]-[VX,ZY]+[VX,VY]+([VX,ZY]+[ZX,VY]-t_{n+1}V[X,Y])+\Tr(Z[X,Y]_Z)B-
 (\Tr(Z[VX,Y])B+\Tr(Z[X,VY])B-t_{n+1}V[X,Y])=T_Z(X,Y)$.
\end{exa}

The next two examples  contain an additional continuous parameter (i.e. moduli of bi-Lie structures with the fixed set of times are appearing here).

\abz\label{140.140}
\begin{exa}\rm Let $\g=\a_n=\sl(n+1,\C)$. Put $T_{\pm(\ep_i-\ep_j)}:=\{t_it_j\}$, if $i<j<n+1$, and $T_{\pm(\ep_i-\ep_{n+1})}=\{t_it_n\},i=1,\ldots,n$. The corresponding admissible $\Z_2^{n-1}$-grading is defined by the model automorphisms of types  $(1,1,0,\ldots,0;1),\lbr (1,0,1,\ldots,0;1),\lbr \ldots,(1,0,0,\ldots,1,0;1)$.
Put $w_n:=aH_{\al_n}, w_{n-1}:=w_n+H_{\al_{n-1}},\ldots,w_1:=w_2+H_{\al_1}$, where $\al_i:=\ep_i-\ep_{i+1}$ are the elements of the standard basis of the root system and $a$ is a complex parameter, and $U(w_i):=t_iw_i$. It is easy to see that any root vector is a linear combination of not more than two eigenvectors of $U$, i.e. the operator $U$ is $R$-admissible, and that the pair $(U,\T)$ is admissible too.
By Theorem \ref{140.50} we get a family of semisimple bi-Lie structures with the set of times $X:=\{t_1,\ldots,t_n\}$ (there are no virtual times) depending on the parameter $a$. It turns out that  these structures are pairwise nonisomorphic if the parameter $a$ is taken from small vicinity of $1$.

Indeed, in order that two such structures related with quasigradings $\g=\bigoplus_{t_i,t_j\in X}\g^a_{\{t_it_j\}},\g=\bigoplus_{t_i,t_j\in X}\g^{a'}_{\{t_it_j\}}$, corresponding to parameters $a$ and $a'$, were isomorphic it is necessary that there existed an automorphism $\phi\in\Aut(\g)$   such that $\phi(\g^a_{t_nt_n})=\g^{a'}_{t_nt_n}$ and $\phi(\g^a_{t_it_i})=\g^{a'}_{t_{\sigma(i)}t_{\sigma(i)}},i=1,\ldots,n-1$, for some permutation $\sigma\in S_{n-1}$. By the construction the subalgebra $\g^a_{t_nt_n}$ is the $3$-dimensional subalgebra generated by the subspace $\g_{\al_n}$, while the remaining subalgebras $\g^a_{t_it_i}$ are one dimensional subalgebras in $\h$ of specific form. It can be shown that such $\phi$ has to preserve not only the subalgebra $\h+\g_{\al_n}+\g_{-\al_n}$ (which is the zero subalgebra of the  admissible $\Z_2^{n-1}$-grading), but also the Cartan subalgebra $\h$ itself. Thus the automorphisms which could realize an isomorphism between these bi-Lie structures belong to discrete group, while the parameter $a$ gives a "continuous degree of freedom".
\end{exa}

\abz\label{140.150}
\begin{exa}\rm Let $\g=\b_n=\so(2n+1)$. Consider the   $\Z_2$-grading is defined by the model automorphism of type $(1,1,0,0,\ldots,0;1)$. Its zero subalgebra $\g_0$ is the direct sum of the subalgebra $\so(2n-1)=:\g_{t_2t_2}$ embedded standardly in the lower right corner and arbitrary $1$-dimensional   subalgebra $\ss=:\g_{t_1t_1}  \subset\h$ such that $\ss$ is not contained in $\g_{t_2t_2}\cap\h$. It remains to put $\g_{t_1t_2}:=\g_1$ as in Example \ref{140.30}. Situation here is similar to the preceding example. In order that the bi-Lie structures corresponding to subalgebras $\ss,\ss'$ were isomorphic an automorphism $\phi\in\Aut(\g)$   such that $\phi(\g_{t_2t_2})=\g_{t_2t_2}$ and $\phi(\ss)=\ss'$ should exist. However such an automorphism should preserve the Cartan subalgebra $\h$ and should belong to a discrete group.

Similar examples exist for the Lie algebras $\d_n,\e_6,\e_7$ (cf. \cite[Table 6]{vinbergOnishIII}.

\end{exa}

The following simple observation the proof of which follows from the definition of toral symmetric pairoid qiasigrading of Class I allows to produce new regular semisimple bi-Lie structures of class I from known ones.

\abz\label{140.160}
\begin{lemm} Let $X:=\{t_1,\ldots,t_n\}$ and let $\g=\bigoplus_{t_i,t_j\in X}\g_{t_it_j}$ be a toral symmetric pairoid quasigrading of Class I on $(\g,[,])$ with the base $X$. Then
for any surjective map $\mu:X\to Y:=\{u_1,\ldots,u_m\}$ the formula $\g=\bigoplus_{u_i,u_j\in X}\g^\mu_{u_iu_j}$, where $\g^\mu_{u_iu_j}:=\bigoplus_{x\in\mu^{-1}(u_i),y\in\mu^{-1}(u_j)}\g_{xy}$, gives a toral symmetric pairoid quasigrading of Class I on $(\g,[,])$ with the base $Y$.

In particular, if $\mu_i:X\to Y:=\{u_1,u_2\},i=1,\ldots,n$, is given by $\mu_i(t_i):=u_1,\mu_i(t_j):=u_2$ for any $j\not=i$, we get quasigradings $\g=\g^{\mu_i}_{u_1u_1}\oplus\g^{\mu_i}_{u_2u_2}\oplus\g^{\mu_i}_{u_1u_2}$.
\end{lemm}

The second part of this result gives a hope that toral symmetric pairoid quasigrading of Class I can be completely classified since they are "weaved" from quasigradings with two times, the structure of which is understood (cf. Example \ref{140.30}).

\abz\label{140.170}
\begin{conj} Any toral symmetric pairoid quasigrading of Class I on a simple Lie algebra is equivalent to one of that from Examples \ref{140.95}--\ref{140.140} or to their reductions by means of procedure from Lemma \ref{140.160}.
\end{conj}

\abz\label{140.190}
\begin{rema}\rm
All WNOs built in Theorem \ref{140.50} by toral symmetric quasigradings of Class I with the base $\{t_1,\ldots,t_n\}$ linearly depend on these parameters. Hence all the examples above are in fact examples of families of Lie brackets parametrized by the linear space $\C^n$ with corresponding $n$ (cf. the result of Kantor and Persits discussed in Introduction).
\end{rema}

\abz\label{140.200}
\begin{rema}\rm
The WNO of the bi-Lie structure from Example \ref{0.30} in the case when the matrix $A$ is diagonal with a simple spectrum is also related to a specific pairoid quasigrading of Class I on $\g=\so(n,\C)$ which, however, is not toral. Let $E_{ij},i<j$, be the "antisymmtric matrix unit"; put $X=\{t_1,\ldots,t_n\}$ and $\g=\bigoplus_{i<j}\g_{t_it_j}$, where $\g_{t_it_j}:=\C E_{ij}$. Now one can proceed as in Theorem \ref{140.50} for constructing the WNO (and the formulae from the proof give  the primitive).
\end{rema}

We conclude this section by one more example explaining the subtlety of the notions of virtual time and virtual element  (cf. Definition \ref{140.20}).

\abz\label{140.180}
\begin{exa}\rm
consider the pairoid quasigrading from Example \ref{140.90} and apply to it the procedure of Lemma \ref{140.160} with $\mu:\{t_1,\ldots,t_{n+1}\}\to\{t_1,\ldots,t_{n}\}$ given by the identification of  $t_{n+1}$ with $t_n$ (the corresponding bi-Lie structure is that from Example \ref{0.30} with $A=\diag(t_1,t_1,t_2,t_2,\ldots,t_n,t_n,t_{n})$). Then the time $t_n$ is not virtual (since $\g_{t_nt_n}$ is generated by $\g_{\al_n}$) but $t_n$ is a virtual element of all the pairs $\{t_i(t_n)\},i=1,\ldots,n-1$.
\end{exa}

\section{Bi-Lie structures of Class II}
\label{130}

In this section we will build a theory related to pairs diagrams of Class II (see Definitions \ref{110.50},  \ref{120.15}) similar to that from the previous section.

We start from auxiliary results. Let $(\g,[,])$ be a simple Lie algebra (cf. Remark \ref{rema23}), $\g_0 \subset\g$ a reductive in $\g$ Lie subalgebra of maximal rank. In particular, there exists a Cartan subalgebra $\h$ in $\g$ such that $\h \subset\g_0$. In what follows we fix a Cartan subalgebra $\h \subset\g$.

Let $R=R(\g,\h) \subset \h_\R^*$ be the corresponding reduced irreducible root system. Recall \cite{bourb4-6} that a set $R_0 \subset R$ is  {\em closed} if for any $\al,\be\in R_0$ such that $\ga:=\al+\be\in R$ the root $\ga$ belongs to $R_0$. The set $R_0$ is {\em symmetric} if $R_0=-R_0$. An empty set is closed symmetric by definition.

In terms of the root system $R$ any reductive in $\g$ Lie subalgebra $\g_0 \supset\h$ can be described as follows. If $\g=\h+\bigoplus_{\al\in R}\g_\al$ is the corresponding root decomposition, then $\g_0=\h+\bigoplus_{\al\in R_0}\g_\al$, where $R_0 \subset R$ is some closed symmetric root set \cite[VII, \S 3]{bourb7-8}.

\abz\label{130.10}
\begin{defi}\rm We say that a subalgebra $\g_0 \subset\g,\g_0\supset\h$, is {\em admissible} if it is a reductive subalgebra  in $\g$ and
it is not  the fixed point subalgebra of an inner automorphism of order 2. Given a pair of reductive in $\g$ subalgebras such that $\g_0^1\oplus\g_0^2=\g_0$, we say that $(\g_0^1,\g_0^2)$ is an {\em admissible pair of subalgebras} (trivial cases $\g_0^1=\{0\}$ or $\g_0^2=\{0\}$ are admitted). The closed symmetric root set $R_0$ corresponding to an admissible subalgebra $\g_0$ is called {\em admissible} too.
\end{defi}

\abz\label{130.15}
\begin{rema}\rm The fixed point subalgebras of an inner automorphism of order 2 of simple Lie algebras  and their root systems are well known \cite[Chapter X]{helgason}.
\end{rema}

\abz\label{130.155}
\begin{defi}\rm Let $X=\{x,y\}$ be a 2-element set of complex numbers.
A decomposition $\g=\g_{\{xx\}}\oplus\g_{\{yy\}}\oplus\g_{\{xy\}},\g_{\{xy\}}\not=\{0\}$ (recall that $\{xy\}$ denotes the unordered pair of the elements $x,y\in X$), is a {\em pairoid quasigrading on $\g$ of Class II} with the base $X$ if
\begin{enumerate}\item
$[\g_{\{xx\}},\g_{\{xx\}}]\subset\g_{\{xx\}}$;
 $[\g_{\{yy\}},\g_{\{yy\}}]\subset\g_{\{yy\}}$;
  $[\g_{\{xx\}},\g_{\{yy\}}]=\{0\}$;
  \item $[\g_{\{xx\}},\g_{\{xy\}}]\subset\g_{\{xy\}},[\g_{\{yy\}},\g_{\{xy\}}]
      \subset\g_{\{xy\}}$;
\item    $[\g_{\{xy\}},\g_{\{xy\}}]\cap\g_{\{xy\}}\not=\{0\}$.
\end{enumerate}

Given a pairoid quasigrading of Class II with a base $X=\{x,y\}$, the {\em opposite} pairoid quasigrading is obtained by interchanging $x,y$.

A pairoid quasigrading is {\em toral} if each $\g_{\{zu\}}$ is the sum of the spaces $\g_\al$ and some subspaces of $\h$ (cf. Definition \ref{ApToral.05}). A toral pairoid quasigrading is {\em symmetric} if $\g_{-\al}\subset\g_{\{zu\}}$ whenever $\g_{\al}\subset\g_{\{zu\}}$.

Two pairoid quasigradings of Class II $\g=\g_{\{xx\}}\oplus\g_{\{yy\}}\oplus\g_{\{xy\}}$ and $\g=\g'_{\{x'x'\}}\oplus\g'_{\{y'y'\}}\oplus\g'_{\{x'y'\}}$ with the bases $X=\{x,y\}$ and $X'=\{x',y'\}$ are {\em equivalent} if there exists $\phi\in\Aut(\g)$ and a bijection $\chi:X\to X'$ such that $\phi(\g_{\{ab\}})=\g'_{\{\chi(a)\chi(b)\}}$ for any $a,b\in X$ (here  $\Aut(\g)$ denotes the set of automorphisms of the Lie algebra $(\g,[,])$). Two pairoid quasigradings of Class II are antiequivalent if one of them is equivalent to the opposite to another.
\end{defi}

The following lemma relates the definitions above.

\abz\label{130.200}
\begin{lemm}
Let a toral symmetric pairoid quasigrading $\g=\bigoplus_{t_i,t_j\in X}\g_{\{t_it_j\}}$ of Class II with the base $X=\{t_1,t_2\}$ be given. Put $\g_0^1:=\g_{\{t_1,t_1\}},\g_0^2:=\g_{\{t_2,t_2\}}$. Then $(\g_0^1,\g_0^2)$ is an  admissible pair of subalgebras.

Vice versa, given a set $X=\{t_1,t_2\}$ of distinct complex numbers and an admissible pair $(\g_0^1,\g_0^2)$ of subalgebras, the formulae $\g_{\{t_1,t_1\}}:=\g_0^1,\g_{\{t_2,t_2\}}:=\g_0^2,\g_{\{t_1,t_2\}}:=\g_0^\bot$, where $\g_0=\g_0^1\oplus\g_0^2$, define a toral symmetric pairoid quasigrading  of Class II with the base $X$.

The built correspondence is one-to-one up to the transposition $(t_1 t_2)$.
\end{lemm}

\noindent For the proof we only have to mention that the definition of a pairs diagram of Class II excludes the case when $\g_0$ is   the fixed point subalgebra of an inner automorphism of order 2. Indeed, assuming the contrary,  we have $(R\setminus R_0)+(R\setminus R_0) \subset R_0$ and there does not exist any triangle (see Definition \ref{110.50}) $\al,\be,\ga\in R\setminus R_0$, i.e. such that
$T_\al=T_\be=T_\ga=\{t_1t_2\}$. It is also easy to see that the case of the fixed point subalgebra of an automorphism of order 2 is the only which should be excluded. \qed

\medskip

The next lemma follows from definitions.

\abz\label{130.20}
\begin{lemm} Let $X:=\{t_1,t_2\},t_1\not=t_2, t_i\in\C$.

Let $(U,\T), \T=\{T_\al\}_{\al\in R}$ be an admissible pair (see Definition \ref{110.50}), where $T_\T=X$, $\T$ is the pairs diagram of Class II, and $U\in\End(\h)$ (we transport the root system $R$ to $\h$ by means of the Killing form, see the discussion before Definition \ref{110.255}). Write   $\h_{t_i}$   for the eigenspace of the operator $U$ corresponding to the eigenvalue $t_i$.

Put $\g_{t_it_j}:=\bigoplus_{\al,T_\al=\{t_it_j\}}\g_\al$ if $i\not=j$ and $\g_{\{t_it_i\}}:= \h_{t_i}\oplus\bigoplus_{\al,T_\al=\{t_it_i\}}\g_\al$. Then $\g=\bigoplus_{t_i,t_j\in X}\g_{\{t_it_j\}}$ is a toral symmetric pairoid quasigrading of Class II with the base $X$.

Vice versa, given a toral symmetric pairoid quasigrading $\g=\bigoplus_{t_i,t_j\in X}\g_{\{t_it_j\}}$ of Class II with the base $X$, one obtains an admissible pair $(U,\T)$ with $T_\T=X$ and the pairs diagram of Class II by putting $T_\al:=\{t_it_j\}$ for all  $\al\in R$ such that $\g_\al \subset\g_{\{t_it_j\}}$ and $U|_{\g_{\{t_it_i\}}\cap\h}=t_i\Id_{\g_{\{t_it_i\}}\cap\h}$.

The built correspondence is one-to-one.
\end{lemm}

\abz\label{130.255}
\begin{rema}\rm We introduced the notion of a toral symmetric pairoid quasigrading of Class II to show the parallelism with the theory built in the previous section. However, since by Lemma \ref{130.200} this notion is almost equivalent to the notion of an admissible pair of subalgebras, we will formulate the further results entirely in terms of this last (see also Remark \ref{130.58})

\end{rema}

The following theorem says that  in the particular case of the bi-Lie structures of Class II the necessary conditions obtained in Section \ref{110} are in fact sufficient for an operator to be a WNO. This theorem takes into account the information  about the restriction $W|_\h$ of this operator to the Cartan subalgebra $\h$ and the symmetric part of the operator $W|_{\h^\bot}$ as well as the information about the  antisymmetric part of this operator described in Corollary \ref{110.70}.

\abz\label{130.50}
\begin{theo}
Let $(\g,[,])$ be a simple Lie algebra.
\begin{enumerate}\item  Let an ordered pair $(t_1,t_2)$ of distinct complex numbers be given and an admissible pair of Lie subalgebras $(\g_0^1,\g_0^2)$. Let $\g_0=\g_0^1\oplus\g_0^2$ and let
\begin{equation}\equ\label{grading1}
\g=\bigoplus_{i\in \Ga(\g_0)}\g_i
 \end{equation}
 be the toral irreducible $\Ga(\g_0)$-grading corresponding to the reductive Lie subalgebra $\g_0$ (see Definition \ref{ApToral.20}). Assume that an  operator $W\in\End(\g)$, which  is scalar on the subspaces $\g_i,i\not=0$, and   $\g_0^j,j=1,2$, satisfies the following conditions:
\begin{enumerate}\item   $W|_{\g_0^j}=t_j\Id_{\g_0^j},j=1,2$;
\item the symmetric part $W_{\mathrm{s}}$ of the operator $W|_{\g_0^\bot}$ is equal to $((t_1+t_2)/2)\Id_{\g_0^\bot}$, here $\g_0^\bot=\bigoplus_{i\in \Ga,i\not=0}\g_i$;
\item The extension  by zero of the antisymmetric part $W_{\mathrm{a}}$ of the operator $W|_{\g_0^\bot}$, which we denote by $\overline{W}_{\mathrm{a}}$,  obeys the $((t_1-t_2)/2)$-triangle rule subject to the grading  (see Definition \ref{110.70}). \end{enumerate}

Then the principal projection $\pr(W)$ (see Theorem \ref{100.260}) of the operator $W$ also satisfies Conditions (a),(b),(c), and the triple $(\g,[,],[,]_W)=(\g,[,],[,]_{\pr(W)})$ is a semisimple bi-Lie structure of Class II such that
\begin{enumerate}
\item[(a')] $\{t_1,t_2\}$ is its set of times, the subalgebras $\g_0^i$ coincide with the subalgebras $\g_0^{t_i}$ from Theorem \ref{110.30}, in particular, its basic subalgebra   is   $\g_0$;
\item[(b')] its central subalgebra $\z$ contains $\h$ and,
moreover, $\z=\z^{t_1}\oplus\z^{t_2}$ ($\z^{t_i}$ being the centre of the exceptional bracket $(\g,[,]'-t_i[,])$),  $\z^{t_i}= \g_0^i\cap\h+\bigoplus_{\al\in\widetilde R^i}\g_\al$, where  $\widetilde R^1:=\{\al\in R^1\mid \al(\g_0^2\cap\h)=0\},\widetilde R^2:=\{\al\in R^2\mid \al(\g_0^1\cap\h)=0\}$, here $R^i$ is a closed symmetric root set corresponding to the subalgebra $\g_0^i$. \end{enumerate}
\item Any semisimple bi-Lie structure of Class II is of the form above.

\end{enumerate}

\end{theo}

\noindent  {\em Ad 1.}  We will first prove that the operator $W$ satisfies the basic equality $T_W(x,y)=[x,y]_P$, where $P\in\End(\g)$ is a symmetric operator preserving the root decomposition (\ref{grading}), $P\h \subset\h,P|_{\g_\al}=\pi_\al\Id_{\g_\al}$, where $\pi_\al$ are some complex numbers, $\pi_\al=\pi_{-\al}$.

Let $t_{1,\al}=t_{2,\al}=t_i$ if $\al\in R^i$ and $t_{1,\al}=t_1,t_{2,\al}=t_2$ if $\al\not\in R^1\cup R^2$. Let $\la_\al$ be the eigenvalue of $W$ corresponding to the eigenspace $\g_\al$ and let $\sigma_\al:=(1/2)(\la_\al+\la_{-\al}),\kappa_\al:=(1/2)(\la_\al-\la_{-\al})$ be the corresponding eigenvalues of $W_{\mathrm{s}},W_{\mathrm{a}}$, in particular, $\si_\al=(t_{1,\al}+t_{2,\al})/2$ for any $\al$ by condition (b), $\ka_\al=0$ for $\al\in R^1\cup R^2$ by condition (c).

Put
 $P|_\h=0, \pi_\al:=((t_{1,\al}-t_{2,\al})/2)^2-\ka_\al^2)/2$ (the last choice is suggested by Item 2 of Theorem \ref{110.20}).  One checks  the following formulae:
$$
\la_\al+\la_{-\al}=t_{1,\al}+t_{2,\al}, \la_\al \la_{-\al}=\sigma^2_\al-\ka^2_\al=t_{1,\al}t_{2,\al}+2\pi_\al.
$$
We claim that
$$
(\la_\al-\la_\ga)(\la_\be-\la_\ga)=\pi_\al+\pi_\be-\pi_\ga
$$
whenever $\al+\be=\ga$ for some $\al,\be,\ga\in R$. This can be proven by direct inspection taking into account conditions (a), (b) and considering the following cases
\begin{itemize}\item
 $\al,\be,\ga\in R^i\Rightarrow\si_\al=\si_\be=\si_\ga=t_i,\ka_\al=\ka_\be=\ka_\ga=0,
 \pi_\al=\pi_\be=\pi_\ga=0$;
\item $\al\in R^i,\be,\ga\not\in R^1\cup R^2\Rightarrow\si_\al=t_i,\si_\be=\si_\ga=(t_1+t_2)/2,\ka_\al=0,\ka_\be=\ka_\ga,
    \pi_\al=0,\pi_\be=
    \pi_\ga=((t_{1}-t_{2})/2)^2-\ka_\be^2)/2$;
\item $\al,\be\not\in R^1\cup R^2,\ga\in R^i,\Rightarrow\si_\al=\si_\be=(t_1+t_2)/2,\si_\ga=t_i,\ka_\al=-\ka_\be,\ka_\ga=0,
 \pi_\al=\pi_\be=((t_1-t_2)/2)^2-\ka_\al^2)/2,\pi_\ga=0$;
\item $\al,\be,\ga\not\in R^1\cup R^2,\Rightarrow\si_\al=\si_\be=\si_\ga=(t_1+t_2)/2,\ka_\al+\ka_\be-\ka_\ga=
    \pm(t_1-t_2)/2, \pi_\delta=((t_1-t_2)/2)^2-\ka_\delta^2)/2$, where $\delta=\al,\be,\ga$.
\end{itemize}

Now the proof of the main equality $T_W(x,y)=[x,y]_P$ is literally the same as in the proof of Theorem \ref{140.50}.

Hence we have proven that $(\g,[,],[,]_W)$ is a semisimple bi-Lie structure. Now let us show that $\pr(W)$ satisfies conditions (a), (b), (c). Indeed, since the projecting does not affect $W|_\h$ and $W_{\mathrm{s}}$, conditions (a), (b)  follow from the assumptions. The antisymmetric parts of $W|_{\h^\bot}$ and $\pr(W)|_{\h^\bot}$ differ by an operator of the form $\ad H,H\in\h$, due to Theorem \ref{100.260}, hence obey the same $((t_1-t_2)/2)$-triangle rule by Lemma \ref{130.40}, below. Finally,
the condition  $\ka_\al'=0,\al\in R^1\cup R^2$, where $\ka_\al'$ is the eigenvalue of the antisymmetric part of $\pr(W)|_{\h^\bot}$ corresponding to $\g_\al$, is a consequence of Item 4 of Theorem \ref{110.30} and the fact that $[,]_W=[,]_{\pr(W)}$.

It is clear that $(\g,[,],[,]_W)$ is of Class II. Property (a') follows from Theorem \ref{110.20}. Property (b') is a consequence of  Item 4 of this theorem.

\medskip

\noindent  {\em Ad 2.} Item 2 follows from Theorems \ref{110.20}, \ref{110.30} and Corollary \ref{110.70}. \qed

\medskip

Later we will used the theorem above to construct a series of examples of bi-Lie structures of Class II. In order  to formulate the second our next theorem (\ref{130.500}), dealing with the uniqueness questions, we need two more lemmata.

\medskip
\abz\label{130.40}
\begin{lemm} Let $\h \subset\g$ be a Cartan subalgebra, $\g=\h+\bigoplus_{\al\in R}\g_\al$ the corresponding root decomposition.
Assume $R_0 \subset R$ is an admissible root set (see Definition \ref{130.10}) and $S:\h^\bot\to \h^\bot$ is an antisymmetric  operator preserving the root spaces $\g_\al$ and obeying the  $a$-triangle rule subject to $R_0$ (see Definition \ref{130.30}). Then
\begin{enumerate}
\item given any basis $B=\{\al_1,\ldots,\al_n\}$  of the root system $R$, the operator $S$ can be reconstructed from $a$, its system of labels $\L_S$, and its eigenvalues $\ka_{\al_1},\ldots,\ka_{\al_n}$ corresponding to the eigenspaces $\g_{\al_1},\ldots,\g_{\al_n}$;
\item an antisymmetric operator $S':\h^\bot\to\h^\bot$ preserving the root spaces $\g_\al$ obeys the $a$-triangle rule subject to the same root set with $\L_{S'}=\L_S$ if and only if $S'=S+(\ad H)|_{\h^\bot}$ for some $H\in\h$ such that $\al(H)=0$ for any $\al\in R_0$;
\item in the case when $S$ is pricipal (i.e. if $S=L|_{\h^\bot}$ for some  principal operator  $L\in\End(\g)$, see Definition \ref{defiMain}), any  other principal  antisymmetric operator $S'\in\End(\h^\bot)$ preserving the root spaces $\g_\al$ and obeying  the $a$-triangle rule subject to the same root set with $\L_{S'}=\L_S$  coincides with $S$;
    moreover,
     $S$ can be uniquely reconstructed   from  $a$ and its system of labels $\L_S$.
\end{enumerate}

\end{lemm}

\noindent {\em Ad 1.} Proceeding by induction assume that $\ka_\al$, the eigenvalue of $S$ corresponding to the eigenspace $\g_\al$, is already reconstructed for all positive roots of height $\hei(\al)=k$. Since any positive root $\al$ of height $k+1$ can be decomposed as a sum $\al=\be+\ga$, where $\be,\ga$ are positive roots with $\hei(\be),\hei(\ga)\le k$, we can put $\ka_\al=\ka_\be+\ka_\ga$ or $\ka_\al=\ka_\be+\ka_\ga\pm a$ depending of the label of the triangle $\al,-\be,-\ga$.

\medskip

\noindent {\em Ad 2.} Obviously, for any $H\in\h$ the operator $\ad H|_{\h^\bot}$ is diagonal: $\ad H|_{\g_\al}=\al(H)\Id_{\g_\al}$. Since for any triangle $\al,\be,\ga$ we have $\al(H)+\be(H)+\ga(H)=0$ (i.e. $\ad H$ satisfies the $0$-triangle rule), adding of $\ad H$ to $S$ will not affect the system of labels $\L_S$. On the other hand, the condition $\al(H)=0,\al\in R_0$, guarantees that $\ad H|_{\g_\al}=0$ and $(S+\ad H)|_{\g_\al}=0$ for any $\al\in R_0$. The other implication will be proven below.

We will first show that for the operator $S:=(\ad H)|_{\h^\bot}$ the vector $(\ka_{\al_1},\ldots,\ka_{\al_n})$ runs through all $(z_1,\ldots,z_n)\in\C^n$ as $H$ runs through $\h$. Indeed, we have $\ka_{\al_i}=\al_i(H)$ and in order to find $H\in\h$ such that $\al_i(H)=z_i,i=1,\ldots,n$, we have to solve a system of linear equations with the matrix $\al_i(H_{\al_j})$ the determinant of which is proportional to that of the Cartan matrix of the base $B$, hence is nonzero.

We can now prove the second implication. The operator $S'':=S'-S$ obeys the $0$-triangle rule.  Let  $B=\{\al_1,\ldots,\al_n\}$ be a basis of the root system $R$ and let $\ka''_\al$ be the corresponding eigenvalues of $S''$. The considerations above show that there exists $H\in \h$ such that $(\ad H)|_{\g_{\al_i}}=\ka_{\al_i}'',i=1,\ldots,n$. However  by Item 1 the operators
$S'',\ad H$ coincide. Moreover  $\al(H)=0$ for any $\al\in R_0$ since $S''|_{\g_\al}=0$.

\medskip

\noindent {\em Ad 3.} By Item 2 we have $S'=S+(\ad H)|_{\h^\bot}$ for some $H\in\h$. However, both $S,S'$ are principal and by definition $H=0$.

Now let $S=L|_{\h^\bot}$ for some  principal operator  $L\in\End(\g)$. Starting the induction process of the proof of Item 1 from arbitrary vector $(c_1,\ldots,c_n)$ we will obtain one of the operators $S+(\ad H)|_{\h^\bot}, H\in\h$. To get $S$ one has to project $L+\ad H$ onto $\tilde\g^\bot$ along $\tilde\g=\ad\g \subset\End(\g)$ and to
restrict the result to $\h^\bot$ (cf. Theorem \ref{100.260}).
\qed

\medskip

 Let   $\phi\in\Aut(\g)$ and let $\g_0 \subset\g$ be an admissible subalgebra (see Definition \ref{130.10}). Put $\g_0'=\phi(\g_0)$. If $\h \subset\g_0$ is a Cartan subalgebra, $Q \subset\h_\R^*$ the corresponding root lattice, and $Q_0 \subset Q$ the sublattice related with the subalgebra $\g_0$, put also $\h':=\phi(\h),Q':=\tilde\phi(Q)$ and $Q_0':=\phi(Q_0)$, here $\tilde\phi:=((\phi|_\h)^*)^T=((\phi|_\h)^{-1})^T$ is the operator transposed to the conjugate with respect to the Killing form to $\phi|_\h$.
Clearly $\h' \subset\g_0'$ is a  Cartan subalgebra, $Q' \subset(\h'_\R)^*$ is the corresponding root lattice, and $Q_0'$ is the sublattice related to the subalgebra $\g_0'$.

It is easy to see that $\tilde\phi$ induces the isomorphism of groups $\bar\phi:\Ga\to\Ga'$, where $\Ga:=\Ga(\h,\g_0)=Q/Q_0$ and $\Ga':=\Ga(\h',\g_0')=Q'/Q'_0$ (see Appendix \ref{ApToral}) and the operator $\phi$ preserves the corresponding toral irreducible gradings $\g=\bigoplus_{i\in\Ga}\g_i=\bigoplus_{i\in\Ga'}\g'_i$, or, more precisely, $\phi(\g_i)=\g'_{\bar\phi(i)}$.

If  $i,j,k\in \Ga$ is a triangle,
the triple $\bar\phi(i),\bar\phi(j),\bar\phi(k)$  is obviously again a triangle (see Definition \ref{110.60}). Moreover, if $\L_S$ is the system of labels of an operator $S\in \End(\g)$ diagonal with respect to the grading $\g=\bigoplus_{i\in\Ga}\g_i$, we can endow the triangle $\bar\phi(i),\bar\phi(j),\bar\phi(k)$ with the same label as the label of $i,j,k$ from the family $\L_S$ and obtain a new system of labels, which will be denoted $\bar\phi\L_S$. Also, by $-\L_S$ we will denote the system of labels opposite to $\L_S$ (i.e. with interchanged pluses and minuses).

\abz\label{130.490}
\begin{lemm}  Let $\g_0,\g_0' \subset\g$ be admissible subalgebras such that $\phi(\g_0)=\g_0'$ for some $\phi\in\Aut(\g)$. Let $S,S'\in \End(\g)$ be principal (see Definition \ref{defiMain})  operators diagonal with respect to the toral irreducible gradings $\g=\bigoplus_{i\in\Ga(\g_0)}\g_i,\g=\bigoplus_{i\in\Ga(\g_0')}\g'_i$ and obeying the $a,a'$-triangle rule subject to these gradings respectively, where $a\not=0,a'\not=0$.

Then the following conditions are equivalent
\begin{enumerate}\item $S'\circ\phi:=\phi\circ  {S}$;
\item either $a=a'$ and $\L_{S'}=\bar\phi\L_S$ or
    $a=-a'$ and $\L_{S'}=-\bar\phi\L_S$.
    \end{enumerate}
\end{lemm}

\noindent Assume $S'\circ\phi:=\phi\circ  {S}$. Let $\ka_i,\ka'_i$ be the eigenvalues of the corresponding operators related to the eigenspaces $\g_i,\g'_i$. Since $\phi(\g_i)=\g'_{\bar\phi(i)}$, then $\ka_i=\ka'_{\bar\phi(i)}$ and, given any triangle $i,j,k$, we have $\ka_i+\ka_j+\ka_k=\ka_{\bar\phi(i)}+\ka_{\bar\phi(j)}+
\ka_{\bar\phi(k)}$. Thus either $a=a'$ and the labels of the corresponding triangles coincide, or $a=-a'$ and they are opposite.

Conversely, let condition 2 hold. Put $S'':=\phi\circ  {S}\circ\phi^{-1}$ and let $\ka''_i$ be the corresponding eigenvalues. Then $\ka_i=\ka''_{\bar\phi(i)}$ and in both cases $S''$ obeys the $a'$-triangle rule with the same labels as ${S}'$.
Now we have to choose a Cartan subalgebra $\h \subset\g_0$ and observe that the operators $S',S''$ preserve the root decomposition of $\g$ with respect to the Cartan subalgebra $\h':=\phi(\h)$ and their restrictions $S'|_{(\h')^\bot},S''|_{(\h')^\bot}$ obey the $a'$-triangle rule subject to the admissible root set $R_0'$ corresponding to the subalgebra $\g_0'$ (see Definition \ref{130.30}).

On the other hand, it is easy to see that $S''$ is principal (indeed, $\Tr(\phi\circ  {S}\circ\phi^{-1}\circ\ad x)=\Tr(S\circ\phi^{-1}\circ\ad x\circ\phi)=\Tr(S\circ\ad(\phi^{-1}x))=0$ since $\Tr(S\circ\ad x)=0$ for any $x\in\g$). By Item 3 of Lemma \ref{130.40} we conclude that $S'=S''$ (cf. Remark \ref{110.65}).  \qed

\abz\label{130.500}
\begin{theo}
Two bi-Lie structures, $(\g,[,],[,]_W)$ and $(\g,[,],[,]_{W'})$, each of which is constructed as in Theorem \ref{130.50} by data $(t_1,t_2,\g_0^1,\g_0^2,\L)$ and $(t_1',t_2', (\g_0^1)', (\g_0^1)',\L')$, where $\L,\L'$ denote the systems of labels related to the antisymmetric parts of the corresponding operators, are strongly isomorphic (see Definition \ref{0.15}) if and only if the following conditions hold:

\noindent either
\begin{enumerate}\item[1a.]
there exists $\phi\in\Aut(\g)$ such that
$\phi(\g_0^1)=(\g_0^1)',\phi(\g_0^2)=(\g_0^2)'$ and $\bar\phi\L=\L'$, see notations introduced before Lemma \ref{130.490};
\item[2a.] $t_1=t'_1,t_2=t'_2$;\end{enumerate}
or
\begin{enumerate}
\item[1b.] there exists $\phi\in\Aut(\g)$ such that
$\phi(\g_0^1)=(\g_0^2)',\phi(\g_0^2)=(\g_0^1)'$ and $\bar\phi\L=-\L'$;
\item[2b.] $t_1=t'_2,t_2=t'_1$.\end{enumerate}
The bi-Lie structures $(\g,[,],[,]_W),(\g,[,],[,]_{W'})$ are isomorphic if and only if either Condition 1a or 1b is satisfied.
\end{theo}

\noindent   The "strong" case follows from the construction of the operators $W,W'$, Lemma \ref{130.490}, and Theorem \ref{40.nonum}. To manage the general case we have only to prove that the {\em principal} WNOs $W$ and $W'$ built by data $(t_1,t_2,\g_0^1,\g_0^2,\L)$ and $(t_1',t_2', (\g_0^1)', (\g_0^1)',\L')$ are related by the formula $W=\la W'+\mu\Id_\g$ for some $\la,\mu\in\C$ if and only if either $\g_0^1=(\g_0^1)',\g_0^2=(\g_0^2)'$ and $\L=\L'$, or $\g_0^1=(\g_0^2)',\g_0^2=(\g_0^1)'$ and $\L=-\L'$.

Indeed, if $W=\la W'+\mu\Id_\g$, the operators $W$ and $W'$ should have the same eigenspaces $\g_0^1,\g_0^2,\g_j,j\in\Ga,j\not=0$. In particular,  either  1) $t_i=\la t'_{i}+\mu$ and $\g_0^i=(\g_0^i)', i=1,2$; or  2) $t_i=\la t'_{i'}+\mu$ and $\g_0^i=(\g_0^{i'})'$, where $i=1,2,i'=2,1$. Since $\Id_\g$ does not affect the antisymmetric part of an operator,  the systems of labels of the corresponding antisymmetric operators  $\overline{W}_{\mathrm{a}},\overline{W}_{\mathrm{a}}'$ (see Theorem \ref{130.50}) should coincide in Case 1) or be opposite in Case 2).

Conversely, assume principal WNOs $W,W'$ are built by data $(t_1,t_2,\g_0^1,\g_0^2,\L)$ and $(t_1',t_2', \g_0^1, \g_0^2,\L)$. Since $t_1'\not=t_2'$, the system of equations $t_i=\la t'_{i}+\mu, i=1,2$, has a unique solution $\la_0,\mu_0$ with $\la_0\not=0$. We have decompositions  $W=W|_{\g_0}+W_{\mathrm{s}}+W_{\mathrm{a}},W'=W'|_{\g_0}+W'_{\mathrm{s}}+W'_{\mathrm{a}}$ and the obvious  equalities $W|_{\g_0}=\la_0 W'|_{\g_0}+\mu_0\Id_{\g_0},W_{\mathrm{s}}=\la_0 W'_{\mathrm{s}}+\mu_0\Id_{\g_0^\bot}$. Thus it remains only to show that $W_{\mathrm{a}}=\la_0W'_{\mathrm{a}}$.

To this end consider the operator $W''\in\End(\h^\bot)$ which is the extension by zero to $\h^\bot\setminus \g_0^\bot$ of the operator  $(1/\la_0)W_{\mathrm{a}}\in\End(\g_0^\bot)$. Then $W''$ is the restriction to $\h^\bot$ of the principal operator $(1/\la_0)\overline{W}_{\mathrm{a}}$, obeys the $((t'_1-t'_2)/2)$-rule subject to the same admissible root set with the same labels as $W'$, hence  $W''=W_{\mathrm{a}}'$ on $\g_0^\bot$ by Item 3 of Lemma \ref{130.40}.

The case when $W,W'$ are built by data $(t_1,t_2,\g_0^1,\g_0^2,\L)$ and $(t_1',t_2', \g_0^2, \g_0^1,-\L)$ is similar.
\qed

\abz\label{130.58}
\begin{rema}\rm One can reformulate Theorem \ref{130.50} by saying that the corresponding bi-Lie structure is built by a toral symmetric pairoid quasigrading of Class II with the base $X=\{t_1,t_2\}$ and by a system of labels $\L$. Theorem \ref{130.500} can be also reformulated in this spirit. For instance, the second part of it says that the bi-Lie structures built by two pairoid quasigradings  and systems of labels $\L,\L'$ are isomorphic if and only if either the pairoid quasigradings are equivalent (see Definition \ref{130.155}) and $\bar\phi\L=\L'$, where $\phi$ is the automorphism realizing the equivalence, or antiequivalent and $\bar\phi\L=-\L'$.

\end{rema}

Now we will use Theorem  \ref{130.50} to construct a series of examples of bi-Lie structures of Class II.

Let $\g_0 \subset\g$ be  a Levi subalgebra.
 Choose  a Cartan subalgebra $\h \subset\g$ such that $\h \subset\g_0$ and a basis $B$ of the  root system $R(\g,\h)$ such that the corresponding root system $R_0$ of $\g_0$ is generated over $\Z$ by some subset $B_0 \subset B$  (see Appendix \ref{ApToral}).
Put $R^\pm$ for the set of positive (negative) roots with respect to $B$. Then $\g_0^\bot=\g_0^+\oplus\g_0^-$, where $\g_0^\pm:=\sum_{\al\in R^\pm\setminus R_0}\g_\al$. The subspace $\p^\pm(\g_0,B):=\g_0\oplus\g_0^\pm$ is a parabolic subalgebra.

In the following theorem we generalize Example \ref{10.35}.

\abz\label{130.60}
\begin{theo}
Let $\g_0$ be an admissible subalgebra which is a Levi subalgebra and let $(\g_0^1,\g_0^2), \g_0^1\oplus\g_0^2=\g_0$, be an admissible pair of subalgebras. Choose  a Cartan subalgebra $\h \subset\g$ such that $\h \subset\g_0$ and a basis $B$ of the  root system $R(\g,\h)$ such that the corresponding root system $R_0$ of $\g_0$ is generated over $\Z$ by some subset $B_0 \subset B$. Let $\g=\g_0^+\oplus\g_0\oplus\g_0^-$ be the corresponding decomposition and  let $t_1,t_2\in\C,t_1\not=t_2$.  Define an operator  $W\in\End(\g)$ by:
\begin{enumerate}\item   $W|_{\g_0^j}=t_j\Id_{\g_0^j},j=1,2$;
\item  $W|_{\g_0^+}=t_1\Id_{\g_0^+},W|_{\g_0^-}=t_2\Id_{\g_0^-}$.
\end{enumerate}

Then the operator $W$ satisfies assumptions of Item 1 of Theorem \ref{130.50} and the triple $(\g,[,],[,]_W)$ is a semisimple bi-Lie structure of Class II satisfying conditions 1 (a'), (b') of this theorem.

Moreover, two bi-Li structures $(\g,[,],[,]_W)$ and $(\g,[,],[,]_{W'})$ constructed by data $(t_1,t_2,\g_0^1,\g_0^2,\h,B)$ and $(t'_1,t'_2,(\g_0^1)',(\g_0^2)',\h',B')$, respectively, are isomorphic if and only if either there exists $\phi\in\Aut(\g)$ such that $\phi(\g_0^j)=(\g_0^j)',j=1,2$, and $\phi(\p^+(\g_0^1\oplus\g_0^2,B))=\p^+((\g_0^1)'\oplus(\g_0^2)',B')$ or there exists $\phi\in\Aut(\g)$ such that $\phi(\g_0^1)=(\g_0^2)',\phi(\g_0^2)=(\g_0^1)'$, and $\phi(\p^+(\g_0^1\oplus\g_0^2,B))=\p^-((\g_0^1)'\oplus(\g_0^2)',B')$.
\end{theo}

\noindent First note that the subspaces $\g_0^\pm$ are invariant with respect to the action of the subalgebra $\g_0$ on $\g_0^\bot=\g_0^+\oplus\g_0^-$, i.e. $\g_0^\pm$ are direct sums of the components of the toral irreducible $\Ga(\g_0)$-grading. Obviously, the operator $W$ is scalar on the components of this grading and the subalgebras $\g_0^j$ and satisfies conditions 1 (a), (b) of Theorem \ref{130.50}.

Thus we have only to show that the operator $\overline{W}_{\mathrm{a}}$  obeys the $((t_1-t_2)/2)$-triangle rule subject to the grading.

To prove the last fact observe that the triangles $i,j,k\in\overline{\Ga(\g_0)}$ are of two kinds:
\begin{itemize}
\item Two of three quasiroots (see Definition \ref{110.60}), say $i,j$, are such that  $\g_i,\g_j \subset\g_0^+$, for the third one, $k$, we have $\g_k \subset\g_0^-$. Then $\ka_i=\ka_j=(t_1-t_2)/2=-\ka_k$, the label is $+$.
\item Two of three quasiroots, say $i,j$, are such that  $\g_i,\g_j \subset\g_0^-$, for the third one, $k$, we have $\g_k \subset\g_0^+$. Then $\ka_i=\ka_j=-(t_1-t_2)/2=-\ka_k$, the label is $-$.
 \end{itemize}
The last statement of the theorem follows from Theorem \ref{130.500}. \qed

\abz\label{130.65}
\begin{rema}\rm It is easy to see that the eigenspaces of the operator $W$ built in Theorem \ref{130.60} are the subalgebras in $(\g,[,])$. Hence the operator $W$ is a Nijenhuis operator, i.e. its torsion $T_W$ vanishes (cf. Example \ref{10.35}).
\end{rema}

Another series of examples (generalizing Example \ref{10.40}, see also \cite[Example 3]{golsok}) will be related to inner automorphisms of finite order and to more general admissible subalgebras.

Recall \cite[Ch. X]{helgason} that any such automorphism  of a simple Lie algebra $\g$ is conjugated to a model one which can be described as follows. Choose a Cartan subalgebra $\h \subset\g$ and let $R$ stand for the corresponding root system. Pick up a basis $B=\{\al_1,\ldots,\al_n\}$ of $R$ and let $a_0=1,a_1,\ldots,a_n$ be the labels of the Dynkin diagram $\widetilde{\Pi}$ of the extended system $\{\al_0,\ldots,\al_n\}$ (here $\al_0$ denotes the lowest root and $-\al_0=\sum_{l=1}^na_l\al_l$). Let $s_0,\ldots,s_n$ be nonnegative integers without nontrivial common factor. Then there exists a system of canonical generators $X_0,\ldots,X_n, X_k\in\g_{\al_k}$, of the Lie algebra $\g$ such that the formula $\sigma(X_k)=e^{s_k(2i\pi/m)}X_k$, where $m=\sum_{k=0}^na_ks_k$, defines uniquely an inner automorphism of $m$-th order, called a {\em model automorphism of type $(s_0,\ldots,s_n;1)$}. Any  conjugated automorphism is said to be an {\em  automorphism of type $(s_0,\ldots,s_n;1)$}.

The fixed point subalgebra $\g_0$ of the automorphism $\sigma$ contains $\h$ and is a direct sum of its $(n-p)$-dimensional centre and the semisimple Lie algebra $[\g_0,\g_0]$ whose Dynkin diagram $\pi$ is the subdiagram of $\widetilde{\Pi}$ consisting of the vertices $i_1,\ldots,i_p$ for which $s_{i_1}= \cdots=s_{i_p}=0$.
In particular, if $s_0\not=0$, the subalgebra $\g_0$ is a Levi subalgebra, and if $s_0=0$, the subalgebra $\g_0$ is a regular subalgebra of Class 2) (see Appendix \ref{ApToral}).

\abz\label{remmm}
\begin{rema}\rm
Any automorphism mentioned induces a $\Z_m  $-grading $\g=\g_{\bar 0}\oplus \cdots\oplus\g_{\overline{m-1}}$, where $\g_{\bar{j}}$ is the eigenspace of $\sigma$ corresponding to the eigenvalue $e^{(2i\pi/m)j(\mod m)}$ (in particular $\g_{\bar{0}}=\g_0$). This grading is toral with respect to $\h$. On the other hand, this grading is a coarsening of the toral irreducible $\Ga(\g_0)$-grading (\ref{grading1}) and there exists an additive epimorphism $\psi:\overline{\Ga(\g_0)}\to \overline{\Z_m}  $ such that $\g_{\bar j}=\bigoplus_{i\in\psi^{-1}(\bar j)}\g_i,j=1,\ldots,m-1$ (recall that we denote by $\overline{\Ga}$ the set of quasiroots of a $\Ga$-grading, see Definition \ref{110.60}).
\end{rema}

\abz\label{130.70}
\begin{theo}  \begin{enumerate}\item Let $\g_0$ be the fixed point subalgebra of an inner automorphism of order $m>2$ of type $(s_0,\ldots,s_n;1)$ (in particular $\g_0$ is admissible) and let $(\g_0^1,\g_0^2)$  be an admissible pair of subalgebras such that $\g_0^1\oplus\g_0^2=\g_0$. Let $\g=\g_{\bar 0}\oplus \cdots\oplus\g_{\overline{m-1}}$ be the corresponding grading.
Given $t_1,t_2\in\C, t_1\not=t_2$,  define an operator $W\in\End(\g)$  by
\begin{enumerate}\item
$W|_{\g_0^j}=t_i\Id_{\g_0^i},i=1,2$;
\item $W|_{\g_{\bar j}}=(((m-j)t_1+jt_2)/m)\Id_{\g_{\bar j}}$ if $j>0$.
    \end{enumerate}
Then the operator $W$ satisfies assumptions of Item 1 of Theorem \ref{130.50}; consequently, the triple $(\g,[,],[,]_W)$ is a semisimple bi-Lie structure of Class II satisfying conditions (a'), (b') of this theorem.
\item If $\g_0$ is a Levi subalgebra (i.e. $s_0\not=0$), the bi-Lie structure $(\g,[,],[,]_W)$ is isomorphic to  that constructed by data $(t_1,t_2,\g_0^1,\g_0^2,\h,B)$ in Theorem \ref{130.60}.
\item Two bi-Lie structures constructed by data $(t_1,t_2,\h,B,(s_0,\ldots,s_n),\g_0^1,\g_0^2)$ and  $(t_1',t_2', \h',B',\linebreak(s_0',\ldots,s_n'),(\g_0^1)',(\g_0^2)')$ are isomorphic if and only if either there exists $\phi\in\Aut(\g)$ such that $\phi(\g_0^j)=(\g_0^j)',j=1,2$, and $\phi(\p^+(\g_0^1\oplus\g_0^2,B))=\p^+((\g_0^1)'\oplus(\g_0^2)',B')$ or there exists $\phi\in\Aut(\g)$ such that $\phi(\g_0^1)=(\g_0^2)',\phi(\g_0^2)=(\g_0^1)'$, and $\phi(\p^+(\g_0^1\oplus\g_0^2,B))=\p^-((\g_0^1)'\oplus(\g_0^2)',B')$.
 \end{enumerate}
\end{theo}

\abz\label{130.71}
\begin{rema}\rm
The formulae (a), (b) appeared in \cite[Example 3]{golsok}.
\end{rema}

\noindent{\em Ad 1.}  Clearly it is enough consider a model automorphism of type $(s_0,\ldots,s_n;1)$. Since the $\Z_m$-grading is a coarsening of the toral irreducible $\Ga(\g_0)$-grading, the operator $W$ is scalar on the components of the latter grading. By Item 4 of Theorem \ref{100.20} we have $W_{\mathrm{s}}|_{\g_{\bar j}}=((t_1+t_2)/2)\Id_{\g_{\bar j}}$ for any $j>0$ and Conditions 1 (a), (b)  of Theorem \ref{130.50} are satisfied.

To check Condition 1 (c) due to Remark \ref{remmm} it is enough to check the $((t_1-t_2)/2)$-triangle rule subject to the $\Z_m  $-grading.
Observe that  $W_{\mathrm{a}}|_{\g_{\bar j}}=\ka_j\Id_{\g_{\bar j}}$, where $\ka_j=((-j+m/2)(t_1-t_2)/m)$,  and that triangles $\bar i,\bar j,\bar k\in\overline{\Z_m  }$, are of two kinds:
 \begin{itemize}
\item $i+j+k=m$ (then
$\ka_i+\ka_j+\ka_k=(t_1-t_2)/2$, the label is $+$);
\item $i+j+k=2m$  (then $\ka_i+\ka_j+\ka_k=-(t_1-t_2)/2$, the label is $-$).
\end{itemize}

\medskip

\noindent{\em Ad 2.}  This item follows from Theorem \ref{130.500} and from the fact that the corresponding operators have identical systems of labels with respect to the toral irreducible $\Ga(\g_0)$-grading.

To prove this fact notice that, if $\bar i,\bar j,\bar k$ is a triangle, then there exist roots $\al,\be,\ga$ such that $\g_{\al}\in\g_{\bar i},\g_{\be}\in\g_{\bar j},\g_{\ga}\in\g_{\bar k}$ and $\al+\be+\ga=0$. We will show that in the case $i+j+k=m$ two of these roots are positive (and the third negative) and vice versa in the case $i+j+k=2m$  (cf. the proof of Theorem \ref{130.60}).

Note that any root $\al$ can be uniquely represented as $\al=\sum_{l=0}^nk^\al_l\al_l$, where $0\le k^\al_l\le a_l$ ($k^\al_0=0$ if $\al>0$ and $k^\al_0=1$ if $\al<0$). Moreover, $\g_{\bar p}$ consists of those roots $\al$ for which $\sum_{l=0}^nk^\al_l s_l=p$ \cite[\S 3.7]{vinbergOnishIII}.

Let $i+j+k=2m$ and assume that two of the three roots, say $\al,\be$, are positive.  Then $\al=\sum_{l=1}^nk^\al_l\al_l,\be=\sum_{l=1}^nk^\be_l\al_l$ and $\al+\be=\sum_{l=1}^n(k^\al_l+k^\be_l)\al_l$, in particular $k^\al_l+k^\be_l\le a_l$ for any $l$. On the other hand, $i+j=\sum_{l=1}^n(k^\al_l+k^\be_l)s_l$ has to be greater than $m=s_0+\sum_{l=1}^na_ls_l$, whence $k^\al_l+k^\be_l>a_l$ for some $l\in\{1,\ldots,n\}$, a contradiction.

Now let $i+j+k=m$ and assume that two of the three roots, say $\al,\be$, are negative.  Then $\al=\al_0+\sum_{l=1}^nk^\al_l\al_l,\be=\al_0+\sum_{l=1}^nk^\be_l\al_l$ and $\al+\be=\al_0+\sum_{l=1}^n(k^\al_l+k^\be_l-a_l)\al_l$, in particular $k^\al_l+k^\be_l\ge a_l$ for any $l$. On the other hand, $i+j=2s_0+\sum_{l=1}^n(k^\al_l+k^\be_l)s_l$ has to be less than $m=s_0+\sum_{l=1}^na_ls_l$, whence $k^\al_l+k^\be_l<a_l$ for some $l\in\{1,\ldots,n\}$, a contradiction.

\medskip

\noindent{\em Ad 3.}  The proof of the previous item shows that the system of labels of the antisymmetric part of the operator $W$ depends only on the subalgebra $\g_0$ and the decomposition of $\g_0^\bot$ to subspaces corresponding to positive and negative roots. Now we can use Theorem \ref{130.500}.
\qed

\medskip

The proofs of Theorems \ref{130.60}, \ref{130.70} show that the   corresponding systems of labels have the following simple property (up to a sign):
\begin{itemize}
\item
any triangle with nonzero label and two positive (negative) roots has label $+$ ($-$).
\end{itemize}
Using this observation we can build new examples of a bi-Lie structures of Class II, in particular, having the most general reductive subalgebras of maximal rank as their basic subalgebras.

\abz\label{130.80}
\begin{exa}\rm Let $\g=\e_7$. We will use notations from \cite[Table VI]{bourb4-6}.  Let $R_0=\Span_\Z\{\al_0,\al_1,\lbr\al_2,\al_4,\al_6,\al_7\}\cap R=\Span_\Z \{\al_1,\al_2,3\al_3,\al_4,3\al_5,\al_6,\al_7\}\cap R$, here $\al_0=-(2\al_1+2\al_2+3\al_3+4\al_4+3\al_5+2\al_6+\al_7)$ is the lowest root. Then the corresponding regular reductive Lie subalgebra $\g_0$ is of Class 3) (see Appendix \ref{ApToral}) since $\Ga(\g_0)=Q/Q_0\cong\Z_3 \times\Z_3 $. The subalgebra $\g_0$ is the intersection of the fixed point subalgebras of two commuting inner automorphisms of order 3 and is not a fixed point subalgebra of any single inner automorphism of finite order (i.e. $\g_0$ is not of Class 2)).

To build the antisymmetric part of a WNO we will use the induction process from the proof of  Lemma \ref{130.40}. Starting this process from any vector $(c_1,\ldots,c_7)$ we will obtain a diagonal  antisymmetric operator $S\in\End(\h^\bot), S|_{\g_\al}=\ka_\al\Id_{\g_\al}$. We have to take care of the condition $\ka_\al=0,\al\in R_0$ (see Definitions \ref{130.30}, \ref{110.60}),  which gives a system of linear  equations on the variables $(c_1,\ldots,c_7)$. Five equations from this system are obvious: $c_1=0,c_2=0,c_4=0,c_6=0,c_7=0$. Thus in fact we can start the induction process from the vector $(0,0,c_3,0,c_5,0,0)$.
In the result the nonzero entries will correspond only to the roots containing the combinations of the roots $\al_3,\al_5$:
$$
\begin{tabular}{|l|l||l|l||l|l||l|l|}\hline
$\al_3$ & $c_3$ & $\al_3+\al_5$ & $c_3+c_5+a$ &  $2\al_3+\al_5$ & $2c_3+c_5+2a$ & $2\al_3+3\al_5$ & $2c_3+3c_5+4a$\\ \hline
$\al_5$ & $c_5$ &  $\al_3+2\al_5$ & $c_2+2c_5+2a$ & $2\al_3+2\al_5$ & $2c_3+2c_5+3a$ & &\\ \hline
\end{tabular}
$$
Note that the combination $3\al_3+2\al_5$ do not entry in any root and the combination $3\al_3+3\al_5$ gives the only relation on the constants $(c_3,c_5)$: $3c_3+3c_5+4a=0$ (the last equation from the above mentioned system). Putting $x:=c_3$ we get $c_5=-x-4a/3$ and the following table:
\begin{center}
\begin{tabular}{|l|l||l|l||l|l||l|l|}\hline
$\al_3$ & $x$ & $\al_3+\al_5$ & $-a/3$ &  $2\al_3+\al_5$ & $x+2a/3$ & $2\al_3+3\al_5$ & $-x$\\ \hline
$\al_5$ & $-x-4a/3$ &  $\al_3+2\al_5$ & $-x-2a/3$ & $2\al_3+2\al_5$ & $a/3$ & &\\ \hline
\end{tabular}
\end{center}
To complete the construction we need to: 1) extend this  table by the skew symmetry to negative roots; 2) split $\g_0$ to an admissible pair $(\g_0^1,\g_0^2)$ (remark that the subalgebra $[\g_0,\g_0]$ is a sum of three simple subalgebras, so such splitting can be done in several ways); 3) choose $t_1,t_2\in\C,t_1\not=t_2$, and put $W_{\mathrm{s}}|_{\g^i_0}=t_i\Id_{\g^i_0}$ for $i=1,2$ and $W_{\mathrm{s}}|_{\g_\al}=((t_1+t_2)/2)\Id_{\g_\al}$ for $\al\not\in R_0$; 4) recall that $a=(t_1-t_2)/2$. We get a one-parameter family of bi-Lie structures of Class II with the set of times $\{t_1,t_2\}$.
\end{exa}

This example can be generalized to the case of arbitrary reductive subalgebra $\g_0$. We can now formulate
\abz\label{130.90}
\begin{conj} For any bi-Lie structure $(\g,[,],[,]')$ of Class II there exists a Cartan subalgebra $\h$ and a basis of the root system $R=R(\g,\h)$ such that the system of labels of the antisymmetric part of corresponding principal WNO restricted to $\h^\bot$ has the property mentioned before Example \ref{130.80}. In particular,
\begin{enumerate}\item
any bi-Lie structure of Class II with the basic subalgebra $\g_0$ being a Levi subalgebra (i.e. a regular reductive subalgebra of Class 1), see Appendix \ref{ApToral}) is isomorphic to one of that built in Theorem \ref{130.60};
\item any bi-Lie structure of Class II with the basic subalgebra $\g_0$ being a regular reductive subalgebra of Class 2) is isomorphic to one of that built in Theorem \ref{130.70}.\end{enumerate}
\end{conj}

We conclude this section by proving this conjecture in a particular case $(\g,[,])=\a_n$ thus obtaining a complete classification of  bi-Lie structures of Class II in this case.

\abz\label{130.100}
\begin{theo} Let $(\g,[,])=\a_n$. Then any semisimple regular bi-Lie structure $(\g,[,],[,]')$ of Class II   is isomorphic to one of that built in Theorem \ref{130.60}.
\end{theo}

\noindent Let $\{t_1,t_2\}$ be the set of times of the bi-Lie structure $(\g,[,],[,]')$. Let $\{\al_1,\ldots,\al_n\}$ be a basis of roots of the corresponding root system $R(\g,\h)$ and let $R_0 \subset R$ be the root set corresponding to the basic subalgebra $\g_0$. By Corollary \ref{110.501} the corresponding principal WNO obeys the $a$-triangle rule subject to $R_0$, where $a=(t_1-t_2)/2$. We will first prove that there exists a WNO $W$ of this bi-Lie structure such that all the eigenvalues of $W_{\mathrm{a}}|_{\h^\bot}$ are equal to zero or $\pm a$.

 Consider the triangle $\{\be_{i-1},\al_i,-\be_i\},i=2,\ldots,n$, where we write $\be_i:=\al_1+\cdots+\al_i$, and put $d_i=0,\pm 1$ for the corresponding label. Consider the vector $(d_1a,d_2a,\ldots,d_na)$, where by definition $d_1=0$ if $\al_1\in R_0$ and $d_1=1$ otherwise, and start the induction process of the proof of Lemma \ref{130.40} from this vector. We claim that the eigenvalues $\ka_\al$ of the resulting antisymmetric operator will satisfy the required condition.

Indeed, let $\al:=\al_k+\al_{k+1}+\cdots+\al_{k+m}=\be_{k+m}-\be_{k-1}$ for some $1\le k\le n-1,m\ge 0$. Then it belongs to the triangle $\{\al,\be_{k-1},-\be_{k+m}\}$ with the label $d=0,\pm 1$ and
\begin{equation}\equ\label{label}
\ka_\al=-\ka_{\be_{k-1}}+\ka_{\be_{k+m}}+d a.
\end{equation}
On the other hand, the induction process gives $\ka_{\be_1}=\ka_{\al_1}=d_1a,\ka_{\be_2}=\ka_{\be_1}+\ka_{\al_2}-d_2a=d_1a+d_2a-d_2a=d_1a,\ldots,
\ka_{\be_{k-1}}=
\ka_{\be_{k-2}}+\ka_{\al_{k-1}}-d_{k-1}a=d_1a+d_{k-1}a-d_{k-1}a=d_1a,\ldots,\ka_{\be_{k+m}}=d_1a$. Hence by (\ref{label}) we get $\ka_\al=da$.

In the next step of the proof consider the operator just built and notice that the set $Y:=\{\al\in R\mid \ka_\al=a\}$ is closed: if $\al,\be\in Y$ and $\ga:=\al+\be\in R$, then $\ka_\ga=\ka_\al+\ka_\be-d'a=a$, here $d'$ is the corresponding label, which necessarily equals $+1$. Now  remark that $Y\cap(-Y)=\emptyset$ and use \cite[Prop. 22, Ch.VI]{bourb4-6} to deduce that there exists a basis of $R$ such that $Y$ is contained in the set of positive roots $R^+$ with respect to this basis.

Finally, observing that any triangle with a nonzero label and with two positive roots by construction has label $+1$ we prove Conjecture \ref{130.90}.

The rest of the proof follows from Theorems \ref{130.50}, \ref{130.500} and Remark \ref{ApRem1}.
 \qed

\section{Appendix I: Toral gradings and regular reductive Lie subalgebras}
\label{ApToral}

\abz\label{ApToral.05}
\begin{defi}\rm
Let $\g$ be a semisimple Lie algebra and let $\g=\h+\sum_{\al\in R}\g_\al$ be the root decomposition of $\g$ with respect to some Cartan subalgebra $\h\subset\g$. We say that a $\Ga$-grading $\g=\bigoplus_{i\in \Ga}\g_i$, where $\Ga$ is an abelian group, is   {\em toral} with respect to $\h$ if the subspaces $\g_i$ are spanned by the subspaces $\h,\g_\al$ and $\g_0\supset\h$.
\end{defi}

Let  $\g_0$ be a reductive subalgebra in $\g$ such that $\g_0\supset\h$. In particular the representation $x\mapsto \ad_\g x$ of the Lie algebra $\g_0$ in $\g_0^\bot$ (here $\g_0^\bot$ is the orthogonal complement to $\g_0$ in $\g$ with respect to the Killing form) is semisimple.

Let $R_0:=\{\al\in R\mid \g_\al \subset[\g_0,\g_0]\} \subset R\subset\h_\R^*$ be the root system of the subalgebra $[\g_0,\g_0]$ considered as the subsytem of $R$ and let $Q_0\subset Q\subset\h_\R^*$ be the corresponding root lattices of $[\g_0,\g_0]$ and $\g$. Given $\al\in Q$, put $\g_{\al+Q_0}:=\bigoplus_{\be\in(\al+Q_0)\cap R}\g_\be$.

\abz\label{ApToral.10}
\begin{theo}\cite[Section 6.2]{ostapenko}
The decomposition $\g=\bigoplus_{\al+Q_0\in Q/Q_0}\g_{\al+Q_0}$ is a toral $Q/Q_0$-grading. Moreover, the subspaces $\g_{\al+Q_0},\al\not\in Q_0$, are the irreducible components of the representation of $\g_0$ in $\g_0^\bot$.
\end{theo}

Let $\h'$ be another Cartan subalgebra such that $\h' \subset\g_0$ and let $Q',Q'_0$ stand for the corresponding root lattices. Then there exists $\phi\in\Aut(\g)$ such that $\phi(\h)=\h'$. It is easy to see that $(\phi|_\h)^T(Q')=Q$ and $(\phi|_\h)^T(Q'_0)=Q_0$, here $(\phi|_\h)^T$ is the operator transposed to $\phi|_\h$. In particular, the groups $\Ga(\h,\g_0):=Q/Q_0$ and $\Ga(\h',\g_0):=Q'/Q'_0$ are isomorphic.

\abz\label{ApToral.20}
\begin{defi}\rm Let $\g_0\subset \g$ be a reductive in $\g$ subalgebra  of maximal rank  (call such a subalgebra {\em regular}) and let $\Ga(\g_0)$ be an abstract group isomorphic to one of the groups $\Ga(\h,\g_0)$, where $\h$ is a Cartan subalgebra such that $\h \subset\g_0$. The grading constructed above will be called the {\em irreducible toral $\Ga(\g_0)$-grading} corresponding to the subalgebra $\g_0$.
\end{defi}

Regular subalgebras can be divided to the following three classes depending on the structure of the group $\Ga(\g_0)$ (cf. \cite{doninOst}): 1) $\Ga(\g_0)$ is free; 2) the torsion component $\Tor(\Ga(\g_0))$ of $\Ga(\g_0)$ is cyclic; 3) $\Tor(\Ga(\g_0))$ is not cyclic.

 In the first case the subalgebra $\g_0$ is the so-called Levi subalgebra, which can be also characterized by the following equivalent conditions  (\cite[Th. 1.3, Ch. 6]{vinbergOnishIII}):
\begin{enumerate}\item $\g_0$ is a centralizer $\z_\g(\c)$ of some subspace $\c \subset\h$, where $\h$ is some Cartan subalgebra.
\item $\g_0$ contains a Cartan subalgebra $\h$ and the corresponding system $R_0 \subset R$ of roots of $\g_0$ is generated over $\Z$ by some subset $B_0$ of some system $B(\g)$ of  simple roots of $\g$. \end{enumerate}

To describe the other two cases recall \cite{dynkin,vinbergOnishIII} that any regular subalgebra $\g_0 \subset\g$ can be obtained as the last element in a nested sequence of subalgebras $\g=\g^0\supset\g^1\supset \cdots\supset\g^k=\g_0$, where each $\g^j,j=1,\ldots,k-1$,    is obtained from  $\g^{j-1}$ by an elementary transformation and $\g_0$ is a Levi subalgebra in $\g^{k-1}$. We say that a semisimple subalgebra $\k \subset\g$ of a semisimple Lie algebra $\g=\bigoplus_{j\in J}\g_j$, where $\g_j$ are the simple components,  is obtained by an {\em elementary transformation} from $\g$ if $\k=\bigoplus_{j\not=j'}\g_j\oplus\tilde{\g}_{j'}$ for some $j'\in J$. Here $\tilde{\g}_{j'}$ is a subalgebra in $\g_{j'}$ such that
its system of simple roots $B(\g_{j'})$ is obtained  by eliminating of some root $\al_i,i>0$, from the system  $B(\g_{j'})\cup\{\al_0\}$ of simple roots of $\g_{j'}$ extended by the lowest root $\al_0$. In other words, $\g_{j'}$ is the fixed point subalgebra of an (inner) automorphism $\sigma\in\Aut(\g_{j'})$ of order $a_i$ of type $(0,\ldots,0,1,0,\ldots,0;1)$ (the first unit is on the $i$-th place and $i>0$) in terminology of \cite[Ch. X]{helgason} (see also the discussion after Remark \ref{130.65}), here $a_i$ is the label (i.e. the coefficient in the highest root) of the root $\al_i$.

The equivalent characterization of regular subalgebras of Class 2) (respectively 3)) is, \cite{doninOst}, that an elementary transformation in the sequence $\g\supset \cdots\supset\g_0$ is made only once (respectively more than once).

\abz\label{ApRem1}
\begin{rema}\rm The structure of the extended Dynkin diagram of the Lie algebra $\g=\a_n$ implies that any regular subalgebra in $\g$ is a Levi subalgebra.
\end{rema}

\section{Appendix II: bi-Lie structures of Class I on compact real forms of complex semisimple Lie algebras}
\label{ApReal}

\abz\label{ApReal.10}
\begin{theo} Let $(\g,[,])$ be a complex semisimple Lie algebra, $\h \subset\g$ a Cartan subalgebra, $\g=\h+\sum_{\al\in R}\g_\al$ the corresponding root grading. Let $\g=\bigoplus_{t_i,t_j\in X}\g_{t_it_j}$  a toral (with respect to $\h$) symmetric pairoid quasigrading   of class I with the base $X=\{t_1,\ldots,t_n\}$, where $t_i\in\R$, such that the subspaces $\g_{t_it_i}\cap\h$ are the complexifications of subspaces in $\h_\R$.

For any $\al\in R$ choose $E_\al\in\g_\al$ such that $B_\g(E_\al,E_{-\al})=1$ and put $H_\al:=[E_\al,E_{-\al}]$ as usual. Let $\uu=\sum_{\al\in R}\R(iH_\al)+\sum_{\al\in R}\R(E_\al-E_{-\al})+\sum_{\al\in R}\R(i(E_\al+E_{-\al})) \subset\g$ be the compact real form of $\g$ related with the root decomposition and the choice of $E_\al$ (see \cite[Th. 6.3, Ch. III]{helgason})

Then the WNO $W$ built   by the quasigrading  (see Theorem \ref{140.50}) has a correct restriction to $\uu$ which is a real WNO and induces a real bi-Lie structure on $\uu$.
\end{theo}

\noindent The fact that $W$ preserves $\uu$ follows from the assumptions and from the construction of $W$ (recall $W|_{\h^\bot}$ is symmetric: $W|_{\g_\al}=W|_{\g_{-\al}}$). The principal primitive $P$ of $W$ is symmetric too (see Definition \ref{prpr} or check this directly using the formulae from the proof of Theorem \ref{140.50}). Thus the main identity $T_N(,)=[,]_P$ implies $T_{N|_\uu}(,)=[,]_{P|_\uu}$. \qed

\abz\label{ApReal.20}
\begin{exa}\rm From Example \ref{140.120} we get the following example of WNO on $\su(n+1)$: $WX:=(1/2)(L_A+R_A)X-\Tr((1/2)(L_A+R_A)X)B$, where $A=\diag(t_1,t_2,\ldots,t_{n+1}),B=\diag(0,0,\ldots,0,1)$ and $t_i\in\R$. Explicitly, if $X=||x_{ij}||\in\su(n+1)$, then $WX=||y_{ij}||$, where $y_{ij}=x_{ij}(t_i+t_j)/2$ for $i\not=j$, $y_{ii}=x_{ii}t_i$ for $i=1,\ldots,n$, and $y_{n+1,n+1}=-\sum_{j=1}^{n}x_{jj}t_j)$.

Analogous example can be obtained from Example \ref{140.140} taking the parameter $a$ to be real.
\end{exa}

\abz\label{ApReal.30}
\begin{exa}\rm Here we will only mention the existence of  a nonstandard bi-Lie structure on $\so(6,\R)$ coming from the isomorphism $\so(6,\C)\cong\sl(4,\C)$ and Examples \ref{140.120}, \ref{140.140}.
\end{exa}

\section{Acknowledgements}

Research partially supported by Polish Ministry of Science and Higher Education, grant N201 607540. The author is indebted to Taras Skrypnyk and Ilya Zakharevich for helpful discussions and to Kirill Mackenzie for suggesting the name "pairoid".

\providecommand{\bysame}{\leavevmode\hbox to3em{\hrulefill}\thinspace}
\providecommand{\MR}{\relax\ifhmode\unskip\space\fi MR }
\providecommand{\MRhref}[2]{%
  \href{http://www.ams.org/mathscinet-getitem?mr=#1}{#2}
}
\providecommand{\href}[2]{#2}

\end{document}